\documentclass{article}
\pdfoutput=1
\usepackage[utf8]{inputenc}
\usepackage[serif]{lindrew}
\usepackage{cancel}
\usepackage{pgf, tikz}
\usepackage{pdfpages}
\usepackage{mathdots}
\usepackage{subcaption}
\usetikzlibrary{arrows, automata}
\usepackage{upgreek}
\usepackage{amsmath, amssymb}
\usepackage{graphicx}
\usepackage{wrapfig}

\usepackage[frozencache,cachedir=_minted]{minted}

\newcommand{\todo}[1]{}

\title{Matrix Calculus \\ (for Machine Learning and Beyond)}
\author{\textbf{Lecturers: Alan Edelman and Steven~G.~Johnson} \\ Notes by Paige Bright, Alan Edelman, and Steven~G.~Johnson}

\usepackage[utf8]{inputenc}
\DeclareUnicodeCharacter{03B5}{$\upvarepsilon$}
\DeclareUnicodeCharacter{03B4}{$\updelta$}
\DeclareUnicodeCharacter{03C9}{$\upomega$}
\DeclareUnicodeCharacter{2032}{$^\prime$}
\DeclareUnicodeCharacter{2248}{$\approx$}
\DeclareUnicodeCharacter{221A}{$\surd$}
\DeclareUnicodeCharacter{03F5}{$\upepsilon$}
\DeclareUnicodeCharacter{2207}{$\nabla$}
\DeclareUnicodeCharacter{03B1}{$\upalpha$}

\renewcommand{\d}{d}
\newcommand{\diag}{\operatorname{diag}}
\newcommand{\cofactor}{\mathrm{cofactor}}
\newcommand{\adj}{\operatorname{adj}}
\newcommand{\vecm}{\operatorname{vec}}
\newcommand{\dotstar}{\mathbin{.*}}
\newcommand{\evalat}[2]{\left. #1 \right|_{#2}}
\newcommand{\red}[1]{\textcolor{red}{#1}}

\usepackage{kbordermatrix}
\renewcommand{\kbldelim}{(} 
\renewcommand{\kbrdelim}{)}

\date{Based on MIT course 18.S096 (now 18.063) in IAP 2023}

\begin{document}

\maketitle

\tableofcontents

\pagebreak

\addcontentsline{toc}{section}{Introduction}
\section*{Introduction}
These notes are based on the class as it was run for the second time in January 2023, taught by Professors Alan Edelman and Steven~G.~Johnson at MIT. The previous version of this course, run in January 2022, can be found \href{https://ocw.mit.edu/courses/18-s096-matrix-calculus-for-machine-learning-and-beyond-january-iap-2022/}{on OCW here}. 

Both Professors Edelman and Johnson use he/him pronouns and are in the Department of Mathematics at MIT; Prof.~Edelman is also a Professor in the MIT Computer Science and Artificial Intelligence Laboratory (CSAIL) running the Julia lab, while Prof.~Johnson is also a Professor in the Department of Physics. 

Here is a description of the course.:
\begin{quote}
    We all know that typical  calculus course sequences 
    begin with  univariate and vector calculus, respectively. Modern applications such as machine learning and large-scale optimization require the next big step, ``matrix calculus'' and calculus on arbitrary vector spaces.


This class covers a coherent approach to matrix calculus showing techniques that allow you to think of a matrix holistically (not just as an array of scalars), generalize and compute derivatives of important matrix factorizations and many other complicated-looking operations, and understand how differentiation formulas must be re-imagined in large-scale computing. We will discuss ``reverse'' (``adjoint'', ``backpropagation'') differentiation and how modern automatic differentiation is more computer science than calculus (it is neither symbolic formulas nor finite differences).
\end{quote}
The class involved numerous example numerical computations using the Julia language, which you can install on your own computer following \href{https://github.com/mitmath/julia-mit#installing-julia-and-ijulia-on-your-own-computer}{these instructions}. The material for this class is also located on GitHub at  \url{https://github.com/mitmath/matrixcalc}.

\pagebreak

\section{Overview and Motivation}
Firstly, where does matrix calculus fit into the MIT course catalog? Well, there are 18.01 (Single-Variable Calculus) and 18.02 (Vector Calculus) that students are required to take at MIT. But it seems as though this sequence of material is being cut off arbitrarily:
\[
\text{Scalar }\to \text{ Vector }\to \text{Matrices}  \to \text{ Higher-Order Arrays?}
\]
After all, this is how the sequence is portrayed in many computer programming languages, including Julia!  Why should calculus stop with vectors?

In the last decade, linear algebra has taken on larger and larger importance in numerous areas, such as machine learning, statistics, engineering, etc. In this sense, linear algebra has gradually taken over a much larger part of today's tools for lots of areas of study---now everybody needs linear algebra. So it makes sense that we would \textit{want} to do calculus on these higher-order arrays, and it won't be a simple/obvious generalization (for instance, $\frac{\d}{\d A} A^2 \neq 2A$ for non-scalar matrices $A$).

More generally, the subjects of \emph{differentiation} and \emph{sensitivity analysis} are much deeper than one might suspect from the simple rules learned in first- or second-semester calculus.  Differentiating functions whose inputs and/or outputs are in more complicated vector spaces (e.g. matrices, functions, or more) is one part of this subject.  Another topic is the \emph{efficient} evaluation of derivatives of functions involving very complicated calculations, from neural networks to huge engineering simulations---this leads to the topic of ``adjoint'' or ``reverse-mode'' differentiation, also known as ``backpropagation.''  \emph{Automatic differentiation (AD)} of computer programs by compilers is another surprising topic, in which the computer does something very different from the typical human process of first writing out an explicit symbolic formula and then passing the chain rule through it.   These are only a few examples: the key point is that differentiation is more complicated than you may realize, and that these complexities are increasingly relevant for a wide variety of applications.
 
Let's quickly talk about some of these applications.

\subsection{Applications}

\subsubsection*{Applications: Machine learning}

Machine learning has numerous buzzwords associated with it, including but not limited to: parameter optimization, stochastic gradient descent, automatic differentiation, and backpropagation. In this whole collage you can see a fraction of how matrix calculus applies to machine learning. It is recommended that you look into some of these topics yourself if you are interested.


\subsubsection*{Applications: Physical modeling}

Large physical simulations, such as engineering-design problems, are increasingly characterized by \emph{huge} numbers of parameters, and the \emph{derivatives} of simulation outputs with respect to these parameters is crucial in order to evaluate sensitivity to uncertainties as well as to apply large-scale optimization.

For example, the shape of an airplane wing might be characterized by thousands of parameters, and if you can compute the derivative of the drag force (from a large fluid-flow simulation) with respect to these parameters then you could optimize the wing shape to minimize the drag for a given lift or other constraints.

An extreme version of such parameterization is known as ``topology optimization,'' in which the material at ``every point'' in space is potentially a degree of freedom, and optimizing over these parameters can discover not only a optimal shape but an optimal \emph{topology} (how materials are connected in space, e.g.~how many holes are present).   For example, topology optimization has been applied in mechanical engineering to design the cross sections of airplane wings, artificial hips, and more into a complicated lattice of metal struts (e.g.~minimizing weight for a given strength).

Besides engineering design, complicated differentiation problems arise in \emph{fitting} unknown parameters of a model to experimental data, and also in \emph{evaluating uncertainties} in the outputs of models with imprecise parameters/inputs.   This is closely related to regression problems in statistics, as discussed below, except that here the model might be a giant set of differential equations with some unknown parameters.

\subsubsection*{Applications: Data science and multivariable statistics}

In multivariate statistics, models are often framed in terms of matrix inputs and outputs (or even more complicated objects such as tensors).  For example, a ``simple'' linear multivariate matrix model might be $Y(X) = XB + U$, where $B$ is an unknown matrix of coefficients (to be determined by some form of fit/regression) and $U$ is unknown matrix of random noise (that prevents the model from exactly fitting the data).  Regression then involves minimizing some function of the error $U(B) = Y - XB$ between the model $XB$ and data $Y$; for example, a matrix norm $\Vert U\Vert_F^2 = \tr U^T U$, a determinant $\det U^T U$, or more complicated functions.  Estimating the best-fit coefficients $B$, analyzing uncertainties, and many other statistical analyses require differentiating such functions with respect to $B$ or other parameters.  A recent review article on this topic is Liu et~al.~(2022): ``Matrix differential calculus with applications in the multivariate linear model and its diagnostics'' (\url{https://doi.org/10.1016/j.sctalk.2023.100274}).


\subsubsection*{Applications: Automatic differentiation}

Typical differential calculus classes are based on symbolic calculus, with students essentially learning to do what Mathematica or Wolfram Alpha can do. Even if you are using a computer to take derivatives symbolically, to use this effectively you need to understand what is going on beneath the hood. But while, similarly, some numerics may show up for a small portion of this class (such as to approximate a derivative using the difference quotient), \textit{today's} automatic differentiation is neither of those two things. It is more in the field of the computer science topic of compiler technology than mathematics. However, the underlying mathematics of automatic differentiation is interesting, and we will learn about this in this class!

Even \emph{approximate} computer differentiation is more complicated than you might expect.  For single-variable functions $f(x)$, derivatives are defined as the limit of a difference $[f(x+\delta x)-f(x)]/\delta x$ as $\delta x\to 0$.  A crude ``finite-difference'' approximation is simply to approximate $f'(x)$ by this formula for a small $\delta x$, but this turns out to raise many interesting issues involving balancing truncation and roundoff errors, higher-order approximations, and numerical extrapolation.

\subsection{First Derivatives}

The derivative of a function of one variable is itself a function of one variable-- it simply is (roughly) defined as the linearization of a function. I.e., it is of the form $(f(x)-f(x_0)) \approx f'(x_0) (x-x_0)$. In this sense, ``everything is easy'' with scalar functions of scalars (by which we mean, functions that take in one number and spit out one number).

There are occasionally other notations used for this linearization:
\begin{itemize}
    \item $\delta y \approx f'(x) \delta x$,
    \item $\d y = f'(x) \d x$,
    \item $(y-y_0) \approx f'(x_0)(x-x_0)$,
    \item and $\d f = f'(x)  \d x$.
\end{itemize}
This last one will be the preferred of the above for this class. One can think of $\d x$ and $\d y$ as ``really small numbers.'' In mathematics, they are called \href{https://en.wikipedia.org/wiki/Infinitesimal}{infinitesimals}, defined rigorously via taking limits. Note that here we do not want to divide by $\d x$. While this is completely fine to do with scalars, once we get to vectors and matrices you can't always divide!

The numerics of such derivatives are simple enough to play around with. For instance, consider the function $f(x) = x^2$ and the point $(x_0,f(x_0)) = (3,9)$. Then, we have the following numerical values near $(3,9)$:
\begin{align*}
    f(3.\mathbf{0001}) &= 9.\mathbf{0006}0001 \\
    f(3.\mathbf{00001}) &= 9.\mathbf{00006}00001 \\
    f(3.\mathbf{000001}) &= 9.\mathbf{000006}000001 \\
    f(3.\mathbf{0000001}) &= 9.\mathbf{0000006}0000001.
\end{align*}
Here, the bolded digits on the left are $\Delta x$ and the bolded digits on the right are $\Delta y$. Notice that $\Delta y = 6\Delta x$. Hence, we have that 
\[
f(3+\Delta x) = 9 + \Delta y = 9+6\Delta x \implies f(3+\Delta x) - f(3) = 6\Delta x \approx f'(3)\Delta x.
\]
Therefore, we have that the linearization of $x^2$ at $x=3$ is the function $f(x)-f(3) \approx 6(x-3)$.

\vspace{.15cm}
\hrule
\vspace{.15cm}

We now leave the world of scalar calculus and enter the world of vector/matrix calculus! Professor Edelman invites us to think about matrices \textit{holistically}---not just as a table of numbers.

The notion of linearizing your function will conceptually carry over as we define the derivative of  functions which take in/spit out more 
than one number.
Of course, this means that the derivative will have a different ``shape'' than a single number. Here is a table on the \textit{shape} of the first derivative. The inputs of the function are given on the left hand side of the table, and the outputs of the function are given across the top. 

\begin{table}[h]
\centering
\begin{tabular}{l|l|l|l|}
\cline{2-4}
                                            
         input $\downarrow$  and output $\rightarrow$   & scalar & vector & matrix \\ \hline
\multicolumn{1}{|l|}{scalar} & scalar & vector (for instance, velocity) & matrix \\ \hline
\multicolumn{1}{|l|}{vector}  & gradient = (column) vector & matrix (called the Jacobian matrix) & higher order array  
\\ \hline
\multicolumn{1}{|l|}{matrix}   & matrix & higher order array & higher order array  \\ \hline
\end{tabular}
\end{table}

You will ultimately learn how to do any of these in great detail eventually in this class! The purpose of this table is to plant the notion of differentials as linearization. Let's look at an example.

\begin{example}
Let $f(x) = x^T x$, where $x$ is a $2\times 1$ matrix and the output is thus a $1\times 1$ matrix.
Confirm that $2x_0^T dx$
is indeed the differential of $f$ at $x_0 = \begin{pmatrix}  3 & 4\end{pmatrix}^T$.
\end{example}

Firstly, let's compute $f(x_0)$:
\[
f(x_0) = x_0^T x_0 = 3^2 + 4^2 = 25.
\]
Then, suppose $\d x = [.001,.002]$. Then, we would have that 
\[
f(x+\d x) = (3.001)^2 + (4.002)^2 = 25.\mathbf{022}005.
\]
Then, notice that $2x_0^T \,\d x = 2\begin{pmatrix}
    3  & 4 
\end{pmatrix}^T \d x = .022$. Hence, we have that 
\[
f(x_0 + \d x) - f(x_0) \approx 2x_0^{T} \d x = .022.
\]
As we will see right now, the $2x_0^T\d x$ didn't come from nowhere!

\subsection{Intro: Matrix and Vector Product Rule}

For matrices, we in fact still have a product rule!  We will discuss this in much more detail in later chapters, but let's begin here with a small taste.

\begin{theorem}[Differential Product Rule]
Let $A,B$ be two matrices. Then, we have the differential product rule for $AB$:
\[
\d(AB) = (\d A)B + A(\d B).
\]
By the differential of the matrix $A$, we think of it as a small (unconstrained) change in the matrix $A.$ Later, constraints may be places on the allowed perturbations.
\end{theorem}
Notice however, that (by our table) the derivative of a matrix is a matrix! So generally speaking, the products will not commute.

If $x$ is a vector, then by the differential product rule we have
\[
\d(x^Tx) = (\d x^T)x + x^T (\d x).
\]
However, notice that this is a dot product, and dot products commute (since $\sum a_i \cdot b_i = \sum b_i \cdot a_i$), we have that 
\[
\d(x^T x) = (2x)^T \d x.
\]

\begin{remark}
The way the product rule works for vectors as matrices is that transposes ``go for the ride.'' See the next example below.
\end{remark}

\begin{example}
    By the product rule. we have 
    \begin{enumerate}
        \item $\d(u^Tv) = (\d u)^T v + u^T (\d v) = v^T \d u + u^T \d v$ since dot products commute.
        \item $\d(uv^T) = (\d u )v^T + u( \d v)^T.$
    \end{enumerate}
\end{example} %

\begin{remark}
    The way to prove these sorts of statements can be seen in Section \ref{sec1S}.
\end{remark}


\pagebreak

\section{Derivatives as Linear Operators} \label{sec1S}
We are now going to revisit the notion of a derivative in a way that we can generalize to higher-order arrays and other vector spaces. We will get into more detail on differentiation as a linear operator, and in particular, will dive deeper into some of the facts we have stated thus far.

\subsection{Revisiting single-variable calculus}

\begin{figure}
    \centering
    \includegraphics[width=0.8\textwidth]{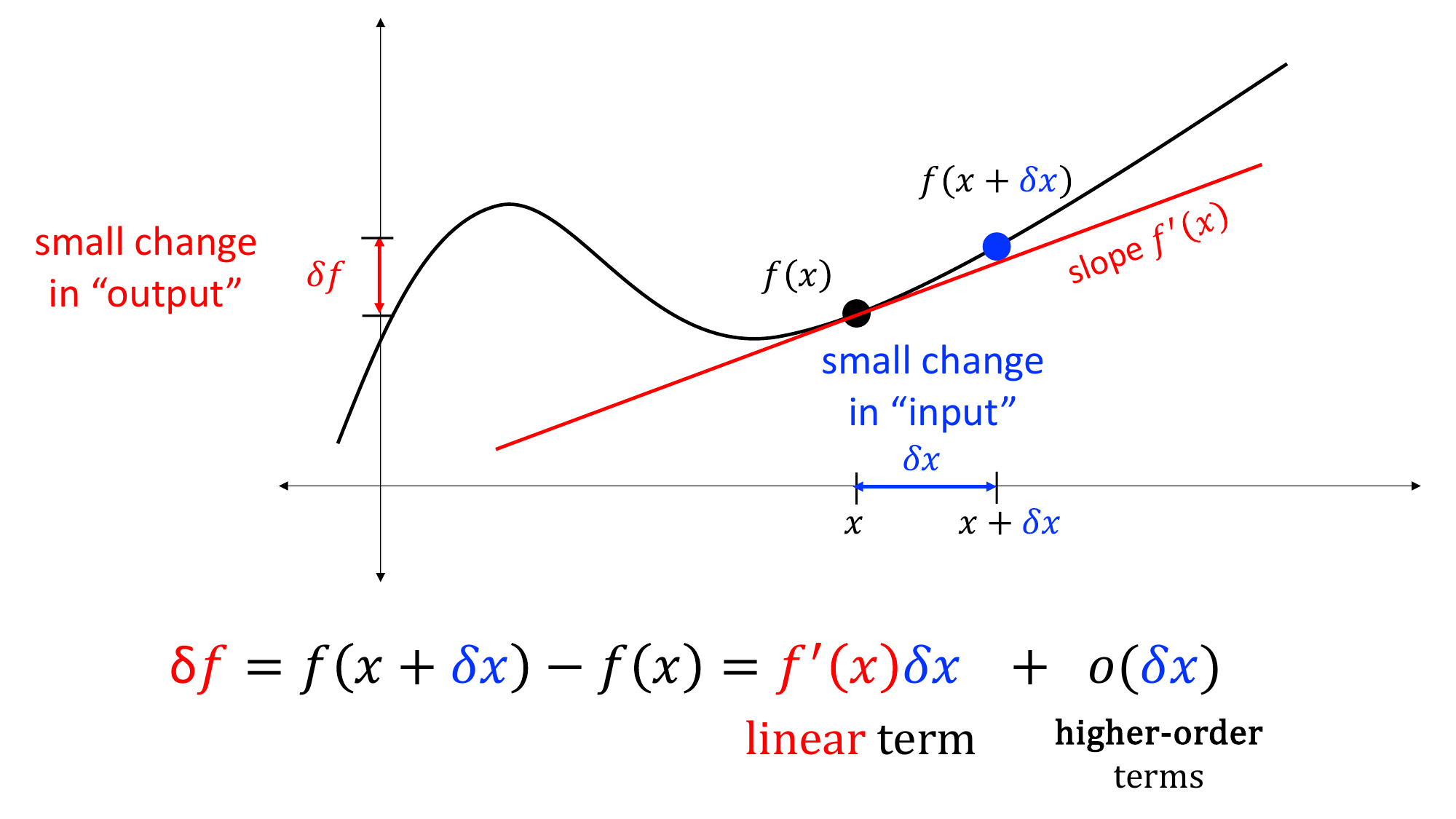}
    \caption{The essence of a derivative is \emph{linearization}: predicting a small change $\delta f$ in the output $f(x)$ from a small change $\delta x$ in the input $x$, to \emph{first order} in $\delta x$.}
    \label{fig:derivative-linearization}
\end{figure}

In a first-semester single-variable calculus course (like 18.01 at MIT), the derivative $f'(x)$ is  introduced as the slope of the tangent line at the point $(x,f(x))$, which can also be viewed as a \emph{linear approximation} of $f$ near $x$. In particular, as depicted in Fig.~\ref{fig:derivative-linearization}, this is equivalent to a prediction of the \emph{change $\delta f$ in the ``output''} of $f(x)$ from a small \emph{change $\delta x$ in the ``input''} to first order (\emph{linear}) in $\delta x$:
\[
\delta f = f(x+\delta x) - f(x) =  f'(x) \, \delta x  + \underbrace{(\text{higher-order terms})}_{o(\delta x)}.
\]
We can more precisely express these higher-order terms using asymptotic ``little-o'' notation ``$o(\delta x)$'', which denotes any function whose magnitude shrinks much faster than $|\delta x|$ as $\delta x \to 0 $, so that for sufficiently small $\delta x$ it is negligible compared to the linear $f'(x) \, \delta x$ term.  (Variants of this notation are commonly used in computer science, and there is a formal definition that we omit here.\footnote{Briefly, a function $g(\delta x)$ is $o(\delta x)$ if $\lim_{\delta x \to 0} \frac{\Vert g(\delta x) \Vert}{\Vert \delta x \Vert} = 0$.  We will return to this subject in Section~\ref{sec:banach}.})  Examples of such higher-order terms include $(\delta x)^2$, $(\delta x)^3$, $(\delta x)^{1.001}$, and $(\delta x)/\log(\delta x)$.

\begin{remark}
Here, $\delta x$ is \emph{not} an infinitesimal but rather a small number.  Note that our symbol ``$\delta$'' (a Greek lowercase ``delta'') is \emph{not} the same as the symbol ``$\partial$'' commonly used to denote partial derivatives.
\end{remark}

This notion of a derivative may remind you of the first two terms in a Taylor series $f(x+\delta x) = f(x) + f'(x) \, \delta x + \cdots$ (though in fact it is much more basic than Taylor series!), and the notation will generalize nicely to higher dimensions and other vector spaces. In differential notation, we can express the same idea as:
\[
\d f = f(x+\d x) - f(x) = f'(x) \, \d x.
\]
In this notation we implicitly drop the $o(\delta x)$ term that vanishes in the limit as $\delta x$ becomes infinitesimally small.

We will use this as the more generalized definition of a derivative.
 In this formulation, we avoid \emph{dividing} by $\d x$, because soon we will allow $x$ (and hence $\d x$) to be something other than a number---if $\d x$ is a vector, we won't be \emph{able} to divide by it!

\subsection{Linear operators}

From the perspective of linear algebra, given a function $f$, we consider the derivative of $f$, to be the \emph{linear operator} $f'(x)$ such that
\[
\d f = f(x+ \d x) - f(x) = f'(x)[\d x].
\]

As above, you should think of the differential notation $\d x$ as representing an \emph{arbitrary small} change in $x$, where we are implicitly dropping any $o(dx)$ terms, i.e.~terms that decay faster than linearly as $dx\to 0$.  Often, we will omit the square brackets and write simply $f'(x)\d x$ instead of $f'(x)[\d x]$, but this should be understood as the linear operator $f'(x)$ \emph{acting on} $\d x$---don't write $\d x\,f'(x)$, which will generally be nonsense!

This definition will allow us to extend differentiation to \emph{arbitrary vector spaces} of inputs $x$ and outputs $f(x)$.  (More technically, we will require vector spaces with a norm $\Vert x \Vert$, called ``Banach spaces,'' in order to precisely define the $o(\delta x)$ terms that are dropped.  We will come back to the subject of Banach spaces later.)

\begin{recall}[Vector Space]
    Loosely, a vector space (over $\R$) is a set of elements in which addition and subtraction between elements is defined, along with multiplication by real scalars. For instance, while it does not make sense to multiply arbitrary vectors in $\R^n$, we can certainly add them together, and we can certainly scale the vectors by a constant factor.
\end{recall}

Some examples of vector spaces include:
\begin{itemize}
    \item $\R^n$, as described in the above. More generally, $\R^{n\times m}$, the space of $n\times m$ matrices with real entries. Notice again that, if $n\neq m$, then multiplication between elements is not defined.
    \item $C^0(\R^n)$, the set of continuous functions over $\R^n$, with addition defined pointwise.
\end{itemize}

\begin{recall}[Linear Operator]
Recall that a linear operator is a map $L$ from a vector $v$ in vector space $V$ to a vector $L[v]$ (sometimes denoted simply $Lv$) in some other vector space. Specifically, $L$ is linear if
\[
L[v_1+v_2] = Lv_1 + Lv_2 \text{~~and~~} L[\alpha v] = \alpha L[v]
\]
for scalars $\alpha \in \mathbb{R}$.

\textit{Remark}: In this course, $f'$ is a map that takes in an $x$ and spits out a linear operator $f'(x)$ (the \textbf{derivative} of $f$ at $x$). Furthermore, $f'(x)$ is a linear map that takes in an input direction $v$ and gives an output vector $f'(x)[v]$ (which we will later interpret as a directional derivative, see Sec.~\ref{sec:directional}). When the direction $v$ is an infinitesimal $d x$, the output $f'(x)[dx] = df$ is the \textbf{differential} of $f$ (the corresponding infinitesimal change in $f$).
\end{recall}

\begin{notation}[Derivative operators and notations]

There are multiple notations for derivatives in common use, along with multiple related concepts of derivative, differentiation, and differentials.  In the table below, we summarize several of these notations, and put $\boxed{\mathrm{boxes}}$ around the notations adopted for this course:

\vspace{.2cm}

\begin{tabular}{| m{2.1cm} | m{5cm}| m{7cm} |}
\hline
name & notations & remark\tabularnewline
\hline
\hline
derivative
 & $\boxed{f'}$, also $\frac{df}{dx}$, $Df$, $f_{x}$, $\partial_{x}f$,
\dots{} & linear operator $f'(x)$ that maps a small change $dx$ in the input to a small
change ${df=f'(x)[dx]}$ in the output \\
& &
\\
&  & In single-variable calculus, this linear operator can be represented by a \emph{single number}, the ``slope,'' e.g. if $f(x) = \sin(x)$ then $f'(x) = \cos(x)$ is the number that we multiply by $dx$ to get $dy = \cos(x) dx$.  In multi-variable calculus, the linear operator $f'(x)$ can be represented by a \emph{matrix}, the {\bf Jacobian} $J$ (see Sec.~\ref{sec:kronecker}), so that $df = f'(x)[dx] = J\, dx$.  But we will see that it is not always convenient to express $f'$ as a matrix, even if we can.
\tabularnewline
\hline
differentiation & $\boxed{^{\prime}}$, $\frac{d}{dx}$, $D$, \dots{} & linear operator that maps a function $f$ to its derivative $f'$ \tabularnewline
\hline
difference & $\boxed{\delta x}$ and $\boxed{\delta f} = f(x+\delta x) - f(x)$  & small (but \emph{not} infinitesimal) change in the input  $x$ and output $f$ (depending implicitly on $x$ and $\delta x$), respectively:
an element of a vector space, \emph{not} a linear operator\tabularnewline
\hline
differential & $\boxed{d x}$ and $\boxed{d f} = f(x+dx) - f(x)$  & arbitrarily small (``infinitesimal''\footnote{Informally, one can think of the vector space of infinitesimals $\d x$ as living in the same space as $x$  (understood as a small change in a vector, but still a vector nonetheless). Formally, one can define a distinct ``vector space of infinitesimals'' in various ways, e.g.~as a cotangent space in differential geometry, though we won't go into more detail here.} --- we drop higher-order terms) change in the input  $x$ and output $f$, respectively:
an element of a vector space, \emph{not} a linear operator\tabularnewline
\hline
gradient & $\boxed{\nabla f}$ & the vector whose inner product $df = \langle \nabla f, dx \rangle$ with a small change $dx$ in the input gives the small change $df$ in the output.  The ``transpose of the derivative'' $\nabla f = (f')^T$.  (See Sec.~\ref{sec:multivarPart1}.)
\tabularnewline
\hline
partial derivative & $\boxed{\frac{\partial f}{\partial x}}$, $f_{x}$, $\partial_{x}f$ & linear operator that maps a small change $dx$ in a \emph{single argument}
of a multi-argument function to the corresponding change in output,
e.g. for $f(x,y)$ we have $df=\frac{\partial f}{\partial x}[dx]+\frac{\partial f}{\partial y}[dy]$.\tabularnewline
\hline
\end{tabular}
\end{notation}

\newpage

Some examples of linear operators include
\begin{itemize}
    \item Multiplication by scalars $\alpha$, i.e.~$Lv = \alpha v$.  Also  multiplication of column vectors $v$ by matrices $A$, i.e. $Lv = Av$.
    \item Some functions like $f(x)=x^2$ are obviously nonlinear.  But what about $f(x)=x+1$?  This may \emph{look} linear if you plot it, but it is \emph{not} a linear operation, because $f(2x)=2x+1\ne 2f(x)$---such functions, which are linear \emph{plus a nonzero constant}, are known as \emph{affine}.
    \item There are also many other examples of linear operations that are not so convenient or easy to write down as matrix--vector products.  For example, if $A$ is a $3\times 3$ matrix, then $L[A]=AB+CA$ is a linear operator given $3\times 3$ matrices $B,C$.  The transpose $f(x)=x^T$ of a column vector $x$ is linear, but is not given by any matrix multiplied by $x$. Or, if we consider vector spaces of \emph{functions}, then the calculus operations of differentiation and integration are linear operators too!
\end{itemize}

\subsubsection{Directional derivatives}
\label{sec:directional}

There is an equivalent way to interpret this linear-operator viewpoint of a derivative, which you may have seen before in multivariable calculus: as a \textbf{directional derivative}.

If we have a function $f(x)$ of arbitrary vectors $x$, then the directional derivative at $x$ in a direction (vector) $v$ is defined as:
\begin{equation}
\left. \frac{\partial}{\partial \alpha} f(x+\alpha v) \right|_{\alpha = 0} = \lim_{\delta\alpha\to 0} \frac{f(x+\delta\alpha \, v) - f(x)}{\delta \alpha}
\end{equation}
where $\alpha$ is a \emph{scalar}.  This transforms derivatives back into \emph{single-variable calculus} from arbitrary vector spaces.  It measures the rate of change of $f$ in the direction~$v$ from~$x$.  But it turns out that this has a very simple relationship to our linear operator $f'(x)$ from above, because (dropping higher-order terms due to the limit $\delta\alpha \to 0$):
$$
f(x + \underbrace{d\alpha \, v}_{d x}) - f(x) = f'(x)[dx] = d\alpha \, f'(x)[v] \, ,
$$
where we have factored out the scalar $d\alpha$ in the last step thanks to $f'(x)$ being a \emph{linear} operator.  Comparing with above, we immediately find that the directional derivative is:
\begin{equation}
\boxed{ \left. \frac{\partial}{\partial \alpha} f(x+\alpha v) \right|_{\alpha = 0} = f'(x)[v]  }\, .
\end{equation}
It is \emph{exactly equivalent} to our $f'(x)$ from before!   (We can also see this as an instance of the chain rule from Sec.~\ref{sec:chainrule}.) One lesson from this viewpoint is that it is perfectly reasonable to input an arbitrary \emph{non-infinitesimal} vector $v$ into $f'(x)[v]$: the result is not a $df$, but is simply a directional derivative.

\subsection{Revisiting multivariable calculus, Part 1: Scalar-valued functions} \label{sec:multivarPart1}

\begin{figure}
    \centering
    \includegraphics[trim={5cm 2cm 5cm 2cm},clip,width=1.0\textwidth]{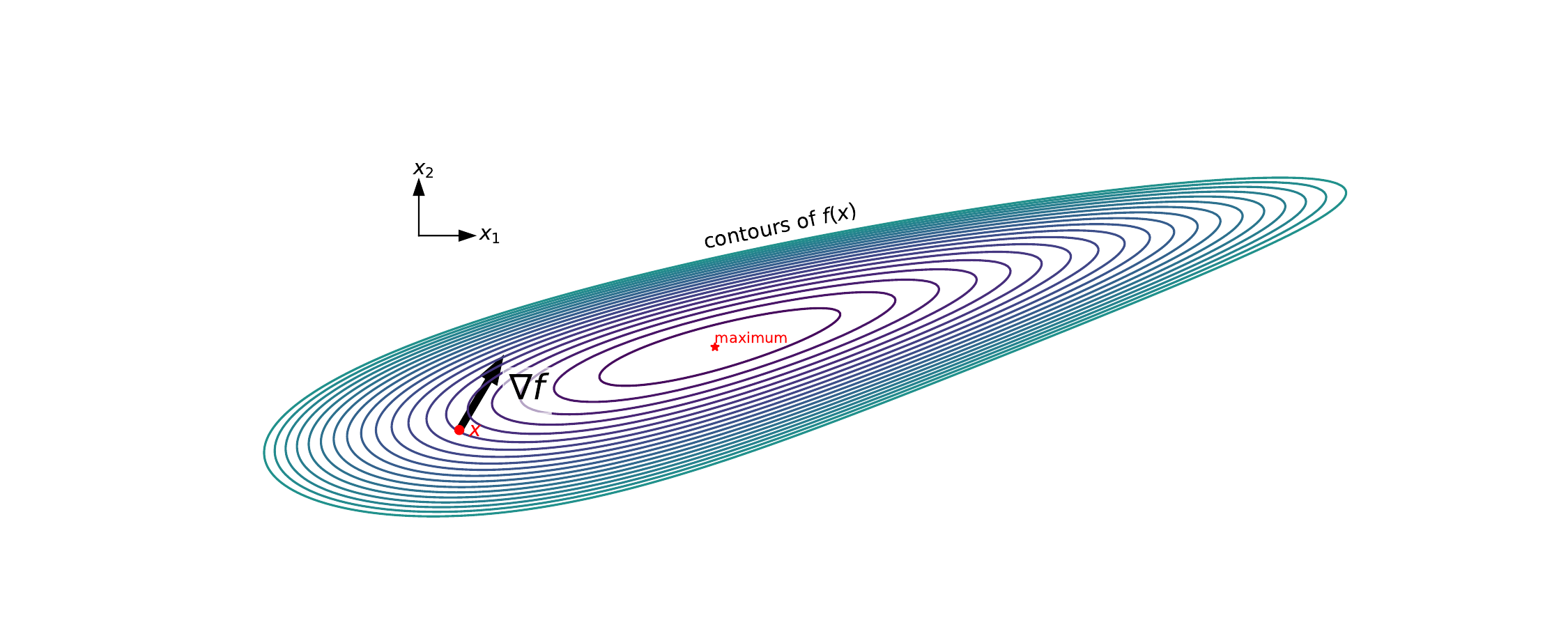}
    \caption{For a real-valued $f(x)$, the gradient $\nabla f$ is defined so that it corresponds to the ``uphill'' direction at a point $x$, which is perpendicular to the contours of $f$.  Although this may not point exactly towards the nearest local maximum of $f$ (unless the contours are circular), ``going uphill'' is nevertheless the starting point for many computational-optimization algorithms to search for a maximum.}
    \label{fig:gradient-uphill}
\end{figure}

Let $f$ be a scalar-valued function, which takes in ``column'' vectors $x \in \R^n$ and produces a scalar (in $\R$). Then,
\[
\d f = f(x+ \d x) - f(x) = f'(x) [\d x] = \text{scalar}.
\]
Therefore, since $\d x$ is a column vector (in an arbitrary direction, representing an arbitrary small change in $x$), the linear operator $f'(x)$ that produces a scalar $\d f$ must be a \textbf{row vector} (a ``1-row matrix'', or more formally something called a \emph{covector} or ``dual'' vector or ``linear form'')! We call this row vector the \emph{transpose of the gradient} $(\nabla f)^T$, so that $df$ is the \emph{dot (``inner'') product of $dx$ with the gradient}. So we have that
\[
\d f = \nabla f \cdot \d x = \underbrace{(\nabla f)^T}_{f'(x)} \d x \hspace{.25cm} \text{where} \hspace{.25cm} \d x = \begin{pmatrix}
\d x_1 \\  \d x_2 \\ \vdots \\ \d x_n.
\end{pmatrix} .
\]
Some authors view the gradient as a row vector (equating it with $f'$ or the Jacobian), but treating it as a ``column vector'' (the transpose of $f'$), as we do in this course, is a common and useful choice. As a column vector, the gradient can be viewed as the ``uphill'' (\emph{steepest-ascent}) direction in the $x$ space, as depicted in Fig.~\ref{fig:gradient-uphill}.  Furthermore, it is also easier to generalize to scalar functions of other vector spaces. In any case, for this class, we will \textit{always} define $\nabla f$ to \emph{have the same ``shape'' as $x$}, so that $\d f$ is a dot product (``inner product'') of $\d x$ with the gradient.

This is perfectly consistent with the viewpoint of the gradient that you may remember from multivariable calculus, in which the gradient was a vector of components
$$
\nabla f = \begin{pmatrix}
\frac{\partial f}{\partial x_1} \\  \frac{\partial f}{\partial x_2} \\ \vdots \\ \frac{\partial f}{\partial x_n}
\end{pmatrix} \, ;
$$
or, equivalently,
$$
\d f = f(x+\d x)-f(x) = \nabla f \cdot dx = \frac{\partial f}{\partial x_1} dx_1 + \frac{\partial f}{\partial x_2} dx_2 + \cdots + \frac{\partial f}{\partial x_n} dx_n \, .
$$
While a component-wise viewpoint may sometimes be convenient, we want to encourage you to view the vector $x$ as a \emph{whole}, not simply a collection of components, and to learn that it is often more convenient and elegant to differentiate expressions \emph{without} taking the derivative component-by-component, a new approach that will generalize better to more complicated inputs/output vector spaces.

Let's look at an example to see how we compute this differential.

\begin{example}
    Consider $f(x) = x^T Ax$ where $x\in \R^n$ and $A$ is a square $n\times n$ matrix, and thus $f(x) \in \R$. Compute $\d f$, $f'(x)$, and $\nabla f$.
\end{example}

We can do this directly from the definition.
\begin{align*}
    \d f &= f(x+ \d x) - f(x) \\
    &=  (x+ \d x)^T A (x+\d x) - x^T A x \\
    &= \cancel{x^T A x} + \d x^T\, Ax + x^T A \d x + \cancelto{\text{higher order}}{\d x^T \,A \d x} - \cancel{x^T A x} \\
    &= \underbrace{x^T(A+A^T)}_{f'(x)=(\nabla f)^T }\d x \implies \nabla f = (A+A^T)x \, .
\end{align*}
Here, we dropped terms with more than one $\d x$ factor as these are asymptotically negligible.  Another trick was to combine $\d x^T A x$ and $x^T A \d x$ by realizing that these are \emph{scalars} and hence equal to their own transpose: $dx^T A x = (\d x^T A x)^T = x^T A^T \d x$. Hence, we have found that $f'(x) = x^T(A+A^T) = (\nabla f)^T$, or equivalently $\nabla f = [x^T(A+A^T)]^T = (A+A^T)x$.

It is, of course, also possible to compute the same gradient component-by-component, the way you probably learned to do in multivariable calculus.  First, you would need to write $f(x)$ explicitly in terms of the components of $x$, as $f(x) = x^T A x = \sum_{i,j} x_i A_{i,j} x_j$.  Then, you would compute $\partial f/\partial x_k$ for each $k$, taking care that $x$ appears twice in the $f$ summation.  However, this approach is awkward, error-prone, labor-intensive, and quickly becomes worse as we move on to more complicated functions.  It is much better, we feel, to get used to treating vectors and matrices \emph{as a whole}, not as mere collections of numbers.

\subsection{Revisiting multivariable calculus, Part 2: Vector-valued functions}

Next time, we will revisit multi-variable calculus (18.02 at MIT) again in a Part 2, where now $f$ will be a vector-valued function, taking in vectors $x\in \R^n$ and giving vector outputs $f(x) \in \R^m$. Then, $\d f$ will be a $m$-component column vector, $\d x$ will be an $n$-component column vector, and we must get a linear operator $f'(x)$ satisfying
\[
\underbrace{\d f}_{m\text{ components}} = \underbrace{f'(x)}_{m \times n} \underbrace{\d x}_{n\text{ components}} \, ,
\]
so $f'(x)$ must be an $m \times n$ \emph{matrix} called the \textit{Jacobian} of $f$!

The Jacobian matrix $J$
represents the linear operator that takes $dx$ to $df$:
\[
\d f = J \d x \, .
\]
The matrix $J$ has entries $J_{ij}=\frac{\partial f_i}{\partial x_j}$ (corresponding to the $i$-th row and the $j$-th column of $J$).

So now, suppose that $f: \R^2 \to \R^2$. Let's understand how we would compute the differential of $f$:
\[
\d f = \begin{pmatrix}
   \frac{\partial f_1}{\partial x_1} & \frac{\partial f_1}{\partial x_2} \\ \frac{\partial f_2}{\partial x_1} & \frac{\partial f_2}{\partial x_2}
\end{pmatrix} \begin{pmatrix}
     \d x_1 \\ \d x_2
\end{pmatrix} = \begin{pmatrix}
    \frac{\partial f_1}{\partial x_1}\d x_1 +  \frac{\partial f_1}{\partial x_2}\d x_2 \\ \frac{\partial f_2}{\partial x_1}\d x_1 +  \frac{\partial f_2}{\partial x_2} \d x_2
\end{pmatrix}.
\]

Let's compute an example.
\begin{example}\label{ex1}
    Consider the function $f(x) = Ax$ where $A$ is a constant $m\times n$ matrix. Then, by applying the distributive law for matrix--vector products, we have
\begin{align*}
    \d f &= f(x+\d x ) - f(x) =  A(x + \d x) - Ax  \\
    &= \cancel{Ax} + A\d x - \cancel{Ax} = A\d x = f'(x) \d x.
\end{align*}
Therefore, $f'(x) = A$.
\end{example}

Notice then that the linear operator $A$ is its own Jacobian matrix!

Let's now consider some derivative rules.

\begin{itemize}
    \item \textbf{Sum Rule}: Given $f(x) = g(x) + h(x)$, we get that
    \[
\d f = \d g + \d h \implies f'(x) \d x = g'(x) \d x + h'(x)\d x.
    \]
    Hence, $f' = g'+h'$ as we should expect.  This is the linear operator $f'(x)[v] = g'(x)[v] + h'(x)[v]$, and note that we can sum linear operators (like $g'$ and $h'$) just like we can sum matrices! In this way, linear operators form a vector space.
    \item \textbf{Product Rule}: Suppose $f(x) = g(x) h(x)$. Then,
    \begin{align*}
        \d f &= f(x+ \d x) - f(x) \\
        &= g(x+\d x) h(x+\d x) - g(x) h(x) \\
        &= (g(x) + \underbrace{g'(x) \d x}_{\d g}) (h(x) + \underbrace{h'(x)\d x}_{\d h}) - g(x)h(x) \\
        &= g h + dg\, h + g\,dh + \cancelto{0}{dg \,dh} - gh \\
        &= \d g \, h + g \,\d h \, ,
    \end{align*}
    where the $\d g \, \d h$ term is higher-order and hence dropped in infinitesimal notation.
    Note, as usual, that $\d g$ and $h$ may not commute now as they may no longer be scalars!
\end{itemize}

Let's look at some short examples of how we can apply the product rule nicely.

\begin{example}\label{ex:Ax}
    Let $f(x) = Ax$ (mapping $\R^n \to \R^m$) where $A$ is a constant $m\times n$ matrix. Then,
    \[
    \d f = d(Ax) =  \cancelto{0}{\d A} \, x + A\d x = A\d x \implies f'(x) = A.
    \]
    We have ${\d A} = 0$ here because $A$ does not change when we change $x$.
\end{example}

\begin{example}
    Let $f(x) = x^T A x$ (mapping $\R^n \to \R$). Then,
    \[
    \d f = \d x^T(Ax) + x^T \d(Ax) = \underbrace{\d x^T \, Ax}_{=\, x^T A^T\d x} + x^T A\d x = x^T(A+A^T)\d x  = (\nabla f)^T \d x\, ,
    \]
    and hence $f'(x) = x^T(A + A^T)$. (In the common case where $A$ is symmetric, this simplifies to $f'(x) = 2x^T A$.) Note that here we have applied Example~\ref{ex:Ax} in computing $d(Ax) = A dx$.
    Furthermore, $f$ is a scalar valued function, so we may also obtain the gradient $\nabla f = (A+A^T) x$ as before (which simplifies to $2Ax$ if $A$ is symmetric).
\end{example}

\begin{example}[Elementwise Products]
    Given $x,y\in \R^m$, define
    $$x\dotstar y = \begin{pmatrix}
        x_1y_1 \\ \vdots \\ x_my_m
    \end{pmatrix} = \underbrace{\begin{pmatrix} x_1 & & & \\ & x_2 & & \\ & & \ddots & \\ & & & x_m \end{pmatrix}}_{\mathrm{diag}(x)} \begin{pmatrix}
        y_1 \\ \vdots \\ y_m
    \end{pmatrix},$$
    the element-wise product of vectors (also called the \href{https://en.wikipedia.org/wiki/Hadamard_product_(matrices)}{Hadamard product}), where for convenience below we also define $\mathrm{diag}(x)$ as the $m \times m$ diagonal matrix with $x$ on the diagonal. Then, given $A \in \R^{m,n}$, define $f:\R^n \to \R^m$ via
    \[
    f(x) = A(x\dotstar x).
    \]

    As an exercise, one can verify the following:
    \begin{itemize}
        \item[(a)] $x\dotstar y = y\dotstar x$,
        \item[(b)] $A(x\dotstar y) = A \, \mathrm{diag}(x)\, y$.
        \item[(c)] $d(x\dotstar y) = (\d x)\dotstar y + x\dotstar(\d y)$.  So if we take $y$ to be a constant and define $g(x) = y \dotstar x$, its Jacobian matrix is $\mathrm{diag}(y)$.
        \item[(d)] $\d f = A(2x\dotstar dx) = 2A \, \mathrm{diag}(x) \, dx = f'(x)[dx]$, so the Jacobian matrix is $J = 2A \, \mathrm{diag}(x)$.
        \item[(e)] Notice that the directional derivative (Sec.~\ref{sec:directional}) of $f$ at $x$ in the direction $v$ is simply given by $f'(x)[v] = 2A (x \dotstar v)$.   One could also check numerically for some arbitrary $A,x,v$ that $f(x + 10^{-8} v) - f(x) \approx 10^{-8} ( 2A (x \dotstar v))$.
    \end{itemize}
\end{example}

\subsection{The Chain Rule}
\label{sec:chainrule}

One of the most important rules from differential calculus is the chain rule, because it allows us to differentiate complicated functions built out of compositions of simpler functions.   This chain rule can also be generalized to our differential notation in order to work for functions on arbitrary vector spaces:
\begin{itemize}
    \item \textbf{Chain Rule}: Let $f(x) = g(h(x))$. Then,
    \begin{align*}
    \d f = f(x + \d x) - f(x) &= g(h( x + \d x)) - g(h(x)) \\
    &= g'(h(x))[h(x + \d x) - h(x)] \\
    &= g'(h(x))[h'(x)[\d x]] \\
    &= g'(h(x)) h'(x) [\d x]
    \end{align*}
    where $g'(h(x))h'(x)$ is a composition of $g'$ and $h'$ as matrices.

    In other words, $f'(x) = g'(h(x)) h'(x)$: the Jacobian (linear operator) $f'$ is simply the \emph{product (composition) of the Jacobians}, $g' h'$.  Ordering matters because linear operators do not generally commute: left-to-right = outputs-to-inputs.
\end{itemize}

Let's look more carefully at the \emph{shapes} of these Jacobian matrices in an example where each function maps a column vector to a column vector:
\begin{example}
    Let $x\in \R^n$, $h(x) \in \R^p$, and $g(h(x)) \in \R^m$. Then, let $f(x) = g(h(x))$ mapping from $\R^n$ to $\R^m$. The chain rule then states that
    \[
    f'(x) = g'(h(x)) h'(x),
    \]
    which makes sense as $g'$ is an $m\times p$ matrix and $h'$ is a $p\times n$ matrix, so that the product gives an $m \times n$ matrix $f'$! However, notice that this is \textit{not} the same as $h'(x) g'(h(x))$ as you cannot (if $n\neq m$) multiply a $p\times n$ and an $m\times p$ matrix together, and even if $n = m$ you will get the wrong answer since they probably won't commute.
\end{example}

Not only does the order of the multiplication matter, but the associativity of matrix multiplication matters \textit{practically}. Let's consider a function
\[
f(x)= a(b(c(x)))
\]
where $c: \R^n \to \R^p$, $b:\R^p \to \R^q$, and $a: \R^q\to \R^m$. Then, we have that, by the chain rule,
\[
f'(x) = a'(b(c(x))) b'(c(x))c'(x).
\]
Notice that this is the same as
\[
f'  = (a' b')c' = a'(b'c')
\]
by associativity (omitting the function arguments for brevity). The left-hand side is multiplication from left to right, and the right-hand side is multiplication from right to left.

But who cares? Well it turns out that associativity is deeply important. So important that the two orderings have names: multiplying  left-to-right is called ``reverse mode''  and multiplying right-to-left is called ``forward mode'' in the field of \textit{automatic differentiation} (AD).  Reverse-mode differentation is also known as an ``adjoint method'' or ``backpropagation'' in some contexts, which we will explore in more detail later. Why does this matter? Let's think about the computational cost of matrix multiplication.

\subsubsection{Cost of Matrix Multiplication}
\label{sec:cost-matrix-mult}

\begin{figure}
    \centering
    \includegraphics[width=0.8\textwidth]{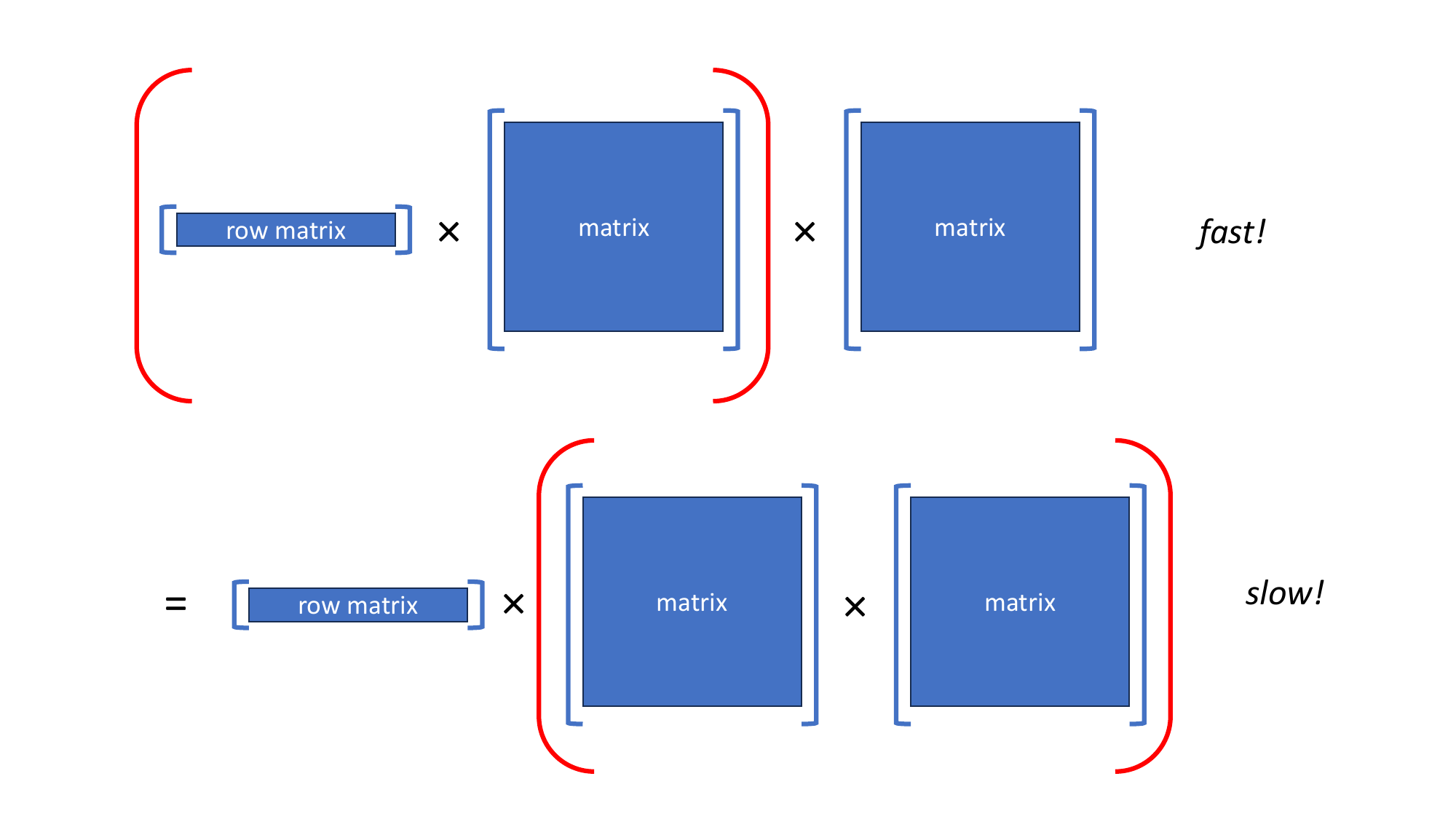}
    \caption{Matrix multiplication is \emph{associative}---that is, $(AB)C = A(BC)$ for all $A,B,C$---but multiplying left-to-right can be much more efficient than right-to-left if the leftmost matrix has only one (or few) rows, as shown here.  Correspondingly, the order in which you carry out the chain rule has dramatic consequences for the computational effort required.  Left-to-right is known as ``reverse mode'' or ``backpropagation'', and is best suited to situations where there are many fewer outputs than inputs.}
    \label{fig:matrix-mult-assoc}
\end{figure}

If you multiply a $m\times q$ matrix by a $q\times p$ matrix, you normally do it by computing $mp$ dot products of length $q$ (or some equivalent re-ordering of these operations). To do a dot product of length $q$ requires $q$ multiplications and $q-1$ additions of scalars. Overall, this is approximately $2mpq$ scalar operations in total. In computer science, you would write that this is ``$\Theta(mpq)$'': the computational effort is \emph{asymptotically proportional} to $mpq$ for large $m,p,q$.

So why does the order of the chain rule matter? Consider the following two examples.
\begin{example}
    Suppose you have a lot of inputs $n \gg 1$, and only one output $m=1$, with lots of intermediate values, i.e. $q=p=n$. Then reverse mode (left-to-right) will cost $\Theta(n^2)$ scalar operations while forward mode (right-to-left) would cost $\Theta(n^3)$! This is a huge cost difference, depicted schematically in Fig.~\ref{fig:matrix-mult-assoc}.

    Conversely, suppose you have a lot of outputs $m \gg 1$ and only one input $n=1$, with lots of intermediate values $q=p=m$.  Then reverse mode would cost $\Theta(m^3)$ operations but forward mode would be only $\Theta(m^2)$!

    Moral: If you have a lot of inputs and few outputs (the usual case in machine learning and optimization), compute the chain rule left-to-right (reverse mode).  If you have a lot of outputs and few inputs, compute the chain rule right-to-left (forward mode).  We return to this in Sec.~\ref{sec:forward-vs-reverse}.
\end{example}

\subsection{Beyond Multi-Variable Derivatives}

Now let's compute some derivatives that go beyond first-year calculus, where the inputs and outputs are in more general vector spaces. For instance, consider the following examples:

\begin{example}
    Let $A$ be an $n\times n$ matrix. You could have the following matrix-valued functions. For example:
\begin{itemize}
    \item $f(A) = A^3$,
    \item $f(A) = A^{-1}$ if $A$ is invertible,
    \item or $U$, where $U$ is the resulting matrix after applying Gaussian elimination to $A$!
\end{itemize}

\noindent You could also have scalar outputs. For example:
\begin{itemize}
    \item $f(A) = \det A$,
    \item $f(A) =$ trace $A$,
    \item or $f(A) = \sigma_1(A)$, the largest singular value of $A$.
\end{itemize}
\end{example}

Let's focus on two simpler examples for this lecture.
\begin{example}
    Let $f(A) = A^3$ where $A$ is a square matrix. Compute $\d f$.
\end{example}

Here, we apply the chain rule one step at a time:
\[
\d f = \d A\, A^2 + A \,\d A \, A + A^2 \,\d A = f'(A) [\d A].
\]
Notice that this is not equal to $3A^2$ (unless $\d A$ and $A$ commute, which won't generally be true since $\d A$ represents an \emph{arbitrary} small change in $A$).   The right-hand side is a linear operator $f'(A)$ acting on $\d A$, but it is not so easy to interpret it as simply a single ``Jacobian'' matrix multiplying $\d A$!

\begin{example}
    Let $f(A) = A^{-1}$ where $A$ is a square invertible matrix. Compute $\d f = \d(A^{-1})$.
\end{example}

Here, we use a slight trick. Notice that $AA^{-1} = \Id$, the identity matrix. Thus, we can compute the differential using the product rule (noting that $\d\Id = 0$, since changing $A$ does not change $\Id$) so
\[
\d(AA^{-1}) = \d A\, A^{-1} + A\,\d(A^{-1}) = \d (\Id) = 0\implies \d(A^{-1}) = -A^{-1} \,\d A \,A^{-1}.
\]

\pagebreak

\section{Jacobians of Matrix Functions}
\label{sec:kronecker}
When we have a function that has \emph{matrices} as inputs and/or
outputs, we have already seen in the previous lectures that we can
still define the derivative as a linear operator by a \emph{formula}
for $f'$ mapping a small change in input to the corresponding small change in output. However,
when you first learned linear algebra, probably most linear operations
were represented by matrices multiplying vectors, and it may take
a while to get used to thinking of linear operations more generally.
In this chapter, we discuss how it is still \emph{possible} to represent
$f'$ by a \textbf{Jacobian matrix} even for matrix inputs/outputs,
and how the most common technique to do this involves \textbf{matrix
``vectorization''} and a new type of matrix operation, a \textbf{Kronecker
product}. This gives us another way to think about our $f'$ linear
operators that is occasionally convenient, but at the same time it
is important to become comfortable with other ways of writing down
linear operators too---sometimes, the explicit Jacobian-matrix approach
can obscure key structure, and it is often computationally inefficient as well.

For this section of the notes, we refer to the linked \href{https://rawcdn.githack.com/mitmath/matrixcalc/3f6758996e40c5c1070279f89f7f65e76e08003d/notes/2x2Jacobians.jl.html}{Pluto Notebook}
for computational demonstrations of this material in Julia, illustrating
multiple views of the derivative of the square $A^{2}$ of $2\times2$
matrices $A$.

\subsection{Derivatives of matrix functions: Linear operators}

As we have already emphasized, the derivative $f'$ is the linear
operator that maps a small change in the input to a small change in
the output. This idea can take an unfamiliar form, however, when applied
to functions $f(A)$ that map matrix inputs $A$ to matrix outputs.
For example, we've already considered the following functions on square
$m\times m$ matrices:
\begin{itemize}
\item $f(A)=A^{3}$, which gives $df=f'(A)[dA]=dA\,A^{2}+A\,dA\,A+A^{2}\,dA$.
\item $f(A)=A^{-1}$, which gives $df=f'(A)[dA]=-A^{-1}\,dA\,A^{-1}$
\end{itemize}
\begin{example}An even simpler example is the \emph{matrix-square}
function:

\[
f(A)=A^{2}\,,
\]
which by the product rule gives 
\[
df=f'(A)[dA]=dA\,A+A\,dA\,.
\]
You can also work this out explicitly from $df=f(A+dA)-f(A)=(A+dA)^{2}-A^{2}$,
dropping the $(dA)^{2}$ term.\end{example}

In all of these examples, $f'(A)$ is described by a simple formula
for $f'(A)[dA]$ that relates an arbitrary change $dA$ in $A$ to
the change $df=f(A+dA)-f(A)$ in $f$, to first order. If the differential
is distracting you, realize that we can plug any matrix $X$ we want
into this formula, not just an ``infinitesimal'' change $dA$, e.g.~in
our matrix-square example we have 
\[
f'(A)[X]=XA+AX
\]
 for an arbitrary $X$ (a directional derivative, from Sec.~\ref{sec:directional}). This is \emph{linear} in $X$: if we scale
or add inputs, it scales or adds outputs, respectively:
\[
f'(A)[2X]=2XA+A\,2X=2(XA+AX)=2f'(A)[X]\,,
\]
\begin{align*}
f'(A)[X+Y] & =(X+Y)A+A(X+Y)=XA+YA+AX+AY=XA+AX+YA+AY\\
 & =f'(A)[X]+f'(A)[Y]\,.
\end{align*}
This is a perfectly good way to define a linear operation! We are
\emph{not} expressing it here in the familiar form $f'(A)[X]=(\text{some matrix?})\times(X\text{ vector?})$,
and that's okay! A formula like $XA+AX$ is easy to write down, easy
to understand, and easy to compute with. 

But sometimes you still may want to think of $f'$ as a single ``Jacobian''
matrix, using the most familiar language of linear algebra, and it
is possible to do that! If you took a sufficiently abstract linear-algebra
course, you may have learned that \emph{any} linear operator can be
represented by a matrix once you choose a basis for the input and
output vector spaces. Here, however, we will be much more concrete,
because there is a conventional ``Cartesian'' basis for matrices
$A$ called ``vectorization'', and in this basis linear operators
like $AX+XA$ are particularly easy to represent in matrix form once
we introduce a new type of matrix product that has widespread applications
in ``multidimensional'' linear algebra.

\subsection{A simple example: The two-by-two matrix-square function}

To begin with, let's look in more detail at our matrix-square function
\[
f(A)=A^{2}
\]
for the simple case of $2\times2$ matrices, which are described by
only four scalars, so that we can look at every term in the derivative
explicitly. In particular,

\begin{example} For a $2\times2$ matrix 
\[
A=\begin{pmatrix}p & r\\
q & s
\end{pmatrix},
\]
the matrix-square function is 
\[
f(A)=A^{2}=\begin{pmatrix}p & r\\
q & s
\end{pmatrix}\begin{pmatrix}p & r\\
q & s
\end{pmatrix}=\begin{pmatrix}p^{2}+qr & pr+rs\\
pq+qs & qr+s^{2}
\end{pmatrix}.
\]
\end{example}

Written out explicitly in terms of the matrix entries $(p,q,r,s)$
in this way, it is natural to think of our function as mapping 4
scalar inputs to 4 scalar outputs. That is, we can think of $f$
as equivalent to a ``vectorized'' function $\tilde{f}:\R^{4}\to\R^{4}$
given by 
\[
\tilde{f}(\left(\begin{array}{c}
p\\
q\\
r\\
s
\end{array}\right))=\left(\begin{array}{c}
p^{2}+qr\\
pq+qs\\
pr+rs\\
qr+s^{2}
\end{array}\right)\,.
\]
Converting a matrix into a column vector in this way is called \textbf{vectorization},
and is commonly denoted by the operation ``$\vecm$'':
\begin{align*}
\vecm A & =\vecm\begin{pmatrix}p & r\\
q & s
\end{pmatrix}=\kbordermatrix{ & \\
A_{1,1} & p\\
A_{2,1} & q\\
A_{1,2} & r\\
A_{2,2} & s
} \, ,\\
\vecm f(A) & =\vecm\begin{pmatrix}p^{2}+qr & pr+rs\\
pq+qs & qr+s^{2}
\end{pmatrix}=\left(\begin{array}{c}
p^{2}+qr\\
pq+qs\\
pr+rs\\
qr+s^{2}
\end{array}\right) \, .
\end{align*}
In terms of $\vecm$, our ``vectorized'' matrix-squaring function
$\tilde{f}$ is defined by 
\[
\tilde{f}(\vecm A)=\vecm f(A)=\vecm(A^{2})\,.
\]
More generally, \begin{definition}The \textbf{vectorization} $\vecm A\in\mathbb{R}^{mn}$
of any $m\times n$ matrix $A\in\mathbb{R}^{m\times n}$ is a defined
by simply \textbf{stacking the columns} of $A$, from left to right,
into a column vector $\vecm A$. That is, if we denote the $n$ columns
of $A$ by $m$-component vectors $\vec{a}_{1},\vec{a}_{2},\ldots\in\mathbb{R}^{m}$,
then
\[
\vecm A=\vecm\underbrace{\left(\begin{array}{cccc}
\vec{a}_{1} & \vec{a}_{2} & \cdots & \vec{a}_{n}\end{array}\right)}_{A\in\mathbb{R}^{m\times n}}=\left(\begin{array}{c}
\vec{a}_{1}\\
\vec{a}_{2}\\
\vdots\\
\vec{a}_{n}
\end{array}\right)\in\mathbb{R}^{mn}
\]
is an $mn$-component column vector containing all of the entries of
$A$.

On a computer, matrix entries are typically stored in a consecutive
sequence of memory locations, which can be viewed a form of vectorization.
In fact, $\vecm A$ corresponds exactly to what is known as ``column-major''
storage, in which the column entries are stored consecutively; this
is the default format in Fortran, Matlab, and Julia, for example,
and the venerable Fortran heritage means that column major is widely used in
linear-algebra libraries. \end{definition}.

\begin{problem}The vector $\vecm A$ corresponds to the coefficients
you get when you express the $m\times n$ matrix $A$ in a \emph{basis}
of matrices. What is that basis? \end{problem}

Vectorization turns unfamilar things (like matrix functions and derivatives
thereof) into familiar things (like vector functions and Jacobians
or gradients thereof). In that way, it can be a very attractive tool,
almost \emph{too} attractive---why do ``matrix calculus'' if you
can turn everything back into ordinary multivariable calculus? Vectorization
has its drawbacks, however: conceptually, it can obscure the underlying
mathematical structure (e.g. $\tilde{f}$ above doesn't look much
like a matrix square $A^{2}$), and computationally this loss of structure
can sometimes lead to severe inefficiencies (e.g. forming huge $m^2\times m^2$
Jacobian matrices as discussed below). Overall, we believe that the
\emph{primary} way to study matrix functions like this should be to
view them as having matrix inputs ($A$) and matrix outputs ($A^{2}$), and that one should likewise generally view the derivatives as linear operators on matrices,
not vectorized versions thereof. However, it is still useful to be
familiar with the vectorization viewpoint in order to have the benefit
of an alternative perspective.

\subsubsection{The matrix-squaring four-by-four Jacobian matrix}

To understand Jacobians of functions (from matrices to matrices),
let's begin by considering a basic question: \begin{question} What is the
\emph{size} of the Jacobian of the matrix-square function? \end{question}

Well, if we view the matrix squaring function via its vectorized equivalent
$\tilde{f}$, mapping $\mathbb{R}^{4}\mapsto\mathbb{R}^{4}$ (4-component
column vectors to 4-component column vectors), the Jacobian would
be a $4\times4$ matrix (formed from the derivatives of each output
component with respect to each input component). Now let's think about
a more general square matrix $A$: an $m\times m$ matrix. If we wanted
to find the Jacobian of $f(A)=A^{2}$, we could do so by the same
process and (symbolically) obtain an $m^{2}\times m^{2}$ matrix (since
there are $m^{2}$ inputs, the entries of $A$, and $m^{2}$ outputs,
the entries of $A^{2}$). Explicit computation of these $m^{4}$ partial derivatives
is rather tedious even for small $m$, but is a task that symbolic computational
tools in e.g.~Julia or Mathematica can handle. In fact, as seen in the
\href{https://rawcdn.githack.com/mitmath/matrixcalc/3f6758996e40c5c1070279f89f7f65e76e08003d/notes/2x2Jacobians.jl.html}{Notebook},
Julia spits out the Jacobian quite easily. For the $m=2$ case that we
wrote out explicitly above, you can either take the derivative of
$\tilde{f}$ by hand or use Julia's symbolic tools to obtain the Jacobian:
\begin{equation}
\tilde{f}'=\kbordermatrix{& (1,1) & (2,1) & (1,2) & (2,2) \\
(1,1) & 2p & r & q & 0\\
(2,1) & q & p+s & 0 & q\\
(1,2) & r & 0 & p+s & r\\
(2,2) & 0 & r & q & 2s
} \,.
\end{equation}
For example, the first row of $\tilde{f}'$ consists of the partial
derivatives of $p^{2}+qr$ (the first output) with respect to the
4 inputs $p,q,r,\mbox{and }s$.   Here, we have labeled the rows by the (row,column) indices $(j_\mathrm{out}, k_\mathrm{out})$ of the entries in the ``output'' matrix $d(A^2)$, and have labeled the columns by the indices $(j_\mathrm{in}, k_\mathrm{in})$ of the entries in the ``input'' matrix $A$.  Although we have written the Jacobian $\tilde{f}'$ as a ``2d'' matrix, you can therefore also imagine it to be a ``4d'' matrix indexed by $j_\mathrm{out}, k_\mathrm{out}, j_\mathrm{in}, k_\mathrm{in}$.

However, the matrix-calculus approach of viewing the derivative $f'(A)$
as a \emph{linear transformation on matrices} (as we derived above),
\[
f'(A)[X]=XA+AX\,,
\]
seems to be much more revealing than writing out an explicit component-by-component
``vectorized'' Jacobian $\tilde{f}'$, and gives a formula for any
$m\times m$ matrix without   laboriously requiring us to take $m^{4}$ partial
derivatives one-by-one. If we really want to pursue the vectorization perspective,
we need a way to recapture some of the structure that is obscured
by tedious componentwise differentiation. A key tool to bridge the
gap between the two perspectives is a type of matrix operation that
you may not be familiar with: \textbf{Kronecker products} (denoted
$\otimes$).

\subsection{Kronecker Products}

A linear operation like $f'(A)[X]=XA+AX$ can be thought of as a ``higher-dimensional
matrix:'' ordinary ``2d'' matrices map ``1d'' column vectors to 1d
column vectors, whereas to map 2d matrices to 2d matrices you might
imagine a ``4d'' matrix (sometimes called a \emph{tensor}). To transform
2d matrices back into 1d vectors, we already saw the concept of vectorization
($\vecm A$). A closely related tool, which transforms ``higher dimensional''
linear operations on matrices back into ``2d'' matrices for the vectorized
inputs/outputs, is the Kronecker product $A\otimes B$. Although they
don't often appear in introductory linear-algebra courses, Kronecker
products show up in a wide variety of mathematical applications where
multidimensional data arises, such as multivariate statistics and
data science or multidimensional scientific/engineering problems.

\begin{definition}If $A$ is an $m\times n$ matrix with entries
$a_{ij}$ and $B$ is a $p\times q$ matrix, then their \textbf{Kronecker
product} $A\otimes B$ is defined by
\[
A=\left(\begin{array}{ccc}
a_{11} & \cdots & a_{1n}\\
\vdots & \ddots & \vdots\\
a_{m1} & \cdots & a_{mn}
\end{array}\right)\Longrightarrow\underbrace{A}_{m\times n}\otimes\underbrace{B}_{p\times q}=\underbrace{\left(\begin{array}{ccc}
a_{11}B & \cdots & a_{1n}B\\
\vdots & \ddots & \vdots\\
a_{m1}B & \cdots & a_{mn}B
\end{array}\right)}_{mp\times nq}\,,
\]
so that $A\otimes B$ is an $mp\times nq$ matrix formed by ``pasting
in'' a copy of $B$ multiplying every element of $A$. \end{definition}
For example, consider $2\times2$ matrices
\[
A=\begin{pmatrix}p & r\\
q & s
\end{pmatrix}\text{~~and~~}B=\begin{pmatrix}a & c\\
b & d
\end{pmatrix} \, .
\]
Then $A\otimes B$ is a $4\times4$ matrix containing all possible
products of entries $A$ with entries of $B$. Note that $A\otimes B\ne B\otimes A$
(but the two are related by a re-ordering of the entries): 
\[
A\otimes B=\begin{pmatrix}p\red{B} & rB\\
qB & sB
\end{pmatrix}=\begin{pmatrix}p\red{a} & p\red{c} & ra & rc\\
p\red{b} & p\red{d} & rb & rd\\
qa & qc & sa & sc\\
qb & qd & sb & sd
\end{pmatrix}\qquad\ne\qquad B\otimes A=\begin{pmatrix}aA & cA\\
bA & dA
\end{pmatrix}=\begin{pmatrix}\red{a}p & ar & \red{c}p & cr\\
aq & as & cq & cs\\
\red{b}p & br & \red{d}p & dr\\
bq & bs & dq & ds
\end{pmatrix} \, ,
\]
where we've colored one copy of $B$ red for illustration.
See the \href{https://rawcdn.githack.com/mitmath/matrixcalc/3f6758996e40c5c1070279f89f7f65e76e08003d/notes/2x2Jacobians.jl.html}{Notebook}
for more examples of Kronecker products of matrices (including some
with pictures rather than numbers!). 

Above, we saw that $f(A)=A^{2}$ at $A=\begin{pmatrix}p & r\\
q & s
\end{pmatrix}$ could be thought of as an equivalent function $\tilde{f}(\vecm A)$
mapping column vectors of 4 inputs to 4 outputs ($\mathbb{R}^{4}\mapsto\mathbb{R}^{4}$),
with a $4\times4$ Jacobian that we (or the computer) laboriously
computed as 16 element-by-element partial derivatives. It turns out
that this result can be obtained \emph{much} more elegantly once we
have a better understanding of Kronecker products. We will find that
the $4\times4$ ``vectorized'' Jacobian is simply
\[
\tilde{f}'=\Id_{2}\otimes A+A^{T}\otimes\Id_{2}\,,
\]
where $\Id_{2}$ is the $2\times2$ identity matrix. That is, the
matrix linear operator $f'(A)[dA]=dA\,A+A\,dA$ is equivalent, after
vectorization, to:

\[
\vecm\underbrace{f'(A)[dA]}_{dA\,A+A\,dA}=\underbrace{(\Id_{2}\otimes A+A^{T}\otimes\Id_{2})}_{\tilde{f}'}\vecm dA=\underbrace{\begin{pmatrix}2p & r & q & 0\\
q & p+s & 0 & q\\
r & 0 & p+s & r\\
0 & r & q & 2s
\end{pmatrix}}_{\tilde{f}'}\underbrace{\begin{pmatrix}dp\\
dq\\
dr\\
ds
\end{pmatrix}}_{\vecm dA}.
\]
In order to understand \emph{why} this is the case, however, we must
first build up some understanding of the algebra of Kronecker products.
To start with, a good exercise is to convince yourself of a few simpler
properties of Kronecker products: \begin{problem} From the definition
of the Kronecker product, derive the following identities:
\begin{enumerate}
\item $(A\otimes B)^{T}=A^{T}\otimes B^{T}$.
\item $(A\otimes B)(C\otimes D)=(AC)\otimes(BD)$.
\item $(A\otimes B)^{-1}=A^{-1}\otimes B^{-1}$.  (Follows from property~2.)
\item $A\otimes B$ is orthogonal (its transpose is its inverse) if $A$
and $B$ are orthogonal. (From properties 1~\&~3.)
\item $\det(A\otimes B)=\det(A)^{m}\det(B)^{n}$, where $A\in\R^{n,n}$
and $B\in\R^{m,m}$.
\item $\tr(A\otimes B)=(\tr A)(\tr B)$.
\item Given eigenvectors/values $Au=\lambda u$ and $Bv=\mu v$ of $A$
and $B$, then $\lambda\mu$ is an eigenvalue of $A\otimes B$ with
eigenvector $u\otimes v$. (Also, since $u\otimes v=\vecm X$ where
$X=vu^{T}$, you can relate this via Prop.~\ref{prop5000} below to the identity $BXA^{T}=Bv(Au)^T=\lambda\mu X$.) 
\end{enumerate}
\end{problem}

\subsubsection{Key Kronecker-product identity}

In order to convert linear operations like $AX+XA$ into Kronecker
products via vectorization, the key identity is:

\begin{proposition}\label{prop5000} Given (compatibly sized) matrices
$A,B,C$, we have
\[
(A\otimes B)\vecm(C)=\vecm(BCA^{T}).
\]
We can thus view $A\otimes B$ as a vectorized equivalent of the linear
operation $C\mapsto BCA^{T}$. We are tempted to introduce a parallel
notation $(A\otimes B)[C]=BCA^{T}$ for the ``non-vectorized'' version
of this operation, although this notation is not standard.

One possible mnemonic for this identity is that
the $B$ is just to the left of the $C$ while the $A$ ``circles around'' to the right and gets transposed.
\end{proposition}

Where does this identity come from? We can break it into simpler pieces by
first considering the cases where either $A$ or $B$ is an identity
matrix $\Id$ (of the appropriate size). To start with, suppose that
$A=\Id$, so that $BCA^{T}=BC$. What is $\vecm(BC)$? If we let $\vec{c}_{1},\vec{c}_{2},\ldots$
denote the columns of $C$, then recall that $BC$ simply multiples
$B$ on the left with each of the columns of $C$: 
\[
BC=B\left(\begin{array}{ccc}
\vec{c}_{1} & \vec{c}_{2} & \cdots\end{array}\right)=\left(\begin{array}{ccc}
B\vec{c}_{1} & B\vec{c}_{2} & \cdots\end{array}\right)\Longrightarrow\vecm(BC)=\left(\begin{array}{c}
B\vec{c}_{1}\\
B\vec{c}_{2}\\
\vdots
\end{array}\right).
\]
Now, how can we get this $\vecm(BC)$ vector as something multiplying
$\vecm C$? It should be immediately apparent that 
\[
\vecm(BC)=\left(\begin{array}{c}
B\vec{c}_{1}\\
B\vec{c}_{2}\\
\vdots
\end{array}\right)=\underbrace{\left(\begin{array}{ccc}
B\\
 & B\\
 &  & \ddots
\end{array}\right)}_{\Id\otimes B}\underbrace{\left(\begin{array}{c}
\vec{c}_{1}\\
\vec{c}_{2}\\
\vdots
\end{array}\right)}_{\vecm C},
\]
but this matrix is exactly the Kronecker product $I\otimes B$! Hence,
we have derived that 
\[
(\Id\otimes B)\vecm C=\vecm(BC).
\]
What about the $A^{T}$ term? This is a little trickier, but again
let's simplify to the case where $B=\Id$, in which case $BCA^{T}=CA^{T}$.
To vectorize this, we need to look at the columns of $CA^{T}$. What
is the first column of $CA^{T}$? It is a linear combination of the
columns of $C$ whose coefficients are given by the first column of
$A^{T}$ (=~first row of $A$): 
\[
\text{column 1 of }CA^{T}=\sum_{j}a_{1j}\vec{c}_{j}\:.
\]
Similarly for column~2, etc, and we then ``stack'' these columns
to get $\vecm(CA^{T})$. But this is exactly the formula for multipling
a matrix $A$ by a vector, if the ``elements'' of the vector were
the columns $\vec{c}_{j}$. Written out explicitly, this becomes:
\[
\vecm(CA^{T})=\left(\begin{array}{c}
\sum_{j}a_{1j}\vec{c}_{j}\\
\sum_{j}a_{2j}\vec{c}_{j}\\
\vdots
\end{array}\right)=\underbrace{\left(\begin{array}{ccc}
a_{11}\Id & a_{12}\Id & \cdots\\
a_{21}\Id & a_{22}\Id & \cdots\\
\vdots & \vdots & \ddots
\end{array}\right)}_{A\otimes\Id}\underbrace{\left(\begin{array}{c}
\vec{c}_{1}\\
\vec{c}_{2}\\
\vdots
\end{array}\right)}_{\vecm C},
\]
and hence we have derived
\[
(A\otimes\Id)\vecm C=\vecm(CA^{T}).
\]
The full identity $(A\otimes B)\vecm(C)=\vecm(BCA^{T})$ can then
be obtained by straightforwardly combining these two derivations: replace $CA^T$ with $BCA^T$ in the second derivation, which replaces $\vec{c}_j$ with $B\vec{c}_j$ and hence $\Id$ with $B$.

\subsubsection{The Jacobian in Kronecker-product notation}

So now we want to use Prop.~\ref{prop5000} to calculate the
Jacobian of $f(A)=A^{2}$ in terms of the Kronecker product. Let $\d A$
be our $C$ in Prop.~\ref{prop5000}. We can now immediately see that
\[
\vecm(A\,dA+\d A\,A)=\underbrace{(\Id\otimes A+A^{T}\otimes\Id)}_{\mbox{Jacobian }\tilde{f}'(\vecm A)}\vecm(\d A) \, ,
\]
where $\Id$ is the identity matrix of the same size as $A$. We can
also write this in our ``non-vectorized'' linear-operator notation:
\[
A\,\d A+\d A\,A=(\Id\otimes A+A^{T}\otimes\Id)[\d A] \, .
\]
In the $2\times2$ example, these Kronecker products can be computed
explicitly: 
\begin{align*}
\underbrace{\begin{pmatrix}1 & \\
 & 1
\end{pmatrix}}_{\Id}\otimes\underbrace{\begin{pmatrix}p & r\\
q & s
\end{pmatrix}}_{A}+\underbrace{\begin{pmatrix}p & q\\
r & s
\end{pmatrix}}_{A^{T}}\otimes\underbrace{\begin{pmatrix}1 & \\
 & 1
\end{pmatrix}}_{\Id} & =\underbrace{\left(\begin{array}{cccc}
p & r &  & \\
q & s &  & \\
 &  & p & r\\
 &  & q & s
\end{array}\right)}_{\Id\otimes A}+\underbrace{\left(\begin{array}{cccc}
p &  & q & \\
 & p &  & q\\
r &  & s & \\
 & r &  & s
\end{array}\right)}_{A^{T}\otimes\Id}\\
 & =\left(\begin{array}{cccc}
2p & r & q & 0\\
q & p+s & 0 & q\\
r & 0 & p+s & r\\
0 & r & q & 2s
\end{array}\right)=\tilde{f}'\,,
\end{align*}
which exactly matches our laboriously computed Jacobian $\tilde{f}'$
from earlier!

\begin{example}For the matrix-cube function $A^{3}$, where $A$
is an $m\times m$ square matrix, compute the $m^{2}\times m^{2}$
Jacobian of the vectorized function $\vecm(A^{3})$.\end{example}
Let's use the same trick for the matrix-cube function. Sure, we could
laboriously compute the Jacobian via element-by-element partial derivatives
(which is done nicely by symbolic computing in the notebook), but
it's much easier and more elegant to use Kronecker products. Recall
that our ``non-vectorized'' matrix-calculus derivative is the linear
operator: 
\[
(A^{3})'[dA]=dA\,A^{2}+A\,dA\,A+A^{2}\,dA,
\]
which now vectorizes by three applications of the Kronecker
identity:
\[
\vecm(dA\,A^{2}+A\,dA\,A+A^{2}\,dA)=\underbrace{\left((A^{2})^{T}\otimes\Id+A^{T}\otimes A+\Id\otimes A^{2}\right)}_{\text{vectorized Jacobian}}\vecm(\d X)\,.
\]
You could go on to find the Jacobians of $A^{4}$, $A^{5}$, and so
on, or any linear combination of matrix powers. Indeed, you could
imagine applying a similar process to the Taylor series of any (analytic)
matrix function $f(A)$, but it starts to become awkward. Later on
(and in homework), we will discuss more elegant ways to differentiate
other matrix functions, not as vectorized Jacobians but as linear
operators on matrices.

\subsubsection{The computational cost of Kronecker products}

One must be cautious about using Kronecker products as a \emph{computational}
tool, rather than as more of a \emph{conceptual} tool, because they
can easily cause the computational cost of matrix problems to explode
far beyond what is necessary.

Suppose that $A$, $B$, and $C$ are all $m\times m$ matrices. The
cost of multiplying two $m\times m$ matrices (by the usual methods)
scales proportional to $\sim m^{3}$, what the computer scientists
call $\Theta(m^{3})$ ``complexity.'' Hence, the cost of the linear
operation $C\mapsto BCA^{T}$ scales as $\sim m^{3}$ (two $m\times m$ multiplications).
However, if we instead compute the \emph{same answer} via $\vecm(BCA^{T})=(A\otimes B)\vecm C$,
then we must:
\begin{enumerate}
\item Form the $m^{2}\times m^{2}$ matrix $A\otimes B$. This requires
$m^{4}$ multiplications (all entries of $A$ by all entries of $B$),
and $\sim m^{4}$ memory storage. (Compare to $\sim m^{2}$ memory
to store $A$ or $B$. If $m$ is 1000, this is a \emph{million} times
more storage, terabytes instead of megabytes!)
\item Multiply $A\otimes B$ by the vector $\vecm C$ of $m^{2}$ entries.
Multiplying an $n\times n$ matrix by a vector requires $\sim n^{2}$
operations, and here $n=m^{2}$, so this is again $\sim m^{4}$ arithmetic
operations.
\end{enumerate}
So, instead of $\sim m^{3}$ operations and $\sim m^{2}$ storage
to compute $BCA^{T}$, using $(A\otimes B)\vecm C$ requires $\sim m^{4}$
operations and $\sim m^{4}$ storage, vastly worse! Essentially, this
is because $A\otimes B$ has a lot of structure that we are not exploiting
(it is a \emph{very special} $m^{2}\times m^{2}$ matrix). 

There are many examples of this nature. Another famous one involves
solving the linear \emph{matrix} equations 
\[
AX+XB=C
\]
for an unknown matrix $X$, given $A,B,C$, where all of these are
$m\times m$ matrices. This is called a ``Sylvester equation.'' These
are \emph{linear }equations in our unknown $X$, and we can convert them
to an ordinary system of $m^{2}$ linear equations by Kronecker products:
\[
\vecm(AX+XB)=(\Id\otimes A+B^{T}\otimes\Id)\vecm X=\vecm C,
\]
which you can then solve for the $m^{2}$ unknowns $\vecm X$ using
Gaussian elimination. But the cost of solving an $m^{2}\times m^{2}$
system of equations by Gaussian elimination is $\sim (m^2)^3 = m^{6}$.
It turns out, however, that there are clever algorithms to solve $AX+XB=C$
in only $\sim m^{3}$ operations (with $\sim m^{2}$ memory)---for $m=1000$, this saves a factor of $\sim m^3 = {10}^9$ (a \emph{billion}) in computational effort.

(Kronecker products can be a more practical computational tool for \emph{sparse} matrices: matrices that are mostly zero, e.g.~having only a few nonzero entries per row.  That's because the Kronecker product of two sparse matrices is also sparse, avoiding the huge storage requirements for Kronecker products of non-sparse ``dense'' matrices. This can be a convenient way to assemble large sparse systems of equations for things like multidimensional PDEs.)

\pagebreak

\section{Finite-Difference Approximations}
\label{sec:finitedifference}
In this section, we will be referring to this \href{https://github.com/mitmath/matrixcalc/blob/main/notes/Finite%20difference%20checks.ipynb}{Julia notebook} for calculations that are not included here.

\subsection{Why compute derivatives approximately instead of exactly?} \label{sec:approximation}

Working out derivatives by hand is a notoriously error-prone procedure for complicated functions.  Even if every individual step is straightforward, there are so many opportunities to make a mistake, either in the derivation or in its implementation on a computer.   Whenever you implement a derivatives, you should \textbf{always double-check} for mistakes by comparing it to an independent calculation.  The simplest such check is a \emph{finite-difference approximation}, in which we \emph{estimate} the derivative(s) by comparing $f(x)$ and $f(x + \delta x)$ for one or more ``finite'' (non-infinitesimal) perturbations $\delta x$.

There are many finite-difference techniques at varying levels of sophistication, as we will discuss below.  They all incur an intrinsic \textbf{truncation error} due to the fact that $\delta x$ is not infinitesimal.  (And we will also see that you can't make $\delta x$ too small, either, or \emph{roundoff errors} start exploding!)   Moreover, finite differences become expensive for higher-dimensional $x$ (in which you need a separate finite difference for each input dimension to compute the full Jacobian).   This makes them an approach of last resort for computing derivatives accurately.  On the other hand, they are the \emph{first} method you generally employ in order to \emph{check} derivatives: if you have a bug in your analytical derivative calculation, usually the answer is completely wrong, so even a crude finite-difference approximation for a single small $\delta x$ (chosen at random in higher dimensions) will typically reveal the problem.

Another alternative is \textbf{automatic differentiation} (AD), software/compilers perform \emph{analytical} derivatives for you. This is extremely reliable and, with modern AD software, can be very efficient. Unfortunately, there is still lots of code, e.g. code calling external libraries in other languages, that AD tools can't comprehend. And there are other cases where AD is inefficient, typically because it misses some mathematical structure of the problem.  Even in such cases, you can often fix AD by defining the derivative of one small piece of your program by hand,\footnote{In some Julia AD software, this is done with by defining a  \href{https://github.com/JuliaDiff/ChainRulesCore.jl}{``ChainRule''}, and in Python autograd/JAX it is done by defining a custom ``vJp'' (row-vector—Jacobian product) and/or ``Jvp'' (Jacobian–vector product).} which is much easier than differentiating the whole thing.  In such cases, you still will typically want a finite-difference check to ensure that you have not made a mistake.

It turns out that finite-difference approximations are a surprisingly complicated subject, with rich connections to many areas of numerical analysis; in this lecture we will just scratch the surface.

\subsection{Finite-Difference Approximations: Easy Version}

The simplest way to check a derivative is to recall that the definition of a differential:
\begin{displaymath}
	\d f = f(x+\d x) - f(x) = f'(x) \d x
\end{displaymath}
came from dropping higher-order terms from a small but finite difference:
\begin{displaymath}
	\delta f = f(x+\delta x) - f(x) = f'(x) \delta x + o(\Vert \delta x \Vert) \, .
\end{displaymath}
So, we can just compare the \textbf{finite difference} $\boxed{f(x+\delta x) - f(x)}$ to our \textbf{(directional) derivative operator} $f'(x) \delta x$ (i.e. the derivative in the direction $\delta x$). $f(x+\delta x) - f(x)$ is also called a \textbf{forward difference} approximation.  The antonym of a forward difference is a \textbf{backward difference} approximation $f(x) - f(x - \delta x) \approx f'(x) \delta x$. If you just want to compute a derivative, there is not much practical distinction between forward and backward differences.  The distinction becomes more important when discretizing (approximating) differential equations. We'll look at other possibilities below.

\begin{remark}
Note that this definition of forward and backward difference is \textbf{not} the same as forward- and backward-mode differentiation---these are \textbf{unrelated} concepts.
\end{remark}

If $x$ is a scalar, we can also divide both sides by $\delta x$ to get an approximation for $f'(x)$ instead of for $df$:
\begin{displaymath}
	f'(x) \approx \frac{f(x+\delta x) - f(x)}{\delta x} + \text{(higher-order corrections)} \, .
\end{displaymath}
This is a more common way to write the forward-difference approximation, but it only works for scalar $x$, whereas in this class we want to think of $x$ as perhaps belonging to some other vector space. 

Finite-difference approximations come in many forms, but they are generally a \textbf{last resort} in cases where it's too much effort to work out an analytical derivative and AD fails.  But they are also useful to \textbf{check} your analytical derivatives and to quickly \textbf{explore}.

\subsection{Example: Matrix squaring}

Let's try the finite-difference approximation for the square function $f(A) = A^2$, where here $A$ is a square matrix in $\R^{m,m}$. By hand, we obtain the product rule 
\[
\d f = A \,\d A + \d A \, A,
\]
i.e. $f'(A)$ is the \textbf{linear operator} $\boxed{f'(A)[\delta A] = A \,\delta A + \delta A \,A.}$ This is \textit{not equal to} $2A \,\delta A$ because \textit{in general} $A$ and $\delta A$ do not commute. So let's check this difference against a finite difference.  We'll try it for a \textit{random} input A and a \textit{random small} perturbation $\delta A$.

Using a random matrix $A$, let $\d A = A \cdot 10^{-8}$. Then, you can compare $f(A+\d A) - f(A)$ to $A \,\d A + \d A \, A$. If the matrix you chose was really random, you would get that the approximation minus the exact equality from the product rule has entries with order of magnitude around $10^{-16}!$ However, compared to $2A \d A$, you'd obtain entries of order $10^{-8}$.

To be more quantitative, we might compute that "norm" $\lVert \text{approx} - \text{exact} \rVert$ which we want to be small. But small \textbf{compared to what}? The natural answer is \textbf{small compared to the correct answer.} This is called the \href{https://en.wikipedia.org/wiki/Approximation_error}{relative error} (or "fractional error") and is computed via 
\[
\text{relative error} = \frac{\lVert \text{approx} - \text{exact} \rVert}{\lVert \text{exact} \rVert}.
\]
Here, $\lVert \cdot \rVert$ is a \href{https://en.wikipedia.org/wiki/Norm_(mathematics)}{norm}, like the length of a vector. This allows us to understand the size of the error in the finite difference approximation, i.e. it allows us to answer how accurate this approximation is (recall Sec.~\ref{sec:approximation}).

So, as above, you can compute that the relative error between the approximation and the exact answer is about $10^{-8}$, where as the relative error between $2A \d A$ and the exact answer is about $10^{0}$. This shows that our exact answer is likely correct! Getting a good match up between a random input and small displacement isn't a proof of correctness of course, but it is always a good thing to check. This kind of randomized comparison will almost always \textbf{catch major bugs} where you have calculated the symbolic derivative incorrectly, like in our $2A \d A$ example.

\begin{definition}  
Note that the norm of a matrix that we are using, computed by \texttt{norm(A)} in Julia, is just the direct analogue of the familiar Euclidean norm for the case of vectors. It is simply the square root of the sum of the matrix entries squared:
\[
\lVert A \rVert := \sqrt{\sum_{i,j} |A_{ij}|^2} = \sqrt{\tr (A^T A)} \, .
\]
This is called the \href{https://mathworld.wolfram.com/FrobeniusNorm.html}{Frobenius norm}.
\end{definition}

\subsection{Accuracy of Finite Differences}

Now how accurate is our finite-difference approximation above? How should we choose the size of $\delta x$?

Let's again consider the example $f(A) = A^2$, and plot the relative error as a function of $\lVert \delta A \rVert$. This plot will be done \textit{logarithmically} (on a log--log scale) so that we can see power-law relationships as straight lines.
\begin{figure}
    \centering
    \includegraphics[width=0.8\textwidth]{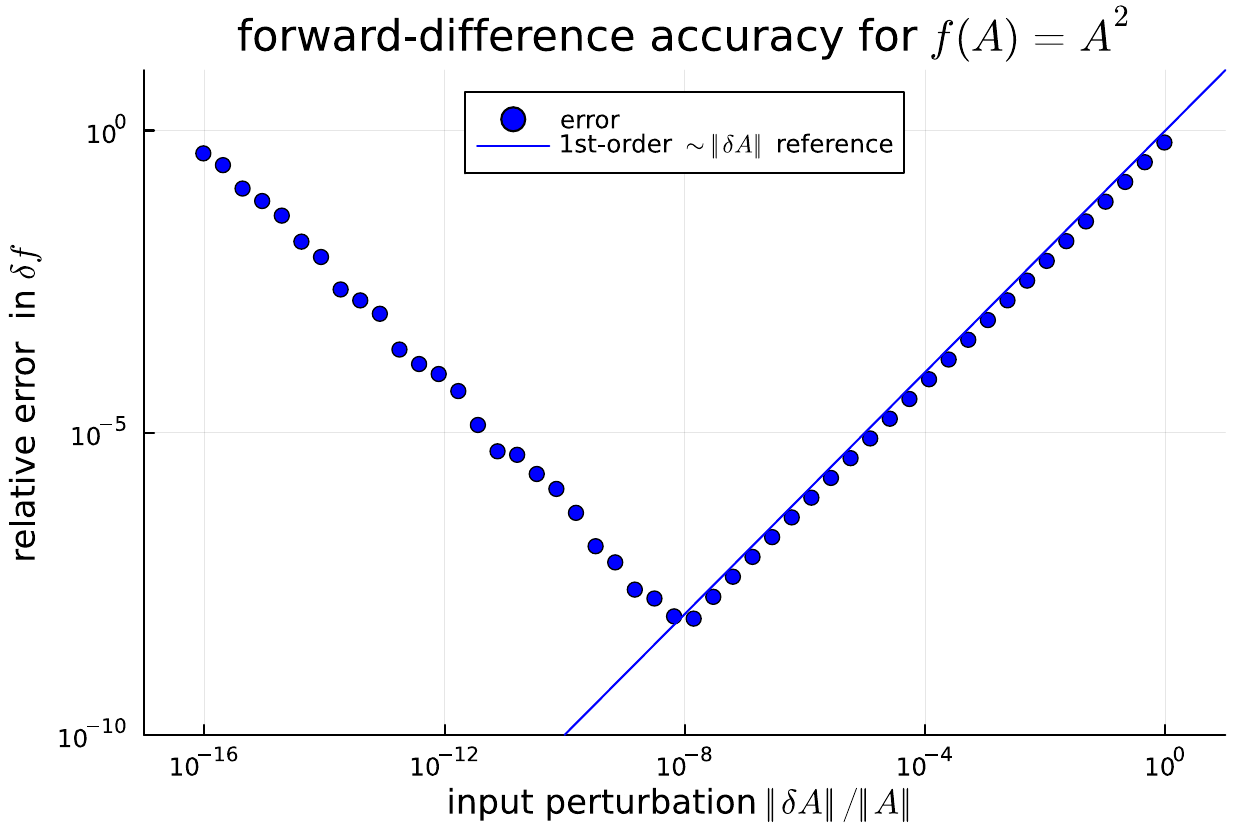}
    \caption{Forward-difference accuracy for $f(A) = A^2$, showing the relative error in $\delta f = f(A + \delta A) - f(A)$ versus the  linearization $f'(A) \delta A$, as a function of the magnitude $\Vert \delta A\Vert$.  $A$ is a $4 \times 4$ matrix with unit-variance Gaussian random entries, and $\delta A$ is similarly a unit-variance Gaussian random perturbation scaled by a factor $s$ ranging from $1$ to $10^{-16}$.}
\end{figure}

We notice two main features as we decrease $\delta A$:
\begin{enumerate}
    \item The relative error at first decreases linearly with $\lVert \delta A\rVert$. This is called \textbf{first-order accuracy}. Why?
    \item When $\delta A$ gets too small, the error increases. Why?
\end{enumerate}

\subsection{Order of accuracy}

The \textbf{truncation error} is the inaccuracy arising from the fact that the input perturbation $\delta x$ is not infinitesimal: we are computing a difference, not a derivative.
If the truncation error in the derivative scales proportional $\Vert \delta x \Vert^n$, we call the approximation \textbf{n-th order accurate}.  For forward differences, here, the order is \textbf{n=1}.  Why?

For any $f(x)$ with a nonzero second derivative (think of the Taylor series), we have
\begin{displaymath}
	f(x + \delta x) = f(x) + f'(x) \delta x + (\text{terms proportional to }\Vert \delta x \Vert^2) + \underbrace{o(\Vert \delta x \Vert^2)}_\text{i.e. higher-order terms}
\end{displaymath}
That is, the terms we \textit{dropped} in our forward-difference approximations are proportional to $\Vert \delta x\Vert^2$.  But that means that the \textbf{relative error is linear}:
\begin{align*}
	\text{relative error} &= \frac{\Vert f(x+\delta x) - f(x) - f'(x) \delta x \Vert}{\Vert f'(x) \delta x \Vert} \\
	&= \frac{(\text{terms proportional to }\Vert \delta x \Vert^2) + o(\Vert \delta x \Vert^2)}{\text{proportional to }\Vert \delta x \Vert} = (\text{terms proportional to }\Vert \delta x \Vert) + o(\Vert \delta x \Vert)
\end{align*}
This is \textbf{first-order accuracy}.
Truncation error in a finite-difference approximation is the \textbf{inherent} error in the formula for \textbf{non-infinitesimal} $\delta x$.  Does that mean we should just make $\delta x$ as small as we possibly can?

\subsection{Roundoff error}

The reason why the error \textit{increased} for very small $\delta A$ was due to \textbf{roundoff errors}.  The computer only stores a \textbf{finite number of significant digits} (about 15 decimal digits) for each real number and rounds off the rest on each operation — this is called  \href{https://en.wikipedia.org/wiki/Floating-point_arithmetic}{floating-point arithmetic}. If $\delta x$ is too small, then the difference $f(x+\delta x) - f(x)$ gets rounded off to zero (some or all of the \textit{significant digits cancel}).  This is called  \href{https://en.wikipedia.org/wiki/Catastrophic_cancellation}{catastrophic cancellation}.

Floating-point arithmetic is much like scientific notation ${*}.{*}{*}{*}{*}{*} \times 10^e$: a finite-precision coefficient ${*}.{*}{*}{*}{*}{*}$ scaled by a power of~$10$ (or, on a computer, a power of~$2$). The number of digits in the coefficient (the ``significant digits'') is the ``precision,'' which in the usual 64-bit floating-point arithmetic is charactized by a quantity $\epsilon = 2^{-52} \approx 2.22 \times 10^{-16}$, called the \href{https://en.wikipedia.org/wiki/Machine_epsilon}{machine epsilon}.  When an arbitrary real number $y \in \mathbb{R}$ is rounded to the closest floating-point value $\tilde{y}$, the roundoff error is bounded by $|\tilde{y} - y| \le \epsilon |y|$.  Equivalently, the computer keeps only about 15--16 $\approx -\log_{10} \epsilon$ decimal digits, or really $53 = 1 - \log_2 \epsilon$ \emph{binary} digits, for each number. 

In our finite-difference example, for $\Vert \delta A \Vert / \Vert A \Vert$ of roughly $10^{-8} \approx \sqrt{\epsilon} \lVert A\rVert$ or larger, the approximation for $f'(A)$ is dominated by the truncation error, but if we go smaller than that the relative error starts increasing due to roundoff. Experience has shown that $\Vert \delta x \Vert \approx \sqrt{\epsilon} \Vert x \Vert $ is often a good rule of thumb---about half the significant digits is the most that is reasonably safe to rely on---but the precise crossover point of minimum error depends on the function $f$ and the finite-difference method. But, like all rules of thumb, this may not always be completely reliable.

\subsection{Other finite-difference methods}

There are more sophisticated finite-difference methods, such as Richardson extrapolation, which consider a sequence of progressively smaller $\delta x$ values in order to adaptively determine the best possible estimate for $f'$ (\emph{extrapolating} to $\delta x \to 0$ using progressively higher degree polynomials).   One can also use higher-order difference formulas than the simple forward-difference method here, so that the truncation error decreases faster than than linearly with $\delta x$.  The most famous higher-order formula is the ``centered difference'' $f'(x) \delta x \approx [f(x + \delta x) - f(x - \delta x)]/2$, which has \emph{second}-order accuracy (relative truncation error proportional to $\Vert \delta x\Vert^2$).  

Higher-dimensional inputs  $x$ pose a fundamental computational challenge for finite-difference techniques, because if you want to know what happens for every possible direction $\delta x$ then you need many finite differences: one for each dimension of $\delta x$. For example, suppose $x \in \mathbb{R}^n$ and $f(x) \in \mathbb{R}$, so that you are computing $\nabla f \in \mathbb{R}^n$; if you want to know the whole gradient, you need $n$ \emph{separate} finite differences.
The net result is that finite differences in higher dimensions are expensive, quickly becoming impractical for high-dimensional optimization (e.g. neural networks) where $n$  might be huge.  On the other hand, if you are just using finite differences as a check for bugs in your code, it is usually sufficient to compare $f(x+\delta x) - f(x)$ to $f'(x)[\delta x]$ in a few random directions, i.e.~for a few random small $\delta x$.

\pagebreak
\section{Derivatives in General Vector Spaces}
\label{sec:generalvectorspaces}
Matrix calculus
requires us  to generalize concepts of derivative and gradient further, to functions whose inputs and/or outputs are not simply scalars or column vectors.  To achieve this, we extend the notion of the ordinary vector \textbf{dot product} and
ordinary Euclidean vector ``length'' to
general \textbf{inner products} and \textbf{norms} on 
\textbf{vector spaces}.
Our first example will consider familiar matrices
from this point of view.

Recall from linear algebra that we can call any set $V$ a ``vector space'' if its elements can be added/subtracted $x \pm y$ and multiplied by scalars $\alpha x$ (subject to some basic arithmetic axioms, e.g.~the  distributive law).  For example, the set of $m \times n$ matrices themselves form a vector space, or even the set of continuous functions $u(x)$ (mapping $\mathbb{R} \to \mathbb{R}$)---the key fact is that we can add/subtract/scale them and get elements of the same set.  It turns out to be extraordinarily useful to extend differentiation to such spaces, e.g.~for functions that map matrices to matrices or functions to numbers.  Doing so crucially relies on our input/output vector spaces $V$ having a \textbf{norm} and, ideally, an \textbf{inner product}.

\subsection{A Simple Matrix Dot Product and Norm}


Recall that for \textit{scalar-valued} functions $f(x) \in \R$ with \textit{vector inputs} $x\in \R^n$ (i.e. $n$-component ``column vectors") we have that 
    \[
    \d f = f(x + \d x) - f(x) = f'(x) [\d x] \in \R.
    \]
    Therefore, $f'(x)$ is a linear operator taking in the vector $\d x$ in and giving a scalar value out. Another way to view this is that $f'(x)$ is the row vector\footnote{The concept of a ``row vector'' can be formalized as something called a ``covector,'' a ``dual vector,'' or an element of a ``\href{https://en.wikipedia.org/wiki/Dual_space}{dual space},'' not to be confused with the \emph{dual numbers} used in automatic differentiation (Sec.~\ref{sec:AD}).}  $(\nabla f)^T$. Under this viewpoint, it follows that $\d f$ is the dot product (or ``inner product''):
    \[
    \d f = \nabla f \cdot \d x
    \]

We can generalize this to any vector space $V$ with inner products! Given $x\in V$, and a scalar-valued function $f$, we obtain the linear operator $f'(x) [\d x] \in \R$, called a ``linear form.'' In order to define the gradient $\nabla f$, we need an inner product for $V$, the vector-space generalization of the familiar dot product!

Given $x,y \in V$, the inner product $\langle x, y \rangle$  is a map ($\cdot$) such that $\langle x, y \rangle \in \R$. This is also commonly denoted $x \cdot y$ or $\langle x \mid y \rangle$. More technically, an inner product is a map  that is 
\begin{enumerate}
    \item \textbf{Symmetric}: i.e. $\langle x, y \rangle = \langle y, x \rangle$ (or conjugate-symmetric,\footnote{Some authors distinguish the ``dot product'' from an ``inner product'' for complex vector spaces, saying that a dot product has no complex conjugation $x \cdot y = y \cdot x$ (in which case $x \cdot x$ need not be real and does not equal $\Vert x \Vert^2$), whereas the inner product must be conjugate-symmetric, via $\langle x, y \rangle = \bar{x} \cdot y$.  Another source of confusion for complex vector spaces is that some fields of mathematics define $\langle x, y \rangle = x \cdot \bar{y}$, i.e.~they conjugate the \emph{right} argument instead of the left (so that it is linear in the left argument and conjugate-linear in the right argument).  Aren't you glad we're sticking with real numbers?} $\langle x, y \rangle = \overline{\langle y, x \rangle}$, if we were using complex numbers), 
    \item \textbf{Linear}: i.e. $\langle x, \alpha y + \beta z\rangle = \alpha \langle x, y \rangle + \beta \langle x, z \rangle$, and 
    \item \textbf{Non-negative}: i.e. $\langle x, x \rangle := \lVert x \rVert^2 \geq 0$, and $=0$ if and only if $x = 0$.
\end{enumerate}
Note that the combination of the first two properties means that it must also be linear in the left vector (or conjugate-linear, if we were using complex numbers).  Another useful consequence of these three properties, which is a bit trickier to derive, is the \emph{Cauchy--Schwarz inequality} $|\langle x, y \rangle| \le \Vert x \Vert \, \Vert y \Vert$.

\begin{definition}[Hilbert Space]
A (complete) vector space with an inner product is called a \textit{Hilbert space}.  (The technical requirement of ``completeness'' essentially means that you can take limits in the space, and is important for rigorous proofs.\footnote{Completeness means that any Cauchy sequence of points in the vector space---any sequence of points that gets closer and closer together---has a limit lying within the vector space.  This criterion usually holds in practice for vector spaces over real or complex scalars, but can get trickier when talking about vector spaces of functions, since e.g.~the limit of a sequence of continuous functions can be a discontinuous function.})
\end{definition}

Once we have a Hilbert space, we can define the gradient for scalar-valued functions. Given $x\in V$ a Hilbert space, and $f(x)$ scalar, then we have the linear form $f'(x) [\d x] \in \R$. Then, under these assumptions, there is a theorem known as the ``Riesz representation theorem'' stating that \emph{any} linear form (including $f'$) must be an inner product with \emph{something}: 
\[
f'(x) [\d x] = \big\langle \underbrace{\text{(some vector)}}_{\text{gradient } \nabla f\bigr|_x} , \d x \big\rangle = \d f.
\]
That is, the gradient $\nabla f$ is \emph{defined} as the thing you take the inner product of $\d x$ with to get $\d f$.
Note that $\nabla f$ always has the ``same shape'' as $x$.

The first few examples we look at involve the usual Hilbert space $V = \R^n$ with different inner products.

\begin{example}
    Given $V = \R^n$ with $n$-column vectors, we have the familiar Euclidean dot product $\langle x, y \rangle = x^T y$. This leads to the usual $\nabla f$.
\end{example}

\begin{example}
    We can have different inner products on $\R^n$. For instance, 
    \[
    \langle x, y\rangle_W = w_1 x_1 y_1 + w_2 x_2 y_2 + \dots w_n x_n y_n = x^T \underbrace{\begin{pmatrix}
        w_1 & & \\
         & \ddots & \\
         & & w_n
    \end{pmatrix}}_{W} y
    \]
    for weights $w_1,\dots, w_n >0$. 
    
    More generally we can define a weighted dot product $\langle x, y\rangle_W= x^T W y$ for any symmetric-positive-definite matrix $W$ ($W = W^T$ and $W$ is positive definite, which is sufficient for this to be a valid inner product).

    If we change the definition of the inner product, then we change the definition of the gradient!  For example, with $f(x) = x^T A x$ we previously found that $\d f = x^T (A + A^T) \d x$.  With the ordinary Euclidean inner product, this gave a gradient $\nabla f = (A+A^T)x$.  However, if we use the weighted inner product $x^T W y$, then we would obtain a different ``gradient'' $\nabla^{(W)} f = W^{-1} (A+A^T)x$ so that $\d f = \langle \nabla^{(W)}  f , \d x \rangle$.
    
    In these notes, we will employ the Euclidean inner product for $x \in \mathbb{R}^n$, and hence the usual $\nabla f$, unless noted otherwise.  However, weighted inner products are useful in lots of cases, especially when the components of $x$ have different scales/units.
\end{example}

We can also consider the space of $m\times n$ matrices $V = \R^{m \times n}$. There, is of course, a vector-space isomorphism from $V \ni A \to \mathrm{vec}(A) \in \R^{mn}$. Thus, in this space we have the analogue of the familiar (``Frobenius") Euclidean inner product, which is convenient to rewrite in terms of matrix operations via the trace: 
\begin{definition}[Frobenius inner product]
The \textbf{Frobenius inner product} of two $m \times n$ matrices $A$ and $B$ is:
\[
\langle A, B \rangle_F = \sum_{ij} A_{ij} B_{ij} = \mathrm{vec}(A)^T \mathrm{vec}(B) = \tr(A^T B) \, .
\]
Given this inner product, we also have the corresponding \textbf{Frobenius norm}: $$\lVert A \rVert_F = \sqrt{\langle A,A \rangle_F} = \sqrt{\tr(A^TA)} = \lVert \mathrm{vec} A\rVert = \sqrt{\sum_{i,j} |A_{ij}|^2} \, .$$ 
Using this, we can now define the gradient of scalar functions with \textit{matrix inputs}!  This will be our default matrix inner product (hence defining our default matrix gradient) in these notes (sometimes dropping the $F$ subscript).
\end{definition}

\begin{example}
    Consider the function 
    \[
    f(A) = \lVert A \rVert_F = \sqrt{\tr(A^T A)}.
    \] What is $\d f$?
\end{example}
Firstly, by the familiar scalar-differentiation chain and power rules we have that 
\[
\d f = \frac{1}{2 \sqrt{\tr(A^T A)}} \d (\tr A^T A).
\]
Then, note that (by linearity of the trace)
\[
\d( \tr B) = \tr(B+\d B) - \tr(B) = \tr(B) + \tr(\d B) - \tr(B) = \tr(\d B).
\]
Hence, 
\begin{align*}
    \d f &= \frac{1}{2\lVert A\rVert_F} \tr(\d (A^T A)) \\
    &= \frac{1}{2\lVert A\rVert_F} \tr( \d A^T\, A + A^T\, \d A) \\
    &= \frac{1}{2 \lVert A\rVert_F} (\tr(\d A^T\, A) + \tr(A^T\, \d A)) \\
    &= \frac{1}{\lVert A\rVert_F} \tr(A^T \, \d A) = \big\langle \frac{A}{\lVert A\rVert_F} , \d A \big\rangle.
\end{align*}
Here, we used the fact that $\tr B = \tr B^T$, and in the last step we connected $\d f$ with a Frobenius inner product. In other words, 
\[
\nabla f = \nabla \lVert A \rVert_F = \frac{A}{\lVert A \rVert_F}.
\]
Note that one obtains exactly the same result for column vectors~$x$, i.e.~$\nabla \Vert x\Vert = x/\Vert x \Vert$ (and in fact this is equivalent via $x = \vecm A$).

Let's consider another simple example:

\begin{example}
Fix some constant $x \in \R^m$, $y\in \R^n$, and consider the function $f:\R^{m\times n} \to \R$ given by
\[
f(A) = x^T A y.
\]
What is $\nabla f$?
\end{example}
We have that 
\begin{align*}
    \d f &= x^T \,\d A \,y \\
    &= \tr( x^T \,\d A\, y) \\
    &= \tr(y x^T \, \d A) \\
    &= \big\langle \underbrace{x y^T}_{\nabla f} , \, \d A \big\rangle.
\end{align*}

More generally, for any scalar-valued function $f(A)$, from the definition of Frobenius inner product it follows that:
$$
\d f = f(A+\d A)-f(A) = \langle \nabla f , \, \d A \rangle = \sum_{i,j} (\nabla f)_{i,j} \, \d A_{i,j} \, ,
$$
and hence the components of the gradient are exactly the elementwise derivatives
$$
(\nabla f)_{i,j} = \frac{\partial f}{\partial A_{i,j}} \, ,
$$
similar to the component-wise definition of the gradient vector from multivariable calculus!  But for non-trivial matrix-input functions $f(A)$ it can be extremely awkward to take the derivative with respect to each entry of $A$ individually.
Using the ``holistic'' matrix inner-product definition, we will soon be able to compute even more complicated matrix-valued gradients, including  $\nabla (\det A)$!

\subsection{Derivatives, Norms, and Banach spaces}
\label{sec:banach}

We have been using the term ``norm'' throughout this class, but what technically is a norm?  Of course, there are familiar examples such as the Euclidean (``$\ell^2$'') norm $\Vert x \Vert = \sqrt{\sum_k x_k^2}$ for $x\in \mathbb{R}^n$, but it is useful to consider how this concept generalizes to other vector spaces.   It turns out, in fact, that norms are crucial to the definition of a derivative!

Given a vector space $V$, a norm $\lVert \cdot \rVert$ on $V$ is a map $\lVert \cdot \rVert: V\to \R$ satisfying the following three properties:
\begin{enumerate}
    \item \textbf{Non-negative}: i.e. $\lVert v \rVert \geq 0$ and $\lVert v \rVert = 0 \iff v = 0$,
    \item \textbf{Homogeneity}: $\lVert \alpha v \rVert = |\alpha |\lVert v \rVert$ for any $\alpha \in \R$, and 
    \item \textbf{Triangle inequality}: $\lVert u + v\rVert \leq \lVert u \rVert + \lVert v \rVert$.
\end{enumerate}

A vector space that has a norm is called an \textit{normed vector space}.  Often, mathematicians technically want a slightly more precise type of normed vector space with a less obvious name: a \textit{Banach} space.

\begin{definition}[Banach Space]
    A (complete) vector space with a norm is called a \textit{Banach space}.   (As with Hilbert spaces, ``completeness'' is a technical requirement for some types of rigorous analysis, essentially allowing you to take limits.)

\end{definition}

For example, given any inner product $\langle u , v \rangle$, there is a corresponding norm $\lVert u \rVert = \sqrt{\langle u , u \rangle}$.   (Thus, every Hilbert space is also a Banach space.\footnote{Proving the triangle inequality for an arbitrary inner product is not so obvious; one uses a result called the Cauchy--Schwarz inequality.})

To define derivatives, we technically need both the input \textit{and} the output to be Banach spaces. To see this, recall our formalism 
\[
f(x + \delta x) - f(x) = \underbrace{f'(x) [\delta x]}_{\mbox{linear}} \; + \underbrace{o(\delta x)}_{\mbox{smaller}}\, .
\]
To precisely define the sense in which the $o(\delta x)$ terms are ``smaller'' or ``higher-order,'' we need norms. In particular, the ``little-$o$'' notation $o(\delta x)$ denotes any function such that 
\[
\lim_{\delta x\to 0} \frac{\lVert o (\delta x) \rVert}{\lVert \delta x\rVert} = 0 \, ,
\]
i.e.~which goes to zero faster than linearly in $\delta x$.
This requires both the input $\delta x$ and the output (the function) to have norms.  This extension of differentiation to arbitrary normed/Banach spaces is sometimes called the \textbf{Fr{\'e}chet derivative}.

\pagebreak

\section[Nonlinear Root-Finding, Optimization, and Adjoint Differentiation]{Nonlinear Root-Finding, Optimization,\\ and Adjoint Differentiation}
The next part is based on these \href{https://docs.google.com/presentation/d/1U1lB5bhscjbxEuH5FcFwMl5xbHl0qIEkMf5rm0MO8uE/edit#slide=id.p}{slides}. Today, we want to talk about why we are computing derivatives in the first place. In particular, we will drill down on this a little bit and then talk about computation of derivatives.

\subsection{Newton's Method}
\label{sec:newton-roots}

One common application of derivatives is to solve nonlinear equations via linearization. 

\subsubsection{Scalar Functions}

\begin{figure}
    \centering
    \includegraphics[width=0.7\textwidth]{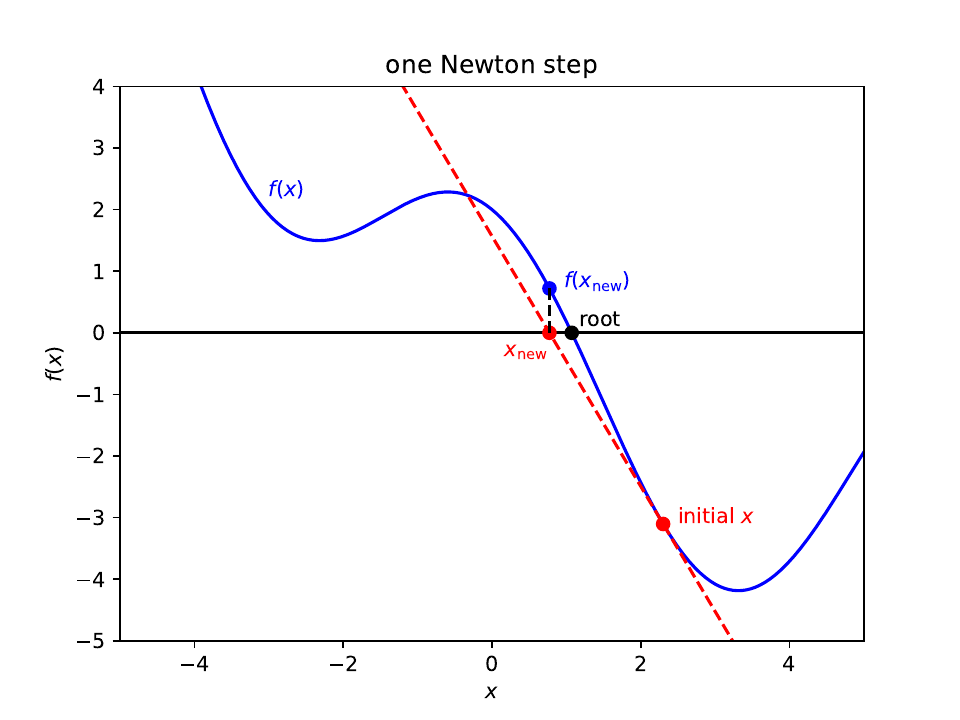}
    \caption{Single step of the scalar Newton's method to solve $f(x)=0$ for an example nonlinear function $f(x) = 2\cos(x) - x + x^2/10$.  Given a starting guess ($x = 2.3$ in this example), we use $f(x)$ and $f'(x)$ to form a linear (affine) approximation of $f$, and then our next step $x_\mathrm{new}$ is the root of this approximation.  As long as the initial guess is not too far from the root, Newton's method converges extremely rapidly to the exact root (black dot).}
    \label{fig:newton-step}
\end{figure}

For instance, suppose we have  a scalar function $f: \R \to \R$ and we wanted to solve $f(x) = 0$ for a root~$x$. Of course, we could solve such an equation explicitly in simple cases, such as when $f$ is linear or quadratic, but if the function is something more arbitrary like $f(x) = x^3 - \sin (\cos x)$ you might not be able to obtain closed-form solutions. However, there is a nice way to obtain the solution approximately to any accuracy you want, as long if you know approximately where the root is. The method we are talking about is known as \emph{Newton's method}, which is really a linear-algebra technique. It takes in the function and a guess for the root, approximates it by a straight line (whose root is easy to find), which is then an approximate root that we can use as a new guess. In particular, the method (depicted in Fig.~\ref{fig:newton-step}) is as follows: 
\begin{itemize}
    \item Linearize $f(x)$ near some $x$ using the approximation 
    \[
    f(x + \delta x) \approx f(x) + f'(x) \delta x,
    \]
    \item solve the linear equation $f(x) + f'(x) \delta x = 0 \implies \delta x = -\frac{f(x)}{f'(x)}$,
    \item and then use this to update the value of $x$ we linearized near---i.e.,~letting the new $x$ be $$x_\text{new} = x - \delta x = x + \frac{f(x)}{f'(x)} \, .$$
\end{itemize}
Once you are close to the root, Newton's method converges amazingly quickly.  As discussed below, it asymptotically \emph{doubles} the number of correct digits on every step!

One may ask what happens when $f'(x)$ is not invertible, for instance here if $f'(x) = 0$. If this happens, then Newton's method may break down! See \href{https://en.wikipedia.org/wiki/Newton%27s_method#Failure_analysis}{here} for examples of when Newton's method breaks down.

\subsubsection{Multidimensional Functions}
We can generalize Newton's method to multidimensional functions! Let $f: \R^n \to \R^n$ be a function which takes in a vector and spits out a vector of the same size~$n$. We can then apply a Newton approach in higher dimensions: 
\begin{itemize}
    \item Linearize $f(x)$ near some $x$ using the first-derivative approximation 
    \[
    f(x + \delta x) \approx f(x) + \underbrace{f'(x)}_\text{Jacobian} \delta x,
    \]
    \item solve the linear equation $f(x) + f'(x) \delta x = 0 \implies \delta x = -\underbrace{f'(x)^{-1}}_\text{inverse Jacobian} f(x)$,
    \item and then use this to update the value of $x$ we linearized near---i.e.,~letting the new $x$ be $$x_\text{new} = x_\text{old} - f'(x)^{-1}f(x)\, .$$
\end{itemize}

That's it! Once we have the Jacobian, we can just solve a linear system on each step. This again converges amazingly fast, doubling the number of digits of accuracy in each step. (This is known as ``quadratic convergence.'') However, there is a caveat: we \textit{need} some starting guess for $x$, and the guess needs to be sufficiently close to the root for the algorithm to make reliable progress. (If you start with an initial $x$ far from a root, Newton's method can fail to converge and/or it can jump around in intricate and surprising ways---google ``Newton fractal'' for some fascinating examples.) This is a widely used and very practical application of Jacobians and derivatives!

\subsection{Optimization}

\subsubsection{Nonlinear Optimization}
\begin{figure}
    \centering
    \includegraphics[width=0.7\textwidth]{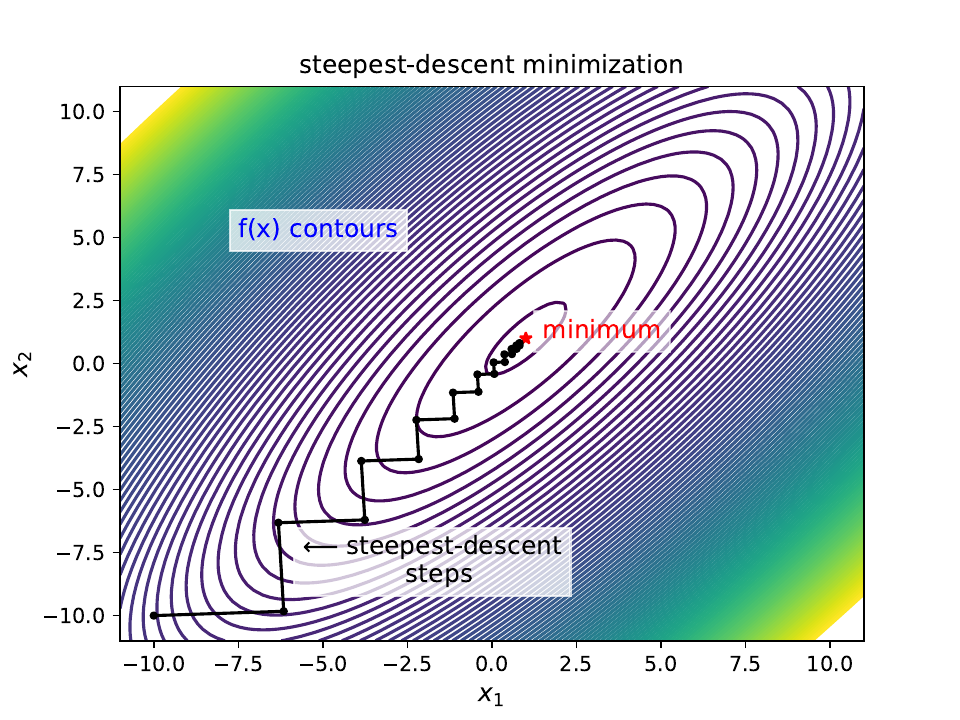}
    \caption{    
    A \emph{steepest-descent algorithm} minimizes a function $f(x)$ by taking successive ``downhill'' steps in the direction $-\nabla f$.  (In the example shown here, we are minimizing a quadratic function in two dimensions $x \in \mathbb{R}^2$, performing an exact 1d minimization in the downhill direction for each step.)  Steepest-descent algorithms can sometimes ``zig-zag'' along narrow valleys, slowing convergence (which can be counteracted in more sophisticated algorithms by ``momentum'' terms, second-derivative information, and so on).}
    \label{fig:steepest-descent}
\end{figure}

A perhaps even more famous application of large-scale differentiation is to nonlinear optimization. Suppose we have a scalar-valued function $f: \R^n \to \R$, and suppose we want to minimize (or maximize) $f$. For instance, in machine learning, we could have a big neural network (NN) with a vector $x$ of a million parameters, and one tries to minimize a ``loss'' function $f$ that compares the NN output to the desired results on ``training'' data. The most basic idea in optimization is to go ``downhill'' (see diagram) to make $f$ as small as possible. If we can take the gradient of this function $f$, to go ``downhill'' we consider $- \nabla f$, the direction of \emph{steepest descent}, as depicted in Fig.~\ref{fig:steepest-descent}. 
    \todo{It would be nice if all the ``dots" were on contours, and we can see the tangents and normals to the contours}

Then, even if we have a million parameters, we can evolve all of them simultaneously in the downhill direction. It turns out that calculating all million derivatives costs about the same as evaluating the function at a point once (using reverse-mode/adjoint/left-to-right/backpropagation methods). Ultimately, this makes large-scale optimization practical for training neural nets, optimizing shapes of airplane wings, optimizing portfolios, etc.

Of course, there are many practical complications that make nonlinear optimization tricky (far more than can be covered in a single lecture, or even in a whole course!), but we give some examples here.

\begin{itemize}
    \item For instance, even though we can compute the ``downhill direction'', how far do we need to step in that direction? (In machine learning, this is sometimes called the ``learning rate.'') Often, you want to take ``as big of a step as you can'' to speed convergence, but you  don't want the step to be too big because $\nabla f$ only tells you a \emph{local} approximation of~$f$. There are many different ideas of how to determine this:
\begin{itemize}
    \item Line search: using a 1D minimization to determine how far to step.
    \item A ``trust region'' bounding the step size (where we trust the derivative-based approximation of $f$). There are many techniques to evolve the size of the trust region as optimization progresses.
\end{itemize}
\item We may also need to consider constraints, for instance 
minimizing $f(x)$ subject to $g_k(x) \leq 0$  or $h_k(x)=0$,
 known as inequality/equality constraints.
Points $x$ satisfying the constraints are called ``feasible''.   One typically uses a combination of $\nabla f$ and $\nabla g_k$ to approximate (e.g.~linearize) the problem and make progress towards the best feasible point.

\item If you just go straight downhill, you might ``zig-zag'' along narrow valleys, making convergence very slow. There are a few options to combat this, such as ``momentum'' terms and conjugate gradients. Even fancier than these techniques, one might estimate second-derivative ``Hessian matrices'' from a sequence of $\nabla f$ values---a famous version of this is known as the BFGS algorithm---and use the Hessian to take approximate Newton steps (for the root $\nabla f = 0$).  (We'll return to Hessians in a later lecture.)
\item Ultimately, there are a lot of techniques and a zoo of competing algorithms that you might need to experiment with to find the best approach for a given problem.  (There are many books on optimization algorithms, and even a whole book can only cover a small slice of what is out there!)
\end{itemize}

Some parting advice: Often the main trick is less about the choice of algorithms than it is about finding the right mathematical formulation of your \emph{problem}---e.g.~what function, what constraints, and what parameters should you be considering---to match your problem to a good algorithm. However, if you have \textit{many} ($\gg 10$) parameters, \textit{try hard} to use an analytical gradient (not finite differences), computed efficiently in reverse mode.

\subsubsection{Engineering/Physical Optimization}

There are many, many applications of optimization besides machine learning (fitting models to data).  It is interesting to also consider engineering/physical optimization. (For instance, suppose you want to make an airplane wing that is as strong as possible.) The general outline of such problems is typically:
\begin{enumerate}
    \item You start with some design parameters $\mathbf{p}$, e.g.~describing the geometry, materials, forces, or other degrees of freedom. 
    \item These $\mathbf{p}$ are then used in some  physical model(s), such as solid mechanics, chemical reactions, heat transport, electromagnetism, acoustics, etc.  For example, you might have a linear model of the form $A(\mathbf{p}) x = b(\mathbf{p})$ for some matrix $A$ (typically very large and sparse).
    \item The solution of the physical model is a solution $x(\mathbf{p})$.  For example, this could be the mechanical stresses, chemical concentrations, temperatures, electromagnetic fields, etc.
    \item The physical solution $x(\mathbf{p})$ is the input into some design objective $f(x(\mathbf{p}))$ that you want to improve/optimize.  For instance, strength, speed power, efficiency, etc. 
    \item To maximize/minimize $f(x(\mathbf{p}))$, one uses the gradient $\nabla_{\mathbf{p}}f$, computed using reverse-mode/``adjoint'' methods, to update the parameters $\mathbf{p}$ and improve the design.
\end{enumerate}
\noindent As a fun example, researchers have even applied ``topology optimization'' to design a chair, optimizing every voxel of the design---the parameters $\mathbf{p}$ represent the material present (or not) in every voxel, so that the optimization discovers not just an optimal shape but an optimal \emph{topology} (how materials are connected in space, how many holes there are, and so forth)---to support a given weight with \textit{minimal material}. To see it in action, watch this \href{https://www.youtube.com/watch?v=bJ_nSSBl040&embeds_referring_euri=https%3A%2F%2Fdocs.google.com%2F&embeds_referring_origin=https%3A%2F%2Fdocs.google.com&source_ve_path=Mjg2NjY&feature=emb_logo}{chair-optimization video}.   (People have applied such techniques to much more practical problems as well, from airplane wings to optical communications.)

\subsection{Reverse-mode ``Adjoint'' Differentiation}
\label{sec:adjoint-method}

But what is adjoint differentiation---the method of differentiating that makes these applications actually feasible to solve? Ultimately, it is yet another example of left-to-right/reverse-mode differentiation, essentially applying the chain rule from outputs to inputs. Consider, for example, trying to compute the gradient $\nabla g$ of the scalar-valued function 
$$
g(p) = f(\underbrace{A(p)^{-1} b}_x) \, .
$$
where $x$ solves $A(p) x  = b$ (e.g.~a parameterized physical model as in the previous section) and $f(x)$ is a scalar-valued function of $x$ (e.g.~an optimization objective depending on our physics solution).  For example, this could arise in an optimization problem 
$$
 \min_p g(p) \Longleftrightarrow \substack{  \text{ \normalsize  $\displaystyle \min_p f(x)$} \\ \text{ subject to } A(p)x = b }\; ,
$$
for which the gradient $\nabla g$ would be helpful to search for a local minimum.
The chain rule for $g$ corresponds to the following conceptual chain of dependencies:
\begin{align*}
\text{change $dg$ in $g$} &\longleftarrow 
\text{change $dx$ in $x = A^{-1} b$}  \\
&\longleftarrow \text{change $d(A^{-1})$ in $A^{-1}$} \\
&\longleftarrow 
\text{change $dA$ in $A(p)$}  \\
&\longleftarrow \text{change $dp$ in $p$}
\end{align*}
which is expressed by the equations:
\begin{align*}
\d g &= f'(x) [\d x] &\text{ } & dg \longleftarrow dx    \\
& = f'(x) [\d (A^{-1}) b] &\text{ } & dx \longleftarrow d(A^{-1})   \\
& = - \underbrace{f'(x) A^{-1}}_{v^T} \d A \, A^{-1} b &\text{ } & dA^{-1} \longleftarrow dA \\
& = - v^T \underbrace{A'(p)[dp]}_{dA} \, A^{-1} b &\text{ } &dA \longleftarrow dp \, .
\end{align*}
Here, we are defining the row vector $v^T = f'(x) A^{-1}$, and we have used the differential of a matrix inverse $\d(A^{-1})=-A^{-1}\,dA\,A^{-1}$ from Sec.~\ref{sec:jacobian-inverse}.

Grouping the terms left-to-right, we first solve the ``adjoint'' (transposed) equation $A^T v = f'(x)^T = \nabla_x f$ for $v$,
and then we obtain $\d g = - v^T \d A \, x$. Because the derivative $A'(p)$ of a matrix with respect to a vector is awkward to write explicitly, it is convenient to examine this object one parameter at a time.  For any given parameter $p_k$, $\partial g/\partial p_k = -v^T (\partial A/\partial p_k) x$ (and in many applications $\partial A /\partial p_k$ is very sparse); here, ``dividing by'' $\partial p_k$ works because this is a scalar factor that commutes with the other linear operations. That is, it takes only \emph{two solves} to get both $g$ and $\nabla g$: one for solving $A x = b$ to find $g(p)=f(x)$, and another with $A^T$ for $v$, after which all of the derivatives $\partial g/\partial p_k$ are just some cheap dot products.

Note that you should \textbf{\textit{not}} use right-to-left ``forward-mode'' derivatives with lots of parameters, because 
\[
\frac{\partial g}{\partial p_k} = - f'(x) \left(A^{-1} \frac{\partial A}{\partial p_k}x\right)
\]
represents one solve per parameter $p_k$! As discussed in Sec.~\ref{sec:forward-vs-reverse}, right-to-left (a.k.a.~forward mode) is better when there is one (or few) input parameters $p_k$ and many outputs, while left-to-right ``adjoint'' differentiation (a.k.a.~reverse mode) is better when there is one (or few) output values and many input parameters. (In Sec.~\ref{sec:dual-AD}, we will discuss using \href{https://en.wikipedia.org/wiki/Dual_number}{dual numbers} for differentiation, and this also corresponds to forward mode.) 

Another possibility that might come to mind is to use finite differences (as in Sec.~\ref{sec:finitedifference}), but you should not use this if you have lots of parameters! Finite differences would involve a calculation of something like
\[
\frac{\partial g}{\partial p_k} \approx [g(p + \epsilon e_k) - g(p)]/\epsilon,
\]
where $e_k$ is a unit vector in the $k$-th direction and $\epsilon$ is a small number.
This, however, requires one  solve for each parameter $p_k$, just like forward-mode differentiation.  (It becomes even more expensive if you use fancier higher-order finite-difference approximations in order to obtain higher accuracy.)

\subsubsection{Nonlinear equations}

You can also apply adjoint/reverse differentiation to nonlinear equations. For instance, consider the gradient of the scalar function $g(p) = f(x(p))$, where $x(p)\in \R^n$ solves some system of $n$ equations $h(p,x) = 0 \in \R^n$. By the chain rule, $$h(p,x) = 0 \implies \frac{\partial h}{\partial p} \d p + \frac{\partial  h}{\partial x} \d x = 0 \implies \d x = -\left(\frac{\partial  h}{\partial x}\right)^{-1}  \frac{\partial h}{\partial p}   \d p \,.$$ (This is an instance of the \href{https://en.wikipedia.org/wiki/Implicit_function_theorem}{Implicit Function Theorem}: as long as 
$\frac{\partial  h}{\partial x}$ is nonsingular, we can locally define a function $x(p)$ from an implicit equation $h=0$, here by linearization.) 
Hence, 
\[
\d g = f'(x) \d x = - \underbrace{ f'(x) \left(\frac{\partial  h}{\partial x}\right)^{-1} }_{v^T} \frac{\partial h}{\partial p} \d p \, .
\]
Associating left-to-right again leads to a single ``adjoint'' equation: $(\partial h/\partial x)^T v = f'(x)^T = \nabla_x f$. In other words, it again only takes two solves to get both $g$ and $\nabla g$---one nonlinear ``forward'' solve for $x$ and one linear ``adjoint'' solve for $v$! Thereafter, all derivatives $\partial g/\partial p_k$ are cheap dot products.   (Note that the linear ``adjoint'' solve involves the transposed Jacobian $\partial h/\partial x$.  Except for the transpose, this is very similar to the cost of a single Newton step to solve $h=0$ for $x$.  So the adjoint problem should be cheaper than the forward problem.)

\subsubsection{Adjoint methods and AD}

If you use automatic differentiation (AD) systems, why do you need to learn this stuff?  Doesn't the AD do everything for you?
In practice, however, it is often helpful to understand adjoint methods even if you use automatic differentiation. Firstly, it helps you understand when to use forward- vs.~reverse-mode automatic differentiation. Secondly, many physical models call large software packages written over the decades in various languages that \emph{cannot be differentiated automatically} by AD.  You can typically correct this by just supplying a ``vector--Jacobian product'' $y^T \d x$ for this physics, or even just \emph{part} of the physics, and then AD will differentiate the rest and apply the chain rule for you. Lastly, often models involve approximate calculations (e.g. for the iterative solution of linear or nonlinear equations, numerical integration, and so forth), but AD tools often don't ``know'' this and spend extra effort trying to differentiate the error in your approximation; in such cases, manually written derivative rules can sometimes be much more efficient.  (For example, suppose your model involves solving a nonlinear system $h(x,p) = 0$ by an iterative approach like Newton's method.  Naive AD will be very inefficient because it will attempt to differentiate through all your Newton steps.  Assuming that you converge your Newton solver to enough accuracy that the error is negligible, it is much more efficient to perform differentiation via the implicit-function theorem as described above, leading to a single linear adjoint solve.)

\subsubsection{Adjoint-method example}

To finish off this section of the notes, we conclude with an example of how to use this ``adjoint method'' to compute a derivative efficiently. Before working through the example, we first state the problem and highly recommend trying it out before reading the solution.

\begin{problem}\label{PSETexample}
Suppose that $A(p)$ takes a vector $p\in\mathbb{R}^{n-1}$ and returns
the $n\times n$ tridiagonal real-symmetric matrix
\[
A(p)=\left(\begin{array}{ccccc}
a_{1} & p_{1}\\
p_{1} & a_{2} & p_{2}\\
 & p_{2} & \ddots & \ddots\\
 &  & \ddots & a_{n-1} & p_{n-1}\\
 &  &  & p_{n-1} & a_{n}
\end{array}\right),
\]
where $a\in\mathbb{R}^{n}$ is some constant vector. Now, define
a scalar-valued function $f(p)$ by 
\[
g(p)=\left(c^{T}A(p)^{-1}b\right)^{2}
\]
for some constant vectors $b,c\in\mathbb{R}^{n}$ (assuming we choose
$p$ and $a$ so that $A$ is invertible). Note that, in practice,
$A(p)^{-1}b$ is \emph{not }computed by explicitly inverting the matrix
$A$---instead, it can be computed in $\Theta(n)$ (i.e., roughly
proportional to $n$) arithmetic operations using Gaussian elimination
that takes advantage of the ``sparsity'' of $A$ (the pattern of
zero entries), a ``tridiagonal solve.''
\begin{enumerate}
\item[(a)] Write down a formula for computing $\partial g/\partial p_{1}$ (in
terms of matrix--vector products and matrix inverses). (Hint: once
you know $dg$ in terms of $dA$, you can get $\partial g/\partial p_{1}$
by ``dividing'' both sides by $\partial p_{1}$, so that $dA$ becomes
$\partial A/\partial p_{1}$.)
\item[(b)] Outline a sequence of steps to compute both $g$ and $\nabla g$ (with
respect to $p$) using only \emph{two} tridiagonal solves $x=A^{-1}b$
and an ``adjoint'' solve $v=A^{-1}\text{(something)}$, plus $\Theta(n)$
(i.e., roughly proportional to $n$) additional arithmetic operations.
\item[(c)] Write a program implementing your $\nabla g$ procedure (in Julia,
Python, Matlab, or any language you want) from the previous part.
(You don't need to use a fancy tridiagonal solve if you don't know
how to do this in your language; you can solve $A^{-1}\text{(vector)}$
inefficiently if needed using your favorite matrix libraries.) Implement
a finite-difference test: Choose $a,b,c,p$ at random, and check that
$\nabla g\cdot\delta p\approx g(p+\delta p)-g(p)$ (to a few digits)
for a randomly chosen small $\delta p$.
\end{enumerate}
\end{problem}

\textbf{\cref{PSETexample}(a) Solution:} From the chain rule and the formula for the differential of a matrix inverse, we have $dg = -2(c^T A^{-1} b) c^T A^{-1} dA\,A^{-1} b$ (noting that $c^T A^{-1} b$ is a scalar so we can commute it as needed).  Hence
\begin{align*}
\frac{\partial g}{\partial p_1} &= \underbrace{-2(c^T A^{-1} b) c^T A^{-1}}_{v^T} \frac{\partial A}{\partial p_1} \underbrace{A^{-1} b}_x \\
&=  v^T \underbrace{\left(\begin{array}{ccccc}
0 & 1\\
1 & 0 & 0\\
 & 0 & \ddots & \ddots\\
 &  & \ddots & 0 & 0\\
 &  &  & 0 & 0
\end{array}\right)}_{\frac{\partial A}{\partial p_1}} x = \boxed{v_1 x_2 + v_2 x_1} \, ,
\end{align*}
where we have simplified the result in terms of $x$ and $v$ for the next part.

\textbf{\cref{PSETexample}(b) Solution:} Using the notation from the previous part, exploiting the fact that $A^T = A$, we can choose $\boxed{v = A^{-1} [-2(c^T x) c]}$, which is a single tridiagonal solve.  Given $x$ and $v$, the results of our two $\Theta(n)$ tridiagonal solves, we can compute each component of the gradient similar to above by $\boxed{\partial g/\partial p_k = v_k x_{k+1} + v_{k+1} x_k}$ for $k=1,\ldots,n-1$, which costs $\Theta(1)$ arithmetic per $k$ and hence $\Theta(n)$ arithmetic to obtain all of $\nabla g$.

\textbf{\cref{PSETexample}(c) Solution:} See the \href{https://nbviewer.org/github/mitmath/matrixcalc/blob/iap2023/psets/pset2sol.ipynb}{Julia solution notebook (Problem~1)} from our IAP 2023 course (which calls the function $f$ rather than $g$).

\pagebreak
\section{Derivative of Matrix Determinant and Inverse}
\subsection{Two Derivations}

This section of notes follows \href{https://rawcdn.githack.com/mitmath/matrixcalc/b08435612045b17745707f03900e4e4187a6f489/notes/determinant_and_inverse.html}{this} Julia notebook. This notebook is a little bit short, but is an important and useful calculation.

\begin{theorem}
    Given $A$ is a square matrix, we have 
    \[
    \nabla (\det A) = \cofactor(A) = (\det A)A^{-T} := \adj(A^T) = \adj(A)^T 
    \]
    where $\adj$ is the ``adjugate''.   (You may not have heard of the matrix adjugate, but this formula tells us that it is simply $\adj(A) = \det(A) A^{-1}$, or $\cofactor(A) = \adj(A^T)$.) Furthermore, 
    \[
    \d (\det A) = \tr(\det(A) A^{-1} \d A) = \tr (\adj(A) \d A) = \tr(\cofactor (A)^T \d A).
    \]
\end{theorem}

You may remember that each entry $(i,j)$ of the cofactor matrix is $(-1)^{i + j}$ times the determinant obtained by deleting row $i$ and column $j$ from $A$. Here are some $2 \times 2$ calculuations to obtain some intuition about these functions: 
\begin{align}
    M &= \begin{pmatrix}
        a & c \\ b & d    
    \end{pmatrix} \\
    \implies \cofactor(M) &= \begin{pmatrix}
         d & -c \\ -b & a
    \end{pmatrix}  \\
    \adj(M) &= \begin{pmatrix}
d & -b \\ -c & a
\end{pmatrix} \\
(M)^{-1}  &= \frac{1}{ad-bc} \begin{pmatrix}
    d & -b \\ -c & a
\end{pmatrix}.
\end{align}

Numerically, as is done in the notebook, you can construct a random $n \times n$ matrix $A$ (say, $9 \times 9$), consider e.g.~$\d A = .00001 A$, and see numerically that 
\[
\det(A + \d A) - \det (A) \approx \tr(\adj(A) \d A),
\]
which numerically supports our claim for the theorem. 

We now prove the theorem in two ways. Firstly, there is a direct proof where you just differentiate the scalar with respect to every input using the \href{https://en.wikipedia.org/wiki/Laplace_expansion}{cofactor expansion} of the determinant based on the $i$-th row. Recall that 
\[
\det (A) = A_{i1} C_{i1} +A_{i2} C_{i2} + \dots + A_{in} C_{in}.
\]
Thus, 
\[
\frac{\partial \det A}{ \partial A_{ij}} = C_{ij} \implies \nabla (\det A) = C, 
\]
the cofactor matrix.  (In computing these partial derivatives, it's important to remember that the cofactor $C_{ij}$ contains no elements of $A$ from row~$i$ or column~$j$.  So, for example, $A_{i1}$ only appears explicitly in the first term, and not hidden in any of the $C$ terms in this expansion.)

There is also a fancier proof of the theorem using linearization near the identity. Firstly, note that it is easy to see from the properties of determinants that $$\det(I + \d A) - 1 = \tr(\d A),$$ and thus 
\begin{align*}
    \det(A + A(A^{-1} \d A)) - \det (A) &= \det(A) (\det (I + A^{-1} \d A) - 1) \\
    &= \det(A) \tr(A^{-1} \d A) = \tr(\det (A) A^{-1} \d A) \\
    &= \tr(\adj(A) \d A).
\end{align*}
This also implies the theorem.

\subsection{Applications}
\subsubsection{Characteristic Polynomial}

We now use this as an application to find the derivative of a characteristic polynomial evaluated at $x$. Let $p(x) = \det(xI -A)$, a scalar function of $x$. Recall that through factorization, $p(x)$ may be written in terms of eigenvalues $\lambda_i$. So we may ask: what is the derivative of $p(x)$, the characteristic polynomial at $x$? Using freshman calculus, we could simply compute 
\[
\frac{\d}{\d x} \prod_i (x-\lambda_i) = \sum_i \prod_{j\neq i} (x-\lambda_j)  = \prod (x-\lambda_i) \{\sum_i (x- \lambda_i)^{-1}\},
\]
as long as $x \neq \lambda_i$.

This is a perfectly good simply proof, but with our new technology we have a new proof:
\begin{align*}
    \d (\det (x I - A)) &= \det(x I - A) \tr((xI - A)^{-1} \d (x I - A)) \\
    &= \det(xI - A) \tr(x I - A)^{-1} \d x.
\end{align*}
Note that here we used that $\d (x I - A) = \d x \, I$ when $A$ is constant and $\tr (A \d x) = \tr(A) \d x$ since $\d x$ is a scalar.

We may again check this computationally as we do in the notebook.

\subsubsection{The Logarithmic Derivative}

We can similarly compute using the chain rule that 
\[
\d (\log (\det (A))) = \frac{\d (\det A)}{ \det A} = \det (A^{-1}) \d (\det (A)) = \tr(A^{-1} \d A).
\]
The logarithmic derivative shows up a lot in applied mathematics. Note that here we use that $\frac{1}{\det A} = \det(A^{-1})$ as $1 = \det(I) = \det(AA^{-1}) = \det (A) \det(A^{-1}).$

For instance, recall Newton's method to find roots $f(x)=0$ of single-variable real-valued functions $f(x)$ by taking a sequence of steps $x \to x + \delta x$.  The key formula in Newton's method is $\delta x = f'(x)^{-1}f(x)$, but this is the same as $\frac{1}{(\log f(x))'}$. So, derivatives of log determinants show up in finding \emph{roots of determinants}, i.e.~for $f(x) = \det M(x)$.  When $M(x) = A - x I$, roots of the determinant are eigenvalues of $A$.  For more general functions $M(x)$, solving $\det M(x) = 0$ is therefore called a \emph{nonlinear eigenproblem}.

\subsection{Jacobian of the Inverse}
\label{sec:jacobian-inverse}

Lastly, we compute the derivative (as both a linear operator and an explicit Jacobian matrix) of the inverse of a matrix. There is a neat trick to obtain this derivative, simply from the property $A^{-1}A = I$ of the inverse.  By the product rule, this implies that
\begin{align*}
\d (A^{-1} A) &= d(I) = 0 = \d (A^{-1}) A + A^{-1} \d A \\
& \implies \boxed{\d (A^{-1}) = (A^{-1})'[dA] = - A^{-1} \, \d A \, A^{-1} }\, .
\end{align*}
Here, the second line defines a perfectly good linear operator for the derivative $(A^{-1})'$, but if we want we can rewrite this as an explicit Jacobian matrix by using Kronecker products acting on the ``vectorized'' matrices as we did in Sec.~\ref{sec:kronecker}:
\[
\vecm\left(\d (A^{-1})\right) = \vecm\left(-A^{-1} (\d A) A^{-1}\right) = \underbrace{- (A^{-T} \otimes A^{-1})}_\mathrm{Jacobian} \vecm(\d A) \, ,
\]
where $A^{-T}$ denotes $(A^{-1})^T = (A^T)^{-1}$.
One can check this formula numerically, as is done in the notebook.

In practice, however, you will probably find that the operator expression $- A^{-1} \, \d A \, A^{-1}$ is more useful than explicit Jacobian matrix for taking derivatives involving matrix inverses.  For example, if you have a matrix-valued function $A(t)$ of a scalar parameter $t \in \mathbb{R}$, you immediately obtain $\frac{d(A^{-1})}{dt} = -A^{-1} \frac{dA}{dt} A^{-1}$.   A more sophisticated application is discussed in Sec.~\ref{sec:adjoint-method}.

\pagebreak

\section{Forward and Reverse-Mode Automatic Differentiation}
\label{sec:AD}
The first time that Professor Edelman had heard about automatic differentiation (AD), it was easy for him to imagine what it was \ldots but what he imagined was wrong! In his head, he thought it was straightforward symbolic differentiation applied to code---sort of like executing Mathematica or Maple, or even just automatically doing what he learned to do in his calculus class. For instance, just plugging in functions and their domains from something like the following first-year calculus table:
\begin{table}[h]
\begin{center}
\begin{tabular}{|l|l|}
\hline
Derivative                    & Domain                                         \\ \hline
$(\sin x)' = \cos x$          & $- \infty< x < \infty$                   \\ \hline
$(\cos x)' = - \sin x$        & $-\infty< x < \infty$                    \\ \hline
$(\tan x)' = \sec^2 x$        & $x\neq \frac{\pi}{2} + \pi n, n \in \ZZ$ \\ \hline
$(\cot x)' = - \csc^2 x$      & $x\neq \pi n, n \in \ZZ$                 \\ \hline
$(\sec x)' = \tan x \sec x$   & $x \neq \frac{\pi}{2} + \pi n, n\in \ZZ$ \\ \hline
$(\csc x)' = - \cot x \csc x$ & $x\neq \pi n, n \in \ZZ$                 \\ \hline
\end{tabular}
\end{center}
\end{table}

\noindent And in any case, if it wasn't just like executing Mathematica or Maple, then  it must be finite differences, like one learns in a numerical computing class (or as we did in Sec.~\ref{sec:finitedifference}).

It turns out that it is definitely \emph{not} finite differences---AD algorithms are generally exact (in exact arithmetic, neglecting roundoff errors), not approximate.  But it also doesn't look much like conventional symbolic algebra: the computer doesn't really construct a big ``unrolled'' symbolic expression and then differentiate it, the way you might imagine doing by hand or via computer-algebra software.  For example, imagine a computer program that computes $\det A$ for an $n \times n$ matrix---writing down the ``whole'' symbolic expression isn't possible until the program runs and $n$ is known (e.g.~input by the user), and in any case a naive symbolic expression would require $n!$ terms.  Thus, AD systems have to deal with computer-programming constructs like loops, recursion, and problem sizes $n$ that are unknown until the program runs, while at the same time avoiding constructing symbolic expressions whose size becomes prohibitively large.  (See Sec.~\ref{sec:babylonian} for an example that looks very different from the formulas you differentiate in first-year calculus.)  Design of AD systems often ends up being more about compilers than about calculus!

\subsection{Automatic Differentiation via Dual Numbers}
\label{sec:dual-AD}

\emph{(This lecture is accompanied by a Julia ``notebook'' showing the results of various computational experiments, which can be found on the course web page.  Excerpts from those experiments are included below.)}

One AD approach that can be explained relatively simply is ``forward-mode'' AD, which is implemented by carrying out the computation of $f'$ in \emph{tandem} with the computation of $f$.  One augments every intermediate value~$a$ in the computer program with another value~$b$ that represents its derivative, along with chain rules to propagate these derivatives through computations on values in the program. 
It turns out that this can be thought of as replacing real numbers (values~$a$) with a new kind of ``dual number'' $D(a,b)$ (values \& derivatives) and corresponding arithmetic rules, as explained below.

\subsubsection{Example: Babylonian square root}  
\label{sec:babylonian}

We start with a simple example, an algorithm for the square-root function, where a practical method of automatic differentiation came as both a mathematical surprise and a computing wonder for Professor Edelman. In particular, we consider the ``Babylonian'' algorithm to compute $\sqrt{x}$, known for millennia (and later revealed as a special case of Newton's method applied to $t^2 - x = 0$): simply repeat $t \leftarrow (t + x/t)/2$ until $t$ converges to $\sqrt{x}$ to any desired accuracy. Each iteration has one addition and two divisions. For illustration purposes, 10 iterations suffice.    Here is a short program in Julia that implements this algorithm, starting with a guess of $1$ and then performing $N$ steps (defaulting to $N=10$):
\begin{minted}{jlcon}
julia> function Babylonian(x; N = 10) 
           t = (1+x)/2   # one step from t=1
           for i = 2:N   # remaining N-1 steps
               t = (t + x/t) / 2
           end    
           return t
       end
\end{minted}
If we run this function to compute the square root of $x=4$, we will see that it converges very quickly: for only $N=3$ steps, it obtains the correct answer~($2$) to nearly~3 decimal places, and well before $N=10$ steps it has converged to $2$ within the limits of the accuracy of computer arithmetic (about 16 digits).  In fact, it roughly doubles the number of correct digits on every step:
\begin{minted}{jlcon}
julia> Babylonian(4, N=1)
2.5

julia> Babylonian(4, N=2)
2.05

julia> Babylonian(4, N=3)
2.000609756097561

julia> Babylonian(4, N=4)
2.0000000929222947

julia> Babylonian(4, N=10)
2.0
\end{minted}

Of course, any first-year calculus student knows the derivative of the square root, $(\sqrt{x})' = 0.5/\sqrt{x}$, which we could compute here via \texttt{0.5 / Babylonian(x)},
but we want to know how we can obtain this derivative \emph{automatically}, directly from the Babylonian algorithm itself.  If we can figure out how to easily and efficiently pass the chain rule through this algorithm, then we will begin to understand how AD can also differentiate much more complicated computer programs for which no simple derivative formula is known.

\subsubsection{Easy forward-mode AD}

The basic idea of carrying the chain rule through a computer program is very simple: replace every number with \emph{two} numbers, one which keeps track of the \emph{value} and one which tracks the \emph{derivative} of that value.  The values are computed the same way as before, and the derivatives are computed by carrying out the chain rule for elementary operations like $+$ and $/$.

In Julia, we can implement this idea by defining a new type of number, which we'll call \texttt{D}, that encapsulates a value \texttt{val} and a derivative \texttt{deriv}. 
\begin{minted}{jlcon}
julia> struct D <: Number
           val::Float64
           deriv::Float64
       end
\end{minted}
(A detailed explanation of Julia syntax can be found \href{https://julialang.org/learning/}{elsewhere}, but hopefully you can follow the basic ideas even if you don't understand every punctuation mark.)
A quantity \texttt{x = D(a,b)} of this new type has two components \texttt{x.val = a} and \texttt{x.deriv = b}, which we will use to represent values and derivatives, respectively.
The \texttt{Babylonian} code only uses two arithmetic operations, $+$ and $/$, so we just need to overload the built-in (``\texttt{Base}'') definitions of these in Julia to include new rules for our \texttt{D} type:
\begin{minted}{jlcon}
julia> Base.:+(x::D, y::D) = D(x.val+y.val, x.deriv+y.deriv)
       Base.:/(x::D, y::D) = D(x.val/y.val, (y.val*x.deriv - x.val*y.deriv)/y.val^2)
\end{minted}
If you look closely, you'll see that the values are just added and divided in the ordinary way, while the derivatives are computed using the sum rule (adding the derivatives of the inputs) and the quotient rule, respectively.   We also need one other technical trick: we need to define \href{https://docs.julialang.org/en/v1/manual/conversion-and-promotion/}{``conversion'' and ``promotion'' rules} that tell Julia how to combine \texttt{D} values with ordinary real numbers, as in expressions like $x+1$ or $x/2$:
\begin{minted}{jlcon}
julia> Base.convert(::Type{D}, r::Real) = D(r,0)
       Base.promote_rule(::Type{D}, ::Type{<:Real}) = D
\end{minted}
This just says that an ordinary real number $r$ is combined with a \texttt{D} value by first converting $r$ to \texttt{D(r,0)}: the value is~$r$ and the derivative is $0$ (the derivative of any constant is zero).

Given these definitions, we can now plug a \texttt{D} value into our \emph{unmodified} \texttt{Babylonian} function, and it will ``magically'' compute the derivative of the square root.  Let's try it for $x = 49 = 7^2$:
\begin{minted}{jlcon}
julia> x = 49
49

julia> Babylonian(D(x,1))
D(7.0, 0.07142857142857142)
\end{minted}
We can see that it correctly returned a value of \texttt{7.0} and a derivative of \texttt{0.07142857142857142}, which indeed matches the square root $\sqrt{49}$ and its derivative $0.5/\sqrt{49}$:
\begin{minted}{jlcon}
julia> (√x, 0.5/√x)
(7.0, 0.07142857142857142)
\end{minted}
Why did we input \texttt{D(x,1)}?  Where did the $1$ come from?   That's simply the fact that the derivative of the input $x$ with respect to \emph{itself} is $(x)' = 1$, so this is the starting point for the chain rule.

In practice, all this (and more) has already been implemented in the \href{https://github.com/JuliaDiff/ForwardDiff.jl}{ForwardDiff.jl} package in Julia (and in many similar software packages in a variety of languages).  That package hides the implementation details under the hood and explicitly provides a function to compute the derivative. For example:
\begin{minted}{jlcon}
julia> using ForwardDiff

julia> ForwardDiff.derivative(Babylonian, 49)
0.07142857142857142
\end{minted}
Essentially, however, this is the same as our little \texttt{D} implementation, but implemented with greater generality and sophistication (e.g.~chain rules for more operations, support for more numeric types, partial derivatives with respect to multiple variables, etc.): just as we did, ForwardDiff augments every value with a second number that tracks the derivative, and propagates both quantities through the calculation.

We could have also implemented the same idea specifically for the Bablylonian algorithm, by writing a new function \texttt{dBabylonian} that tracks both the variable $t$ and its derivative $t' = dt/dx$ through the course of the calculation:
\begin{minted}{jlcon}
julia> function dBabylonian(x; N = 10) 
           t = (1+x)/2
           t′ = 1/2
           for i = 1:N
               t = (t+x/t)/2
               t′= (t′+(t-x*t′)/t^2)/2
           end    
           return t′
       end

julia> dBabylonian(49)
0.07142857142857142
\end{minted}
This is doing \emph{exactly} the same calculations as calling \texttt{Babylonian(D(x,1))} or \texttt{ForwardDiff.derivative(Babylonian, 49)}, but needs a lot more human effort---we'd have to do this for every computer program we write, rather than implementing a new number type \emph{once}.

\subsubsection{Dual numbers}

There is a pleasing algebraic way to think about our new number type $D(a,b)$ instead of the ``value \& derivative'' viewpoint above.   Remember how a complex number $a + bi$ is formed from two real numbers $(a,b)$ by defining a special new quantity $i$ (the imaginary unit) that satisfies $i^2 = -1$, and all the other complex-arithmetic rules follow from this?   Similarly, we can think of $D(a,b)$ as $a + b \epsilon$, where $\epsilon$ is a new ``infinitesimal unit'' quantity that satisfies $\epsilon^2 = 0$.   This viewpoint is called a \textbf{dual number}.

Given the elementary rule $\epsilon^2 = 0$, the other algebraic rules for dual numbers immediately follow:
\begin{align*}
    (a + b \epsilon) \pm (c + d \epsilon) &= (a \pm c) + (b \pm d) \epsilon \\
    (a + b \epsilon) \cdot (c + d\epsilon) &= (ac) + (bc + ad) \epsilon \\
    \frac{a + b \epsilon}{c + d \epsilon} &= \frac{a + b \epsilon}{c + d \epsilon} \cdot \frac{c - d \epsilon}{c - d \epsilon} = \frac{(a + b \epsilon)(c - d \epsilon)}{c^2}  =     
    \frac{a}{c} + \frac{bc - ad}{c^2}\epsilon.
\end{align*}
The $\epsilon$ coefficients of these rules correspond to the sum/difference, product, and quotient rules of differential calculus!

In fact, these are \emph{exactly} the rules we implemented above for our \texttt{D} type. We were only missing the rules for subtraction and multiplication, which we can now include:
\begin{minted}{jlcon}
julia> Base.:-(x::D, y::D) = D(x.val - y.val, x.deriv - y.deriv)
       Base.:*(x::D, y::D) = D(x.val*y.val, x.deriv*y.val + x.val*y.deriv)
\end{minted}
It's also nice to add a \href{https://docs.julialang.org/en/v1/manual/types/#man-custom-pretty-printing}{``pretty printing''} rule to make Julia display dual numbers as \texttt{a + bϵ} rather than as \texttt{D(a,b)}:
\begin{minted}{jlcon}
julia> Base.show(io::IO, x::D) = print(io, x.val, " + ", x.deriv, "ϵ")
\end{minted}
Once we implement the multiplication rule for dual numbers in Julia, then $\epsilon^2 = 0$ follows from the special case $a = c=0$ and $b=d=1$:
\begin{minted}{jlcon}
julia> ϵ = D(0,1)
0.0 + 1.0ϵ

julia> ϵ * ϵ 
0.0 + 0.0ϵ

julia> ϵ^2
0.0 + 0.0ϵ
\end{minted}
(We didn't define a rule for powers $D(a,b)^n$, so how did it compute \texttt{ϵ}$^2$? The answer is that Julia implements $x^n$ via repeated multiplication by default, so it sufficed to define the $*$ rule.)  Now, we can compute the derivative of the Babylonian algorithm at $x = 49$ as above by:
\begin{minted}{jlcon}
julia> Babylonian(x + ϵ)
7.0 + 0.07142857142857142ϵ
\end{minted}
with the ``infinitesimal part'' being the derivative $0.5/\sqrt{49} = 0.0714\cdots$.

A nice thing about this dual-number viewpoint is that it corresponds directly to our notion of a derivative as linearization:
$$
f(x + \epsilon) = f(x) + f'(x) \epsilon + \mbox{(higher-order terms)} \, ,
$$
with the dual-number rule $\epsilon^2 = 0$ corresponding to dropping the higher-order terms.

\subsection{Naive symbolic differentiation}

Forward-mode AD implements the exact analytical derivative by propagating chain rules, but it is completely different from what many people \emph{imagine} AD might be: evaluating a program \emph{symbolically} to obtain a giant symbolic expression, and \emph{then} differentiating this giant expression to obtain the derivative.   A basic issue with this approach is that the size of these symbolic expressions can quickly explode as the program runs.  Let's see what it would look like for the Babylonian algorithm.

Imagine inputting a ``symbolic variable'' $x$ into our \texttt{Babylonian} code, running the algorithm, and writing a big algebraic expression for the result.   After only one step, for example, we would get $(x + 1)/2$.   After two steps, we would get $((x+1)/2 + 2x/(x+1))/2$, which simplifies to a ratio of two polynomials (a ``rational function''):
$$
\frac{x^2 + 6x + 1}{4(x+1)} \, .
$$
Continuing this process by hand is quite tedious, but fortunately the computer can do it for us (as shown in the accompanying Julia notebook). Three Babylonian iterations yields:
$$\frac{x^{4} + 28 x^{3} + 70 x^{2} + 28 x + 1}{8 \left(x^{3} + 7 x^{2} + 7 x + 1\right)} \, ,
$$
four iterations gives
$$
\frac{x^{8} + 120 x^{7} + 1820 x^{6} + 8008 x^{5} + 12870 x^{4} + 8008 x^{3} + 1820 x^{2} + 120 x + 1}{16 \left(x^{7} + 35 x^{6} + 273 x^{5} + 715 x^{4} + 715 x^{3} + 273 x^{2} + 35 x + 1\right)} \, ,
$$
and five iterations produces the enormous expression:
\begin{equation*}
\resizebox{1.0\hsize}{!}{$\frac{x^{16} + 496 x^{15} + 35960 x^{14} + 906192 x^{13} + 10518300 x^{12} + 64512240 x^{11} + 225792840 x^{10} + 471435600 x^{9} + 601080390 x^{8} + 471435600 x^{7} + 225792840 x^{6} + 64512240 x^{5} + 10518300 x^{4} + 906192 x^{3} + 35960 x^{2} + 496 x + 1}{32 \left(x^{15} + 155 x^{14} + 6293 x^{13} + 105183 x^{12} + 876525 x^{11} + 4032015 x^{10} + 10855425 x^{9} + 17678835 x^{8} + 17678835 x^{7} + 10855425 x^{6} + 4032015 x^{5} + 876525 x^{4} + 105183 x^{3} + 6293 x^{2} + 155 x + 1\right)}$} \, .
\end{equation*}
Notice how quickly these grow---in fact, the degree of the polynomials doubles on every iteration!  Now, if we take the symbolic derivatives of these functions using our ordinary calculus rules, and simplify (with the help of the computer), the derivative of one iteration is $\frac{1}{2}$, of two iterations is
$$
\frac{x^{2} + 2 x + 5}{4 \left(x^{2} + 2 x + 1\right)} \, ,
$$
of three iterations is
$$
\frac{x^{6} + 14 x^{5} + 147 x^{4} + 340 x^{3} + 375 x^{2} + 126 x + 21}{8 \left(x^{6} + 14 x^{5} + 63 x^{4} + 100 x^{3} + 63 x^{2} + 14 x + 1\right)} \, ,
$$
of four iterations is
\begin{equation*}
\resizebox{1.0\hsize}{!}{$\frac{x^{14} + 70 x^{13} + 3199 x^{12} + 52364 x^{11} + 438945 x^{10} + 2014506 x^{9} + 5430215 x^{8} + 8836200 x^{7} + 8842635 x^{6} + 5425210 x^{5} + 2017509 x^{4} + 437580 x^{3} + 52819 x^{2} + 3094 x + 85}{16 \left(x^{14} + 70 x^{13} + 1771 x^{12} + 20540 x^{11} + 126009 x^{10} + 440986 x^{9} + 920795 x^{8} + 1173960 x^{7} + 920795 x^{6} + 440986 x^{5} + 126009 x^{4} + 20540 x^{3} + 1771 x^{2} + 70 x + 1\right)}$} \, ,
\end{equation*}
and of five iterations is a monstrosity you can only read by zooming in:
\begin{equation*}
\resizebox{1.0\hsize}{!}{$\frac{x^{30} + 310 x^{29} + 59799 x^{28} + 4851004 x^{27} + 215176549 x^{26} + 5809257090 x^{25} + 102632077611 x^{24} + 1246240871640 x^{23} + 10776333438765 x^{22} + 68124037776390 x^{21} + 321156247784955 x^{20} + 1146261110726340 x^{19} + 3133113888931089 x^{18} + 6614351291211874 x^{17} + 10850143060249839 x^{16} + 13883516068991952 x^{15} + 13883516369532147 x^{14} + 10850142795067314 x^{13} + 6614351497464949 x^{12} + 3133113747810564 x^{11} + 1146261195398655 x^{10} + 321156203432790 x^{9} + 68124057936465 x^{8} + 10776325550040 x^{7} + 1246243501215 x^{6} + 102631341330 x^{5} + 5809427001 x^{4} + 215145084 x^{3} + 4855499 x^{2} + 59334 x + 341}{32 \left(x^{30} + 310 x^{29} + 36611 x^{28} + 2161196 x^{27} + 73961629 x^{26} + 1603620018 x^{25} + 23367042639 x^{24} + 238538538360 x^{23} + 1758637118685 x^{22} + 9579944198310 x^{21} + 39232152623175 x^{20} + 122387258419860 x^{19} + 293729420641881 x^{18} + 546274556891506 x^{17} + 791156255418003 x^{16} + 894836006026128 x^{15} + 791156255418003 x^{14} + 546274556891506 x^{13} + 293729420641881 x^{12} + 122387258419860 x^{11} + 39232152623175 x^{10} + 9579944198310 x^{9} + 1758637118685 x^{8} + 238538538360 x^{7} + 23367042639 x^{6} + 1603620018 x^{5} + 73961629 x^{4} + 2161196 x^{3} + 36611 x^{2} + 310 x + 1\right)}$} \, .
\end{equation*}
This is a terrible way to compute derivatives!  (However, more sophisticated approaches to efficient symbolic differentiation exist, such as the \href{https://www.microsoft.com/en-us/research/publication/the-d-symbolic-differentiation-algorithm/}{``$D^*$'' algorithm}, that avoid explicit giant formulas by exploiting repeated subexpressions.)

To be clear, the dual number approach (absent rounding errors) computes an answer exactly as if it evaluated
these crazy expressions at some particular $x$, but the words ``as if'' are very important here.  As you can see,
we do not form these expressions, let alone evaluate them.  We merely compute results that are equal to the values
we would have gotten if we had.

\pagebreak

\subsection{Automatic Differentiation via Computational Graphs}

Let's now get into automatic differentiation via computational graphs.  For this section, we consider the following simple motivating example.

\begin{example}
    \label{ex:compute-graph}
    Define the following functions: 
    \[
    \begin{cases}
        a(x,y) = \sin x \\
        b(x,y) = \frac{1}{y}\cdot a(x,y) \\
        z(x,y) = b(x,y) + x.
    \end{cases}
    \]
    Compute $\frac{\partial z}{\partial x}$ and $\frac{\partial z}{\partial y}.$
\end{example}

There are a few ways to solve this problem. Firstly, of course, one can compute this symbolically, noting that 
\[
z(x,y) = b(x,y) + x = \frac{1}{y} a(x,y) + x = \frac{\sin x}{y} + x,
\]
which implies 
\[
\frac{\partial z}{\partial x} = \frac{\cos x}{y} + 1 \hspace{.25cm} \text{and} \hspace{.25cm} \frac{\partial z}{\partial y} = -\frac{\sin x}{y^2}.
\]

However, one can also use a Computational Graph (see Figure of Computational Graph below) where the edge from node $A$ to node $B$ is labelled with $\frac{\partial B}{\partial A}$.

\begin{figure}[ht]
\begin{center}
 \begin{tikzpicture}[
            > = stealth, 
            shorten > = 1pt, 
            auto,
            node distance = 3cm, 
            semithick, 
        ]

        \tikzstyle{every state}=[
            draw = black,
            thick,
            fill = white,
            minimum size = 15mm
        ]
        \node[state] (v2) {$a(x,y)$};
        \node[state] (v1) [left of=v2,draw=red] {$x$};
        \node[state] (v3) [right of=v2] {$b(x,y)$};
        \node[state] (v4) [right of=v3,draw=blue] {$z(x,y)$};
        \node[state] (v5) [below of=v1,draw=red] {$y$};

        \draw[->] (v1) edge["$\cos x$"] (v2);
        \draw[->] (v2) edge["$\frac{1}{y}$"] (v3);
        \draw[->] (v3) edge["$1$"] (v4);
        \draw[->] (v1) edge[bend left, "$1$"] (v4);
        \draw[->] (v5) edge[bend right, "$-\frac{a(x,y)}{y^2}$"] (v3);
    \end{tikzpicture}
\caption{A computational graph corresponding to example~\ref{ex:compute-graph}, representing the computation of an \textcolor{blue}{output} $z(x,y)$ from two \textcolor{red}{inputs} $x,y$, with intermediate quantities $a(x,y)$ and $b(x,y)$.  The nodes are labelled by \emph{values}, and edges are labelled with the \emph{derivatives} of the values with respect to the preceding values.}
\end{center}
\end{figure}
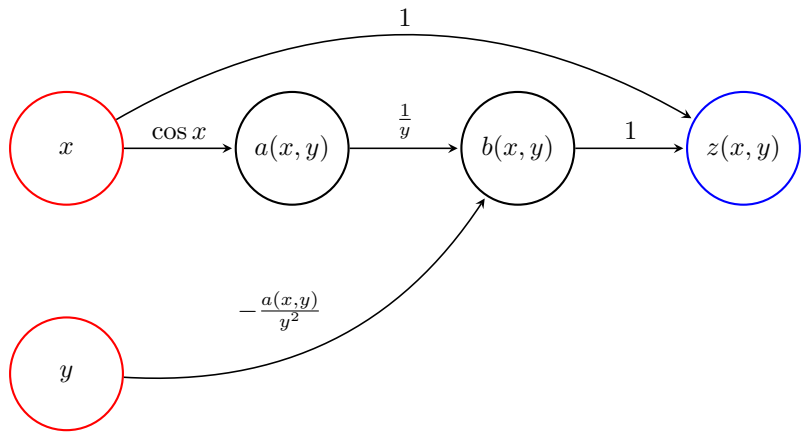
\noindent

Now how do we use this directed acyclic graph (DAG) to find the derivatives? Well one view (called the ``forward view'') is given by following the paths from the inputs to the outputs and (left) multiplying as you go, adding together multiple paths. For instance, following this procedure for paths from $x$ to $z(x,y)$, we have
\[
\frac{\partial z}{\partial x} = 1 \cdot \frac{1}{y} \cdot \cos x + 1 = \frac{\cos x}{y} + 1.
\]
Similarly, for paths from $y$ to $z(x,y)$, we have 
\[
\frac{\partial z}{\partial y} = 1 \cdot \frac{-a(x,y)}{ y^2} = \frac{- \sin x}{y^2},
\]
and if you have numerical derivatives on the edges, this algorithm works. Alternatively, you could follow a reverse view and follow the paths backwards (multiplying right to left), and obtain the same result.
Note that there is nothing magic about these being scalar here-- you could imagine these functions are the type that we are seeing in this class and do the same computations! The only thing that matters here fundamentally is the associativity. However,  when considering vector-valued functions, the order in which you multiply the edge weights is vitally important (as vector/matrix valued functions are not generally commutative).

The graph-theoretic way of thinking about this is to consider ``path products''. A path product is the product of edge weights as you traverse a path. In this way, we are interested in the sum of path products from inputs to outputs to compute derivatives using computational graphs. Clearly, we don't particularly care which order we traverse the paths as long as the \textit{order} we take the product in is correct. In this way, forward and reverse-mode automatic differentiation is not so mysterious.

Let's take a closer view of the implementation of forward-mode automatic differentiation. Suppose we are at a node $A$ during the process of computing the derivative of a computational graph,  as shown in the figure below:

\begin{figure}[ht]
\begin{center}
 \begin{tikzpicture}[
            > = stealth, 
            shorten > = 1pt, 
            auto,
            node distance = 3cm, 
            semithick 
        ]

        \tikzstyle{every state}=[
            draw = black,
            thick,
            fill = white,
            minimum size = 10mm
        ]
        \node[state] (v1) {$A$};
        \node[state] (v2) [right of=v1] {$f(A)$};
        \node[state] (v3) [above left of = v1] {$B_1$};
        \node[state] (v4) [left of=v1] {$B_2$};
        \node[state] (v5) [below left of=v1] {$B_3$};

        \draw[->] (v1) edge["$\frac{\partial f(A)}{\partial A}$"] (v2);
        \draw[->] (v3) -- (v1);
        \draw[->] (v4) -- (v1);
        \draw[->] (v5) -- (v1);
    \end{tikzpicture}
\end{center}
\end{figure}

Suppose we know the path product $P$ of all the edges up to and including the one from $B_2$?
to $A$. Then what is the new path product as we move to the right from $A$?  It is  $f'(A)\cdot P$! So we need a data structure that maps in the following way: 
\[
(\text{value}, \text{path product}) \mapsto (f(\text{value}), f'\cdot \text{path product}).
\]
In some sense, this is another way to look at the Dual Numbers-- taking in our path products and spitting out values. In any case, we overload our program which can easily calculate $f(\text{value})$ and tack-on $f'\cdot \text{(path product)}$.

One might ask how our program starts-- this is how the program works in the ``middle'', but what should our starting value be? Well the only thing it can be for this method to work is $(x, 1)$. Then, at every step you do the following map listed above: 
\[
(\text{value}, \text{path product}) \mapsto (f(\text{value}), f'\cdot \text{path product}),
\]
and at the end we obtain our derivatives.

Now how do we combine arrows? In other words, suppose at the two notes on the LHS we have the values $(a,p)$ and $(b,q)$, as seen in the diagram below:
\begin{figure}[ht]
\begin{center}
 \begin{tikzpicture}[
            > = stealth, 
            shorten > = 1pt, 
            auto,
            node distance = 3cm, 
            semithick 
        ]

        \tikzstyle{every state}=[
            draw = black,
            thick,
            fill = white,
            minimum size = 10mm
        ]
        \node[state] (v1) {$z = f(a,b)$};
        \node[state] (v2) [rectangle, above left of = v1] {$(a,p)$};
        \node[state] (v3) [rectangle, below left of=v1] {$(b,q)$};
        
        \draw[->] (v2) edge["$\frac{\partial z}{\partial a}$"] (v1);
        \draw[->] (v3) edge["$\frac{\partial z}{\partial b}$"] (v1);
    \end{tikzpicture}
\end{center}
\end{figure}
So here, we aren't thinking of $a,b$ as numbers, but as variables. What should the new output value be? We want to add the two path products together, obtaining 
\[
\left(f(a,b), \frac{\partial z}{\partial a} p + \frac{\partial z}{\partial b} q\right).
\]
So really, our overloaded data structure looks like this: 
\begin{figure}[ht]
\begin{center}
 \begin{tikzpicture}[
            > = stealth, 
            shorten > = 1pt, 
            auto,
            node distance = 3cm, 
            semithick 
        ]

        \tikzstyle{every state}=[
            draw = black,
            thick,
            fill = white,
            minimum size = 10mm
        ]
        \node[state] (v1) [rectangle] {$\left(f(a,b), \frac{\partial z}{\partial a} p + \frac{\partial z}{\partial b} q\right)$};
        \node[state] (v2) [rectangle, above left of = v1] {$(a,p)$};
        \node[state] (v3) [rectangle, below left of=v1] {$(b,q)$};
        
        \draw[->] (v2) edge (v1);
        \draw[->] (v3) edge (v1);
    \end{tikzpicture}
\end{center}
\end{figure}

\noindent This diagram of course generalizes if we may many different nodes on the left side of the graph.

If we come up with such a data structure for all of the simple computations (addition/subtraction, multiplication, and division), and if this is all we need for our computer program, then we are set! Here is how we define the structure for addition/subtraction, multiplication, and division.

\newpage \textbf{Addition/Subtraction:} See figure.
\begin{figure}[ht]
\begin{center}
 \begin{tikzpicture}[
            > = stealth, 
            shorten > = 1pt, 
            auto,
            node distance = 3cm, 
            semithick 
        ]

        \tikzstyle{every state}=[
            draw = black,
            thick,
            fill = white,
            minimum size = 10mm
        ]
        \node[state] (v1) [rectangle] {$\left(z = a_1 \pm a_2, \frac{\partial z}{\partial a_1}\cdot 1 + \frac{\partial z}{\partial a_2}\cdot (\pm 1)\right)$};
        \node[state] (v2) [rectangle, above left of = v1] {$(a_1,p= 1)$};
        \node[state] (v3) [rectangle, below left of=v1] {$(a_2, q=\pm 1)$};
        
        \draw[->] (v2) edge (v1);
        \draw[->] (v3) edge (v1);
    \end{tikzpicture}
\end{center}
\caption{Figure of Addition/Subtraction Computational Graph}
\end{figure}
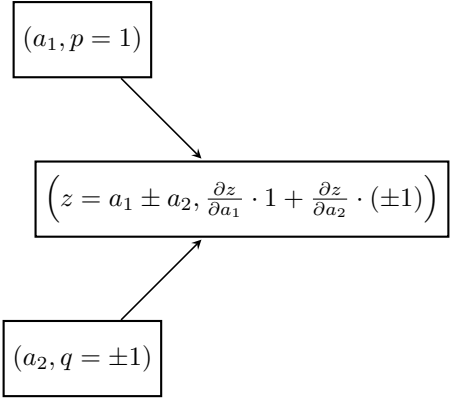

\textbf{Multiplication:} See figure.
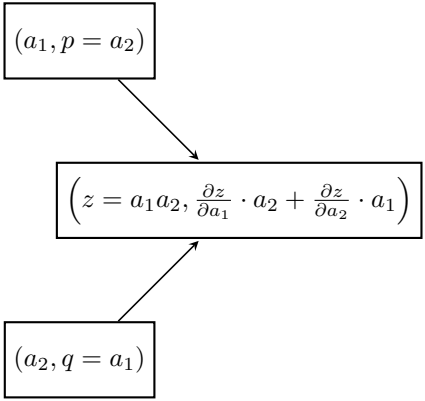
\begin{figure}[ht]
\begin{center}
 \begin{tikzpicture}[
            > = stealth, 
            shorten > = 1pt, 
            auto,
            node distance = 3cm, 
            semithick 
        ]

        \tikzstyle{every state}=[
            draw = black,
            thick,
            fill = white,
            minimum size = 10mm
        ]
        \node[state] (v1) [rectangle] {$\left(z = a_1 a_2, \frac{\partial z}{\partial a_1}\cdot a_2 + \frac{\partial z}{\partial a_2} \cdot a_1\right)$};
        \node[state] (v2) [rectangle, above left of = v1] {$(a_1,p= a_2)$};
        \node[state] (v3) [rectangle, below left of=v1] {$(a_2, q=a_1)$};
        
        \draw[->] (v2) edge (v1);
        \draw[->] (v3) edge (v1);
    \end{tikzpicture}
    \caption{Figure of Multiplication Computational Graph}
\end{center}
\end{figure}

\textbf{Division:} See figure.
\begin{figure}[ht]
\begin{center}
 \begin{tikzpicture}[
            > = stealth, 
            shorten > = 1pt, 
            auto,
            node distance = 3cm, 
            semithick 
        ]

        \tikzstyle{every state}=[
            draw = black,
            thick,
            fill = white,
            minimum size = 10mm
        ]
        \node[state] (v1) [rectangle] {$\left(z = a_1/a_2, \frac{\partial z}{\partial a_1}\cdot \frac{1}{a_2} - \frac{\partial z}{\partial a_2} \cdot \frac{a_1}{a_2^2}\right)$};
        \node[state] (v2) [rectangle, above left of = v1] {$(a_1,p= a_2/a_2^2)$};
        \node[state] (v3) [rectangle, below left of=v1] {$(a_2, q=-a_1/a_2^2)$};
        
        \draw[->] (v2) edge (v1);
        \draw[->] (v3) edge (v1);
    \end{tikzpicture}
    \caption{Figure of Division Computational Graph}
\end{center}
\end{figure}
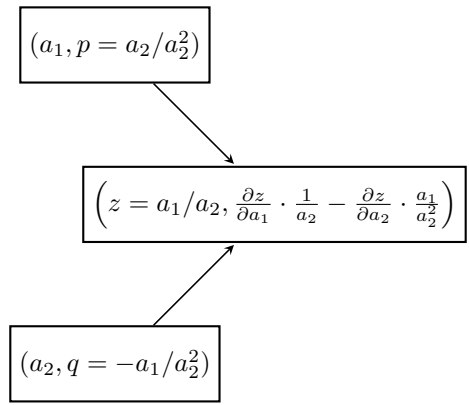

In theory, these three graphs are all we need, and we can use Taylor series expansions for more complicated functions. But in practice, we throw in what the derivatives of more complicated functions are so that we don't waste our time trying to compute something we already know, like the derivative of sine or of a logarithm.

\newpage
\subsubsection{Reverse Mode Automatic Differentiation on Graphs}

When we do reverse mode, we have arrows going the other direction, which we will understand in this section of the notes. In forward mode it was all about ``what do we depend on,'' i.e. computing the derivative on the right hand side of the above diagram using the functions in the nodes on the left. In reverse mode, the question is really ``what are we influenced by?'' or ``what do we influence later?''

When going ``backwards,'' we need know what nodes a given node influences. For instance, given a node A, we want to know the nodes $B_i$ that is influenced by, or depends on, node $A$. So now our diagram looks like this: 
\begin{figure}[ht]
\begin{center}
 \begin{tikzpicture}[
            > = stealth, 
            shorten > = 1pt, 
            auto,
            node distance = 3cm, 
            semithick 
        ]

        \tikzstyle{every state}=[
            draw = black,
            thick,
            fill = white,
            minimum size = 10mm
        ]
        \node[state] (v1) {$\left(a, \frac{\partial z}{\partial a}\right)$};
        \node[state] (v2) [left of=v1] {$(x,\partial z/\partial x)$};
        \node[state] (v3) [above right of = v1] {$\left(b_1, \frac{\partial z}{\partial b_1}\right)$};
        \node[state] (v4) [right of=v1] {$\left(b_2, \frac{\partial z}{\partial b_2}\right)$};
        \node[state] (v5) [below right of=v1] {$\left(b_3, \frac{\partial z}{\partial b_3}\right)$};
        \node[state] (v6) [right of=v4] {$(z,1)$};

        \draw[->] (v1) -- (v3);
        \draw[dotted] (v2) -- (v1);
        \draw[->] (v1) -- (v4);
        \draw[->] (v1) -- (v5); 
        \draw[dotted] (v3) -- (v6);
        \draw[dotted] (v4) -- (v6);
        \draw[dotted] (v5) -- (v6);
    \end{tikzpicture}
\end{center}
\end{figure}

So now, we eventually have a final node $(z,1)$ (far on the right hand side) where everything starts. This time, all of our multiplications take place from right to left as we are in reverse mode. Our goal is to be able to calculate the node $(x,\partial z/\partial x)$. So if we know how to fill in the $\frac{\partial z}{\partial a}$ term, we will be able to go from right to left in these computational graphs (i.e., in reverse mode). In fact, the formula for getting $\frac{\partial z}{\partial a}$ is given by 
\[
\frac{\partial z}{\partial a} = \sum_{i=1}^s \frac{\partial b_i}{\partial a} \frac{\partial z}{\partial b_i}
\]
where the $b_i$s come from the nodes that are influenced by the node $A$. This is again just another chain rule like from calculus, but you can also view this as multiplying the sums of all the weights in the graph influenced by $A$.

\begin{figure}[ht]
\begin{center}
 \begin{tikzpicture}[
            > = stealth, 
            shorten > = 1pt, 
            auto,
            node distance = 3cm, 
            semithick, 
        ]

        \tikzstyle{every state}=[
            draw = black,
            thick,
            fill = white,
            minimum size = 15mm
        ]
        \node[state] (v2) {$p$};
        \node[state] (v1) [above left of=v2] {$x$};
        \node[state] (v3) [right of=v2] {$q$};
        \node[state] (v4) [right of=v3] {$z$};
        \node[state] (v5) [below left of=v2] {$y$};

        \draw[->] (v1) edge["$a$"] (v2);
        \draw[->] (v2) edge["$c$"] (v3);
        \draw[->] (v3) edge["$d$"] (v4);
        \draw[->] (v5) edge["$b$"] (v2);
    \end{tikzpicture}
\end{center}
\end{figure}

Why can reverse mode be more efficient than forward mode? One reason it because it can save data and use it later. Take, for instance, the following sink/source computational graph.

If $x,y$ here are our sources, and $z$ is our sink, we want to compute the sum of products of weights on paths from sources to sinks. If we were using forward mode, we would need to compute the paths $dca$ and $dcb$, which requires four multiplications (and then you would add them together). If we were using reverse mode, we would only need compute $a\underline{cd}$ and $b\underline{cd}$ and sum them; notice reverse mode (since we need only compute $cd$ once), only takes 3 multiplications. In general, this can more efficiently resolve certain types of problems, such as the source/sink one.

\subsection{Forward- vs. Reverse-mode Differentiation}
\label{sec:forward-vs-reverse}

In this section, we briefly summarize the relative benefits and drawbacks of these two approaches to computation of derivatives (whether worked out by hand or using AD software).  From a mathematical point of view, the two approaches are mirror images, but from a computational point of view they are quite different, because computer programs normally proceed ``forwards'' in time from inputs to outputs.

Suppose we are differentiating a function $f: \mathbb{R}^n \mapsto \mathbb{R}^m$, mapping $n$ scalar inputs (an $n$-dimensional input) to $m$ scalar outputs (an $m$-dimensional output).   The first key distinction of forward- vs.~reverse-mode is how the computational cost scales with the number/dimension of inputs and outputs:
\begin{itemize}

\item The cost of forward-mode differentiation (inputs-to-outputs) scales proportional to~$n$, the number of \emph{inputs}.  This is ideal for functions where $n \ll m$ (few inputs, many outputs).

\item The cost of reverse-mode differentiation (outputs-to-inputs) scales proportional to~$m$, the number of \emph{outputs}.  This is ideal for functions where $m \ll n$ (few outputs, many inputs).

\end{itemize}
Before this chapter, we first saw these scalings in Sec.~\ref{sec:cost-matrix-mult}, and again in Sec.~\ref{sec:adjoint-method}; in a future lecture, we'll see it yet again in Sec.~\ref{sec:ODE-sensitivity}.  The case of few outputs is extremely common in large-scale optimization (whether for machine learning, engineering design, or other applications), because then one has many optimization parameters ($n\gg 1$) but only a single output ($m=1$) corresponding to the objective (or ``loss'') function, or sometimes a few outputs corresponding to objective and constraint functions.   Hence, reverse-mode differentiation (``backpropagation'') is the dominant approach for large-scale optimization and applications such as training neural networks.

There are other practical issues worth considering, however:
\begin{itemize}

\item Forward-mode differentiation proceeds in the same order as the computation of the function itself, from inputs to outputs.  This seems to make forward-mode AD easier to implement (e.g. our sample implementation in Sec.~\ref{sec:dual-AD}) and efficient.

\item Reverse-mode differentiation proceeds in the \emph{opposite} direction to ordinary computation.  This makes reverse-mode AD much more complicated to implement, and adds a lot of \emph{storage overhead} to the function computation.  First you evaluate the function from inputs to outputs, but you (or the AD system) keep a \emph{record} (a ``tape'') of all the \emph{intermediate steps} of the computation; then, you run the computation in \emph{reverse} (``play the tape backwards'') to backpropagate the derivatives.

\end{itemize}
As a result of these practical advantages, even for the case of many ($n >1$) inputs and a single ($m = 1$) output, practitioners tell us that they've found forward mode to be more efficient until $n$ becomes sufficiently large (perhaps even until~$n > 100$, depending on the function being differentiated and the AD implementation).  (You may also be interested in the blog post \href{https://www.stochasticlifestyle.com/engineering-trade-offs-in-automatic-differentiation-from-tensorflow-and-pytorch-to-jax-and-julia/}{Engineering Trade-offs in AD} by Chris~Rackauckas, which is mainly about reverse-mode implementations.)

If $n = m$, where neither approach has a scaling advantage, one typically prefers the lower overhead and simplicity of forward-mode differentiation.   This case arises in computing explicit Jacobian matrices for nonlinear root-finding (Sec.~\ref{sec:newton-roots}), or Hessian matrices of second derivatives (Sec.~\ref{sec:hessians}), for which one often uses forward mode\ldots or even a \emph{combination} of forward and reverse modes, as discussed below.  

Of course, forward and reverse are not the only options.  The chain rule is associative, so there are many possible orderings (e.g.~starting from both ends and meeting in the middle, or vice versa).   A difficult\footnote{In fact, extraordinarily difficult: ``NP-complete'' \href{https://dl.acm.org/doi/abs/10.5555/3114201.3114717}{(Naumann, 2006)}.} problem that may often require hybrid schemes is to compute Jacobians (or Hessians) in a minimal number of operations, exploiting any problem-specific structure (e.g.~sparsity: many entries may be zero).  Discussion of this and other AD topics can be found, in vastly greater detail than in these notes, in the book \textit{Evaluating Derivatives} (2nd ed.) by Griewank and Walther~(2008).

\subsubsection{Forward-over-reverse mode: Second derivatives}
\label{sec:forward-over-reverse}

Often, a \emph{combination} of forward- and reverse-mode differentiation is advantageous when computing \emph{second} derivatives, which arise in many practical applications.

\textbf{Hessian computation:}
For example, let us consider a function $f(x): \mathbb{R}^n \to \mathbb{R}$ mapping $n$ inputs $x$ to a single scalar.  The first derivative $f'(x) = (\nabla f)^T$ is best computed by reverse mode if $n \gg 1$ (many inputs).   Now, however, consider the \emph{second} derivative, which is the derivative of $g(x) = \nabla f$, mapping $n$ inputs $x$ to $n$ outputs $\nabla f$.   It should be clear that $g'(x)$ is therefore an $n \times n$ Jacobian matrix, called the \textbf{Hessian} of~$f$, which we will discuss much more generally in Sec.~\ref{sec:hessians}.  Since $g(x)$ has the same number of inputs and outputs, neither forward nor reverse mode has an inherent scaling advantage, so typically forward mode is chosen for $g'$ thanks to its practical simplicity, while still computing $\nabla f$ in reverse-mode.  That is, we compute $\nabla f$ by reverse mode, but then compute $g' = (\nabla f)'$ by applying forward-mode differentiation to the $\nabla f$ algorithm.  This is called a \textbf{forward-over-reverse} algorithm.

An even more clear-cut application of forward-over-reverse differentiation is to \textbf{Hessian--vector products}.  In many applications, it turns out that what is required is only the \emph{product} $(\nabla f)' v$ of the Hessian $(\nabla f)'$ with an arbitrary vector~$v$.  In this case, one can completely avoid computing (or storing) the Hessian matrix explicitly, and incur computational cost proportional only to that of  a single function evaluation $f(x)$.  The trick is to recall (from Sec.~\ref{sec:directional}) that, for \emph{any} function $g$, the linear operation $g'(x)[v]$ is a \emph{directional derivative}, equivalent to a \emph{single-variable} derivative $\frac{\partial}{\partial\alpha} g(x+\alpha v)$ evaluated at $\alpha = 0$.  Here, we simply apply that rule to the function $g(x) = \nabla f$, and obtain the following formula for a Hessian--vector product:
$$
(\nabla f)' v = \evalat{\frac{\partial}{\partial\alpha} \left( \evalat{\nabla f}{x + \alpha v} \right)}{\alpha=0} \, .
$$
Computationally, the inner evaluation of the gradient $\nabla f$ at an arbitrary point $x + \alpha v$ can be accomplished efficiently by a reverse/adjoint/backpropagation algorithm.  In contrast, the \emph{outer} derivative with respect to a \emph{single} input $\alpha$ is best performed by forward-mode differentiation.\footnote{\href{https://jax.readthedocs.io/en/latest/notebooks/autodiff_cookbook.html}{The Autodiff Cookbook}, part of the JAX documentation, discusses this algorithm in a section on Hessian--vector products.  It notes that one could also interchange the $\partial/\partial \alpha$ and $\nabla_x$ derivatives and employ reverse-over-forward mode, but suggests that this is less efficient in practice: ``because forward-mode has less overhead than reverse-mode, and since the outer differentiation operator here has to differentiate a larger computation than the inner one, keeping forward-mode on the outside works best.'' It also presents another alternative: using the identity $(\nabla f)'v = \nabla (v^T \nabla f)$, one can apply reverse-over-reverse mode to take the gradient of $v^T \nabla f$, but this has even more computational overhead.}   
Since the Hessian matrix is symmetric (as discussed in great generality by Sec.~\ref{sec:hessians}), the same algorithm works for \textbf{vector--Hessian products} $v^T (\nabla f)' = [(\nabla f)' v]^T$, a fact that we employ in the next example.

\textbf{Scalar-valued functions of gradients:} There is another common circumstance in which one often combines forward and reverse differentiation, but which can appear somewhat more subtle, and that is in differentiating a scalar-valued function of a gradient of another scalar-valued function.  Consider the following example: 
\todo{Let's talk through whats in chapter 12 first , and then see
I think we need to see how Hessians come up in the book, i've lost track}

\begin{example}

Let $f(x): \mathbb{R}^n \mapsto \mathbb{R}$ be a scalar-valued function of $n \gg 1$ inputs with gradient $\evalat{\nabla f}{x} = f'(x)^T$, and let $g(z): \mathbb{R}^n \mapsto \mathbb{R}$ be another such function with gradient $\evalat{\nabla g}{z} = g'(z)^T$.  Now, consider the scalar-valued function $h(x) = g(\evalat{\nabla f}{x}): \mathbb{R}^n \mapsto \mathbb{R}$ and compute $\evalat{\nabla h}{x} = h'(x)^T$.

Denote $z = \evalat{\nabla f}{x}$. By the chain rule, $h'(x) = g'(z)(\nabla f)'(x)$, but we want to avoid explicitly computing the large $n \times n$ Hessian matrix $(\nabla f)'$.  Instead, as discussed above, we use the fact that such a vector--Hessian product is equivalent (by symmetry of the Hessian) to the transpose of a Hessian--vector product multiplying the Hessian $(\nabla f)'$ with the vector $\nabla g = g'(z)^T$, which is equivalent to a directional derivative:
$$
\evalat{\nabla h}{x} = h'(x)^T =
\evalat{\frac{\partial}{\partial\alpha} \left( \evalat{\nabla f}{x + \alpha \evalat{\nabla g}{z}}\right)}{\alpha = 0} \, ,
$$
involving differentiation with respect to a single scalar $\alpha \in \mathbb{R}$.   As for any Hessian--vector product, therefore, we can evaluate $h$ and $\nabla h$ by:
\begin{enumerate}
    \item Evaluate $h(x)$: evaluate $z = \evalat{\nabla f}{x}$ by reverse mode, and plug it into $g(z)$.
    \item Evaluate $\nabla h$:
    \begin{enumerate}
      \item Evaluate $\evalat{\nabla g}{z}$ by reverse mode.
      \item Implement $\evalat{\nabla f}{x + \alpha \evalat{\nabla g}{z} }$ by reverse mode, and then differentiate with respect to $\alpha$ by \emph{forward} mode, evaluated at $\alpha = 0$.
    \end{enumerate}
\end{enumerate}
This is a ``forward-over-reverse'' algorithm, where forward mode is used efficiently for the single-input derivative with respect to $\alpha \in \mathbb{R}$, combined with reverse mode to differentate with respect to $x,z \in \mathbb{R}^n$.
\end{example}

Example Julia code implementing the above ``forward-over-reverse'' process for just such a $h(x)=g(\nabla f)$ function is given below.  Here, the forward-mode differentiation with respect to $\alpha$ is implemented by the ForwardDiff.jl package discussed in Sec.~\ref{sec:dual-AD}, while the reverse-mode differentiation with respect to $x$ or $z$ is performed by the \href{https://fluxml.ai/Zygote.jl/stable/}{Zygote.jl} package. First, let's import the packages and define simple example functions $f(x) = 1/\Vert x \Vert$ and $g(z) = (\sum_k z_k)^3$, along with the computation of $h$ via Zygote:
\begin{minted}{jlcon}
julia> using ForwardDiff, Zygote, LinearAlgebra
julia> f(x) = 1/norm(x)
julia> g(z) = sum(z)^3
julia> h(x) = g(Zygote.gradient(f, x)[1])
\end{minted}
Now, we'll compute $\nabla h$ by forward-over-reverse:
\begin{minted}{jlcon}
julia> function ∇h(x)
           ∇f(y) = Zygote.gradient(f, y)[1]
           ∇g = Zygote.gradient(g, ∇f(x))[1]
           return ForwardDiff.derivative(α -> ∇f(x + α*∇g), 0)
       end
\end{minted}
We can now plug in some random numbers and compare to a finite-difference check:
\begin{minted}{jlcon}
julia> x = randn(5); δx = randn(5) * 1e-8;

julia> h(x)
-0.005284687528953334

julia> ∇h(x)
5-element Vector{Float64}:
 -0.006779692698531759
  0.007176439898271982
 -0.006610264199241697
 -0.0012162087082746558
  0.007663756720005014

julia> ∇h(x)' * δx    # directional derivative
-3.0273434457397667e-10

julia> h(x+δx) - h(x)  # finite-difference check
-3.0273433933303284e-10
\end{minted}
The finite-difference check matches to about 7~significant digits, which is as much as we can hope for---the forward-over-reverse code works!

\begin{problem}
A common variation on the above procedure, which often appears in machine learning, involves a function $f(x,p) \in \mathbb{R}$ that maps input ``data'' $x \in \mathbb{R}^n$ and ``parameters'' $p \in \mathbb{R}^N$ to a scalar.  Let $\nabla_x f$ and $\nabla_p f$ denote the gradients with respect to $x$ and $p$.

Now, suppose we have a function $g(z): \mathbb{R}^n \mapsto \mathbb{R}$ as before, and define $h(x,p) = g(\evalat{\nabla_x f}{x,p})$.  We want to compute $\nabla_p h = (\partial h / \partial p)^T$, which will involve ``mixed'' derivatives of $f$ with respect to \emph{both} $x$ and $p$.

Show that you can compute $\nabla_p h$ by:
$$
\evalat{\nabla_p h}{x,p} =
\evalat{\frac{\partial}{\partial\alpha} \left( \evalat{\nabla_p f}{x + \alpha \evalat{\nabla g}{z},p}\right)}{\alpha = 0} \, ,
$$
where $z = \evalat{\nabla_x f}{x,p}$.  (Crucially, this avoids ever computing an $n \times N$ mixed-derivative matrix of~$f$.)

Try coming up with simple example functions $f$ and $g$, implementing the above formula by forward-over-reverse in Julia similar to above (forward mode for $\partial/\partial \alpha$ and reverse mode for the $\nabla$'s), and checking your result against a finite-difference approximation.
\end{problem}

\pagebreak

\section{Differentiating ODE solutions}
In this lecture, we will consider the problem of differentiating the \emph{solution} of ordinary differential equations (ODEs) with respect to parameters that appear in the equations and/or initial conditions.  This is as important topic in a surprising number of practical applications, such as evaluating the effect of uncertainties, fitting experimental data, or machine learning (which is increasingly combining ODE models with neural networks).  As in previous lectures, we will find that there are crucial practical distinctions between ``forward'' and ``reverse'' (``adjoint'') techniques for computing these derivatives, depending upon the number of parameters and desired outputs.

Although a basic familiarity with the concept of an ODE will be helpful to readers of this lecture, we will begin with a short review in order to establish our notation and terminology.

The video lecture on this topic for IAP 2023 was given by Dr. Frank Sch\"afer (MIT).  These notes follow the same basic approach, but differ in some minor notational details.

\subsection{Ordinary differential equations (ODEs)}
An \textbf{ordinary differential equation} (\textbf{ODE}) is an equation
for a function $u(t)$ of ``time''\footnote{Of course, the independent variable need not be time, it just needs to be a real scalar.  But in a generic context it is convenient to imagine ODE solutions as evolving in time.} $t\in\mathbb{R}$ in terms of
one or more derivatives, most commonly in the \textbf{first-order}
form 
\[
\frac{du}{dt}=f(u,t)
\]
for some right-hand-side function $f$. Note that $u(t)$ need not
be a scalar function---it could be a column vector $u\in\mathbb{R}^{n}$,
a matrix, or any other differentiable object. One could also write
ODEs in terms of higher derivatives $d^{2}u/dt^{2}$ and so on, but
it turns out that one can write any ODE in terms of first derivatives
alone, simply by making $u$ a vector with more components.\footnote{For example, the second-order ODE $\frac{d^{2}v}{dt^{2}}+\frac{dv}{dt}=h(v,t)$
could be re-written in first-order form by defining $u=\left(\begin{array}{c}
u_{1}\\
u_{2}
\end{array}\right)=\left(\begin{array}{c}
v\\
dv/dt
\end{array}\right)$, in which case $du/dt=f(u,t)$ where $f=\left(\begin{array}{c}
u_{2}\\
h(u_{1},t)-u_{2}
\end{array}\right)$. } To uniquely determine a solution of a first-order ODE, we need some
additional information, typically an \textbf{initial value} $u(0)=u_{0}$
(the value of $u$ at $t=0$), in which case it is called an \textbf{initial-value
problem}. These facts, and many other properties of ODEs, are reviewed
in detail by many textbooks on differential equations, as well as
in classes like 18.03 at MIT.

ODEs are important for a huge variety of applications, because the
behavior of many realistic systems is defined in terms of rates of
change (derivatives). For example, you may recall Newton's laws of
mechanics, in which acceleration (the derivative of velocity) is related
to force (which may be a function of time, position, and/or velocity),
and the solution $u=[\text{position},\text{velocity}]$ of the corresponding
ODE tells us the trajectory of the system. In chemistry, $u$ might
represent the concentrations of one or more reactant molecules, with
the right-hand side $f$ providing reaction rates. In finance, there
are ODE-like models of stock or option prices. \emph{Partial} differential
equations (PDEs) are more complicated versions of the same idea, for
example in which $u(x,t)$ is a function of space $x$ as well as
time $t$ and one has $\frac{\partial u}{\partial t}=f(u,x,t)$ in which $f$ may involve some spatial derivatives of $u$. 

In linear algebra (e.g.~18.06 at MIT), we often consider initial-value
problems for \emph{linear} ODEs of the form $du/dt=Au$ where $u$ is a column
vector and $A$ is a square matrix; if $A$ is a constant matrix (independent
of $t$ or $u$), then the solution $u(t)=e^{At}u(0)$ can be described in terms of a matrix exponential $e^{At}$.
More generally, there are many tricks to find explicit solutions of
various sorts of ODEs (various functions $f$). However, just as one
cannot find explicit formulas for the integrals of most functions,
there is no explicit formula for the solution of \emph{most} ODEs,
and in many practical applications one must resort to approximate
numerical solutions. Fortunately, if you supply a computer program
that can compute $f(u,t)$, there are mature and sophisticated software
libraries\footnote{For a modern and full-featured example, see the DifferentialEquations.jl
suite of ODE solvers in the Julia language.} which can compute $u(t)$ from $u(0)$ for any desired set of times
$t$, to any desired level of accuracy (for example, to 8 significant
digits).

For example, the most basic numerical ODE method computes the solution
at a sequence of times $t_{n}=n\Delta t$ for $n=0,1,2,\ldots$ simply
by approximating $\frac{du}{dt}=f(u,t)$ using the finite difference
$\frac{u(t_{n+1})-u(t_{n})}{\Delta t}\approx f(u(t_{n}),t_{n})$,
giving us the ``explicit'' timestep algorithm:
\[
u(t_{n+1})\approx u(t_{n})+\Delta t\,f(u(t_{n}),t_{n}).
\]
Using this technique, known as ``Euler's method,'' we can march
the solution forward in time: starting from our initial condition
$u_{0}$, we compute $u(t_1) = u(\Delta t)$, then $u(t_2) = u(2\Delta t)$ from $u(\Delta t)$, and so
forth. Of course, this might be rather inaccurate unless we make $\Delta t$
very small, necessitating many timesteps to reach a given time $t$,
and there can arise other subtleties like ``instabilities'' where the
error may accumulate exponentially rapidly with each timestep. It turns
out that Euler's method is mostly obsolete: there are much more sophisticated algorithms that robustly produce accurate solutions with far less
computational cost. However, they all resemble Euler's method in the
conceptual sense: they use evaluations of $f$ and $u$ at a few nearby
times $t$ to ``extrapolate'' $u$ at a subsequent time somehow, and
thus march the solution forwards through time.

Relying on a computer to obtain numerical solutions to ODEs is practically essential, but it can also make ODEs a lot more fun to work with.  If you ever took a class on ODEs, you may remember a lot of tedious labor (tricky integrals, polynomial roots, systems of equations, integrating factors, etc.) to obtain solutions by hand.  Instead, we can focus here on simply setting up the correct ODEs and integrals and trust the computer to do the rest.

\subsection{Sensitivity analysis of ODE solutions}
\label{sec:ODE-sensitivity}

\begin{figure}
    \centering
\includegraphics[width=0.6\textwidth]{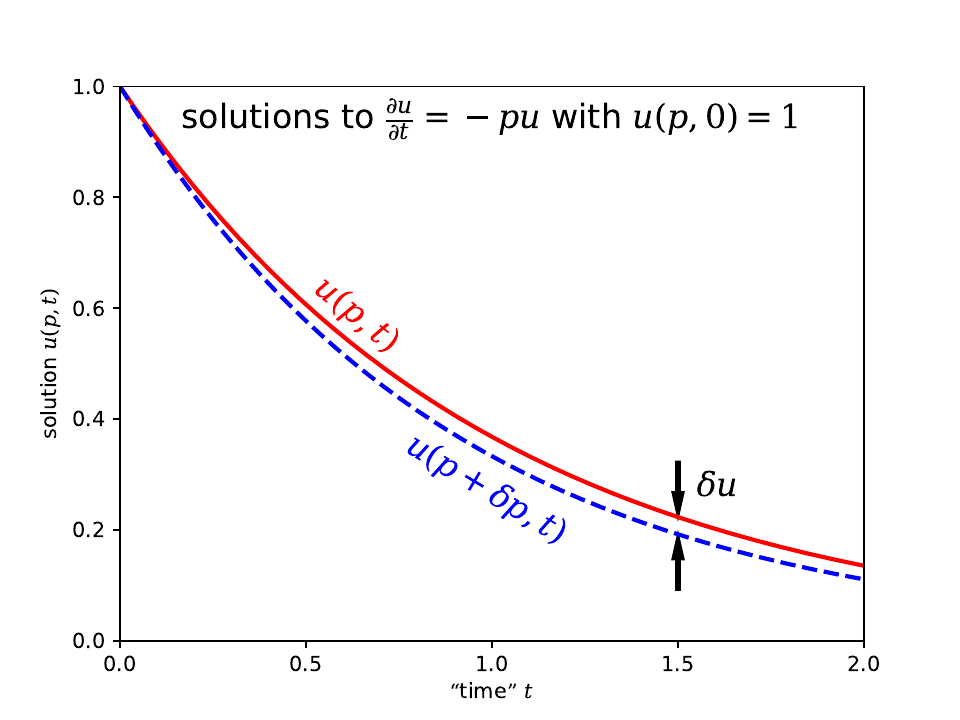}
    \caption{If we have an ordinary differential equation (ODE) $\frac{\partial u}{\partial t} = f(u, p, t)$ whose solution $u(p, t)$ depends on parameters $p$, we would like to know the change $d u = u(p + d p, t) - u(p, t)$ in the solution due to changes in $p$.  Here, we show a simple example $\frac{\partial u}{\partial t} = -pu$, whose solution $u(p,t) = e^{-p t} u(p,0)$ is known analytically, and show the change $\delta u$ from changing $p=1$ to by $\delta p = 0.1$.}
    \label{fig:ode}
\end{figure}

Often, ODEs depend on some additional parameters $p\in\mathbb{R}^{N}$
(or some other vector space). For example, these might be reaction-rate
coefficients in a chemistry problem, the masses of particles in a
mechanics problem, the entries of the matrix $A$ in a linear ODE,
and so on. So, you really have a problem of the form
\[
\frac{\partial u}{\partial t}=f(u,p,t),
\]
where the solution $u(p,t)$ depends both on time $t$ and the parameters
$p$, and in which the initial condition $u(p,0)=u_{0}(p)$ may also
depend on the parameters.

The question is, how can we compute the derivative $\partial u/\partial p$
of the solution with respect to the parameters of the ODE? By this,
as usual, we mean the linear operator that gives the first-order change
in $u$ for a change in $p$, as depicted in Fig.~\ref{fig:ode}:
\[
u(p+dp,t)-u(p,t)=\frac{\partial u}{\partial p}[dp] \qquad \mbox{(an }n\mbox{-component infinitesimal vector)},
\]
where of course $\partial u/\partial p$ (which can be thought of as an $n \times N$ Jacobian matrix) depends on $p$ and $t$.
This kind of question is commonplace. For example, it is important
in:
\begin{itemize}
\item Uncertainty quantification (UQ): if you have some uncertainty in the
parameters of your ODE (for example, you have a chemical reaction
in which the reaction rates are only known experimentally $\pm$ some
measurement errors), the derivative $\partial u/\partial p$ tells
you (to first order, at least) how sensitive your answer is to
each of these uncertainties.
\item Optimization and fitting: often, you want to choose the parameters
$p$ to maximize or minimize some objective (or ``loss'' in machine
learning). For example, if your ODE models some chemical reaction
with unknown reaction rates or other parameters $p$, you might want
to \emph{fit }the parameters $p$ to minimize the difference between
$u(p,t)$ and some experimentally observed concentrations.
\end{itemize}
In the latter case of optimization, you have a \textbf{scalar objective
function} of the solution, since to minimize or maximize something
you need a real number (and $u$ might be a vector). For example,
this could take on one of the following two forms:
\begin{enumerate}
\item A real-valued function $g(u(p,T),T) \in \mathbb{R}$ that depends on the solution $u(p,T)$ at
a particular time $T$. For example, if you have an experimental solution
$u_{*}(t)$ that you are are trying to match at $t=T$, you might
minimize $g(u(p,T),T)=\Vert u(p,T)-u_{*}(T)\Vert^{2}$. 
\item A real-valued function $G(p)=\int_{0}^{T}g(u(p,t),t)dt$ that depends on an average (here scaled by $T$)
over many times $t\in(0,T)$ of our time-dependent $g$. In the example
of fitting experimental data $u_{*}(t)$, minimizing $G(p)=\int_{0}^{T}\Vert u(p,t)-u_{*}(t)\Vert^{2}dt$
corresponds to a least-square fit to minimize the error averaged over
a time $T$ (e.g.~the duration of your experiment).
\end{enumerate}

More generally, you can give more weight to certain times than others
by including a non-negative weight function $w(t)$ in the integral:
$$G_w(p)=\int_0^\infty \underbrace{\|u(p,t)-u_*(t)\|^2}_{g(u(p,t),t)}  \, w(t) \, dt , .$$
The two cases above are simply the choices $w(t) = \delta(t-T)$ (a Dirac delta function) and $w(t) = \begin{cases} 1 & t \le T \\ 0 & \text{otherwise} \end{cases}$ (a~step function), respectively.
As discussed in Problem~\ref{prob:discrete-data}, you can let $w(t)$ be a sum of delta functions to represent data at a sequence of discrete times.

In both cases, since these are scalar-valued functions, for optimization/fitting
one would like to know the gradient $\nabla_{p}g$ or $\nabla_{p}G$,
such that, as usual,  
\[
g(u(p+dp,t),t)-g(u(p,t),t)=\left(\nabla_{p}g\right)^{T}dp
\]
so that $\pm\nabla_{p}g$ is the steepest ascent/descent direction
for maximization/minimization of $g$, respectively.
It is worth emphasizing that gradients (which we only define for scalar-valued functions) have the same shape as their inputs~$p$, so 
$\nabla_p g $ is a vector of length
$N$ (the number of parameters)  that depends on $p$ and $t$.

These are ``just derivatives,'' but probably you can see the difficulty:
if we don't have a formula (explicit solution) for $u(p,t)$, only
some numerical software that can crank out numbers for $u(p,t)$ given
any parameters~$p$ and~$t$, how do we apply differentiation rules to find $\partial u/\partial p$
or $\nabla_{p}g$? Of course, we could use finite differences as in Sec.~\ref{sec:finitedifference}---just
crank through numerical solutions for $p$ and $p+\delta p$ and subtract
them---but that will be quite slow if we want to differentiate with
respect to many parameters ($N\gg1$), not to mention giving potentially
poor accuracy. In fact, people often have \emph{huge} numbers of parameters
inside an ODE that they want to differentiate. Nowadays, our right-hand-side
function $f(u,p,t)$ can even contain a \emph{neural network} (this
is called a ``neural ODE'') with thousands or millions ($N$) of
parameters $p$, and we need all $N$ of these derivatives $\nabla_{p}g$
or $\nabla_{p}G$ to minimize the ``loss'' function $g$ or $G$.
So, not only do we need to find a way to differentiate our ODE solutions
(or scalar functions thereof), but these derivatives must be obtained
\emph{efficiently}. It turns out that there are two ways to do this,
and both of them hinge on the fact that the derivative is obtained
by \emph{solving another ODE}:
\begin{itemize}
\item \textbf{Forward} mode: $\frac{\partial u}{\partial p}$ turns out
to solve \emph{another} ODE that we can integrate with the same numerical
solvers for $u$. This gives us all of the derivatives we could want,
but the drawback is that the ODE for $\frac{\partial u}{\partial p}$
is larger by a factor of $N$ than the original ODE for $u$, so it
is only practical for small $N$ (few parameters).
\item \textbf{Reverse} (``\textbf{adjoint}'') mode: for scalar objectives,
it turns out that $\nabla_{p}g$ or $\nabla_{p}G$ can be computed
by solving a different ODE for an ``adjoint'' solution $v(p,t)$
of the \emph{same size} as $u$, and then computing some simple integrals
involving $u$ (the ``forward'' solution) and $v$.
This has the advantage of giving us all $N$ derivatives with only
about \emph{twice} the cost of solving for $u$, regardless of the
number $N$ of parameters. The disadvantage is that, since it turns
out that $v$ must be integrated ``backwards'' in time (starting
from an ``initial'' condition at $t=T$ and working back to $t=0$)
and depends on $u$, it is necessary to store $u(p,t)$ for all $t\in[0,T]$
(rather than marching $u$ forwards in time and discarding values
from previous times when they are no longer needed), which can require
a vast amount of computer memory for large ODE systems integrated
over long times.
\end{itemize}
We will now consider each of these approaches in more detail.

\subsubsection{Forward sensitivity analysis of ODEs}

Let us start with our ODE $\frac{\partial u}{\partial t}=f(u,p,t)$,
and consider what happens to $u$ for a small change $dp$ in $p$:
\begin{align*}
d\underbrace{\left(\frac{\partial u}{\partial t}\right)}_{=f(u,p,t)} &= \frac{\partial}{\partial t} (du)
=
\frac{\partial}{\partial t}
\left( \frac{\partial u}{\partial p} [ dp ]\right)=\frac{\partial}{\partial t}\left(\frac{\partial u}{\partial p}\right)[dp] \\
&=d(f(u,p,t)) = \left(\frac{\partial f}{\partial u}\frac{\partial u}{\partial p}+\frac{\partial f}{\partial p}\right)[dp],
\end{align*}
where we have used the familiar rule (from multivariable calculus)
of interchanging the order of partial derivatives---a property that
we will re-derive explicitly for our generalized linear-operator derivatives
in our lecture on Hessians and second derivatives. Equating the right-hand sides of the two lines, we see that we
have an ODE 
\[
\boxed{\frac{\partial}{\partial t}\left(\frac{\partial u}{\partial p}\right)=\frac{\partial f}{\partial u}\frac{\partial u}{\partial p}+\frac{\partial f}{\partial p}}
\]
for the derivative $\frac{\partial u}{\partial p}$, whose initial
condition is obtained simply by differentiating the initial condition
$u(p,0)=u_{0}(p)$ for $u$:
\[
\left.\frac{\partial u}{\partial p}\right|_{t=0}=\frac{\partial u_{0}}{\partial p}.
\]
We can therefore plug this into any ODE solver technique (usually
numerical methods, unless we are extremely lucky and can solve this
ODE analytically for a particular $f$) to find $\frac{\partial u}{\partial p}$
at any desired time $t$. Simple, right?

The only thing that might seem a little weird here is the \emph{shape}
of the solution: $\frac{\partial u}{\partial p}$ is a linear operator,
but how can the solution of an ODE be a linear operator? It turns
out that there is nothing wrong with this, but it is helpful to think
about a few examples:
\begin{itemize}
\item If $u,p\in\mathbb{R}$ are scalars (that is, we have a single scalar
ODE with a single scalar parameter), then $\frac{\partial u}{\partial p}$
is just a (time-dependent) number, and our ODE for $\frac{\partial u}{\partial p}$
is an ordinary scalar ODE with an ordinary scalar initial condition.
\item If $u\in\mathbb{R}^{n}$ (a ``system'' of $n$ ODEs) and $p\in\mathbb{R}$
is a scalar, then $\frac{\partial u}{\partial p}\in\mathbb{R}^{n}$
is another column vector and our ODE for $\frac{\partial u}{\partial p}$
is another system of $n$ ODEs. So, we solve two ODEs of the same
size $n$ to obtain $u$ and $\frac{\partial u}{\partial p}$.
\item If $u\in\mathbb{R}^{n}$ (a ``system'' of $n$ ODEs) and $p\in\mathbb{R}^{N}$
is a vector of $N$ parameters, then $\frac{\partial u}{\partial p}\in\mathbb{R}^{n\times N}$
is an $n\times N$ Jacobian \emph{matrix. }Our ODE for $\frac{\partial u}{\partial p}$
is effectivly system of $nN$ ODEs for all the components of this matrix,
with a matrix $\frac{\partial u_{0}}{\partial p}$ of $nN$ initial
conditions! Solving this ``matrix ODE'' with numerical methods poses
no conceptual difficulty, but will generally require about $N$ times
the computational work of solving for $u$, simply because there are
$N$ times as many unknowns. This could be expensive if $N$ is large!
\end{itemize}
This reflects our general observation of forward-mode differentiation:
it is expensive when the number $N$ of ``input'' parameters being
differentiated is large. However, forward mode is straightforward
and, especially for $N\lesssim100$ or so, is often the first method
to try when differentiating ODE solutions. Given $\frac{\partial u}{\partial p}$
, one can then straightforwardly differentiate scalar objectives by
the chain rule: 
\begin{align*}
\left.\nabla_{p}g\right|_{t=T} & =\left.\underbrace{\frac{\partial u}{\partial p}^{T}}_{\text{Jacobian}^{T}}\underbrace{\frac{\partial g}{\partial u}^{T}}_{\text{vector}}\right|_{t=T},\\
\nabla_{p}G & =\int_{0}^{T}\nabla_{p}g\,dt.
\end{align*}
The left-hand side $\nabla_p G$ is gradient of a scalar function of $N$ parameters, and hence the gradient is a vector of $N$ components.  Correspondingly, the right-hand side is an integral of an $N$-component gradient $\nabla_p g$ as well, and the integral of a vector-valued function can be viewed as simply the elementwise integral (the vector of integrals of each component).

\subsubsection{Reverse/adjoint sensitivity analysis of ODEs}

For large $N\gg1$ and scalar objectives $g$ or $G$ (etc.), we can
in principle compute derivatives \emph{much} more efficiently, with
about the same cost as computing $u$, by applying a ``reverse-mode''
or ``adjoint'' approach. In other lectures, we've obtained analogous
reverse-mode methods simply by evaluating the chain rule left-to-right
(outputs-to-inputs) instead of right-to-left. Conceptually, the process
for ODEs is similar,\footnote{This ``left-to-right'' picture can be made very explicit if we imagine
discretizing the ODE into a recurrence, e.g.~via Euler's method for
an arbitrarily small $\Delta t$, as described in the MIT course notes
\emph{\href{https://math.mit.edu/~stevenj/18.336/recurrence2.pdf}{Adjoint methods and sensitivity analysis for recurrence relations}}
by S.~G.~Johnson (2011).} but algebraically the derivation is rather trickier and less direct.
The key thing is that, if possible, we want to avoid computing $\frac{\partial u}{\partial p}$
explicitly, since this could be a prohibitively large Jacobian matrix
if we have many parameters ($p$ is large), especially if we have
many equations ($u$ is large). 

In particular, let's start with our forward-mode sensitivity analysis,
and consider the derivative $G'=(\nabla_{p}G)^{T}$ where $G$ is
the integral of a time-varying objective $g(u,p,t)$ (which we allow
to depend explicitly on $p$ for generality). By the chain rule,
\[
G'=\int_{0}^{T}\left(\frac{\partial g}{\partial p}+\frac{\partial g}{\partial u}\frac{\partial u}{\partial p}\right)dt,
\]
which involves our unwanted factor $\frac{\partial u}{\partial p}$.
To get rid of this, we're going to use a ``weird trick''
\todo{I don't think this needs
to look like a weird trick, but i have
to remember why, EVEN Lagrange
multipliers don't have to be taught
as a weird trick, but I know they often are}
(much like Lagrange
multipliers) of adding \emph{zero} to this equation:
\[
G'=\int_{0}^{T}\left[\left(\frac{\partial g}{\partial p}+\frac{\partial g}{\partial u}\frac{\partial u}{\partial p}\right)+v^{T}\underbrace{\left(\frac{\partial}{\partial t}\left(\frac{\partial u}{\partial p}\right)-\frac{\partial f}{\partial u}\frac{\partial u}{\partial p}-\frac{\partial f}{\partial p}\right)}_{=0}\right]dt
\]
for some function $v(t)$ of the \emph{same shape as u} that multiplies our ``forward-mode'' equation for $\partial u/ \partial p$. (If $u\in\mathbb{R}^{n}$
then $v\in\mathbb{R}^{n}$; more generally, for other vector spaces,
read $v^{T}$ as an inner product with $v$.) The new term $v^T (\cdots)$
is zero because the parenthesized expression is precisely the ODE satisfied by $\frac{\partial u}{\partial p}$,
as obtained in our forward-mode analysis above,  \emph{regardless} of
$v(t)$. This is important because it allows us the freedom to \emph{choose}
$v(t)$ to \emph{cancel} the unwanted $\frac{\partial u}{\partial p}$
term. In particular, if we first \emph{integrate by parts} on the
$v^{T}\frac{\partial}{\partial t}\left(\frac{\partial u}{\partial p}\right)$
term to change it to $-\left(\frac{\partial v}{\partial t}\right)^{T}\frac{\partial u}{\partial p}$
plus a boundary term, then re-group the terms, we find:
\[
G'=\left.v^{T}\frac{\partial u}{\partial p}\right|_{0}^{T}+\int_{0}^{T}\left[\frac{\partial g}{\partial p}-v^{T}\frac{\partial f}{\partial p}+\underbrace{\left(\frac{\partial g}{\partial u}-v^{T}\frac{\partial f}{\partial u}-\left(\frac{\partial v}{\partial t}\right)^{T}\right)}_{\text{want to be zero!}}\frac{\partial u}{\partial p}\right]dt\:.
\]
If we could now set the $(\cdots)$ term to zero, then the unwanted
$\frac{\partial u}{\partial p}$ would vanish from the integral calculation
in $G'$. We can accomplish this by \emph{choosing} $v(t)$ (which
could be \emph{anything} up to now) to satisfy the \textbf{``adjoint''
ODE}:
\[
\boxed{\frac{\partial v}{\partial t}=\left(\frac{\partial g}{\partial u}\right)^{T} - \left(\frac{\partial f}{\partial u}\right)^{T}v}.
\]
What initial condition should we choose
for $v(t)$? Well, we can use this choice to get rid of the boundary
term we obtained above from integration by parts: 
\[
\left.v^{T}\frac{\partial u}{\partial p}\right|_{0}^{T}=v(T)^{T}\underbrace{\left.\frac{\partial u}{\partial p}\right|_{T}}_{\text{unknown}}-v(0)^{T}\underbrace{\frac{\partial u_{0}}{\partial p}}_{\text{known}}.
\]
Here, the unknown $\left.\frac{\partial u}{\partial p}\right|_{T}$
term is a problem---to compute that, we would be forced to go back
to integrating our big $\frac{\partial u}{\partial p}$ ODE from forward
mode. The other term is okay: since the initial condition $u_{0}$
is always given, we should know its dependence on $p$ explicitly
(and we will simply have $\frac{\partial u_{0}}{\partial p}=0$ in
the common case where the initial conditions don't depend on $p$).
To eliminate the $\left.\frac{\partial u}{\partial p}\right|_{T}$
term, therefore, we make the choice 
\[
\boxed{v(T)=0}.
\]
Instead of an \emph{initial} condition, our adjoint ODE has a \textbf{final
condition}. That's no problem for a numerical solver: it just means
that the \textbf{adjoint ODE is integrated }\textbf{\emph{backwards}}\textbf{
in time}, starting from $t=T$ and working down to $t=0$. Once we
have solved the adjoint ODE for $v(t)$, we can plug it into our equation
for $G'$ to obtain our gradient by a simple integral:
\[
\nabla_{p}G=\left(G'\right)^{T}=-\left(\frac{\partial u_{0}}{\partial p}\right)^{T}v(0)+\int_{0}^{T}\left[\left(\frac{\partial g}{\partial p}\right)^{T}-\left(\frac{\partial f}{\partial p}\right)^{T}v\right]dt\:.
\]
(If you want to be fancy, you can compute this $\int_{0}^{T}$ simultaneously
with $v$ itself, by augmenting the adjoint ODE with an additional
set of unknowns and equations representing the $G'$ integrand. But
that's mainly just a computational convenience and doesn't change
anything fundamental about the process.)

The only remaining annoyance is that the adjoint ODE depends on $u(p,t)$
for all $t\in[0,T]$. Normally, if we are solving the ``forward''
ODE for $u(p,t)$ numerically, we can ``march'' the solution $u$
forwards in time and only store the solution at a few of the most
recent timesteps. Since the adjoint ODE starts at $t=T$, however,
we can only start integrating $v$ after we have completed the calculation
of $u$. This requires us to save essentially \emph{all} of our previously
computed $u(p,t)$ values, so that we can evaluate $u$ at arbitrary
times $t\in[0,T]$ during the integration of $v$ (and $G'$). This
can require a lot of computer memory if $u$ is large (e.g.~it could
represent \emph{millions} of grid points from a spatially discretized
PDE, such as in a heat-diffusion problem) and many timesteps $t$
were required. To ameliorate this challenge, a variety of strategies
have been employed, typically centered around ``checkpointing''
techniques in which $u$ is only saved at a subset of times $t$,
and its value at other times is obtained during the $v$ integration
by \emph{re-computing} $u$ as needed (numerically integrating the
ODE starting at the closest ``checkpoint'' time). A detailed discussion
of such techniques lies outside the scope of these notes, however.

\subsection{Example}

Let us illustrate the above techniques with a simple example. Suppose
that we are integrating the scalar ODE 
\[
\frac{\partial u}{\partial t}=f(u,p,t)=p_{1}+p_{2}u+p_{3}u^{2}=p^{T}\left(\begin{array}{c}
1\\
u\\
u^{2}
\end{array}\right)
\]
for an initial condition $u(p,0)=u_{0}=0$ and three parameters $p \in \mathbb{R}^3$. (This is probably simple
enough to solve in closed form, but we won't bother with that here.) We will also consider the scalar function 
\[
G(p)=\int_{0}^{T}\underbrace{\left[u(p,t)-u_{*}(t)\right]^{2}}_{g(u,p,t)}dt
\]
that (for example) we may want to minimize for some given $u_{*}(t)$ (e.g.~experimental
data or some given formula like $u_{*}=t^{3}$), so we are hoping to compute $\nabla_{p}G$.

\subsubsection{Forward mode}

The Jacobian matrix $\frac{\partial u}{\partial p}=\left(\begin{array}{ccc}
\frac{\partial u}{\partial p_{1}} & \frac{\partial u}{\partial p_{2}} & \frac{\partial u}{\partial p_{3}}\end{array}\right)$ is simply a row vector, and satisfies our ``forward-mode'' ODE:
\[
\frac{\partial}{\partial t}\left(\frac{\partial u}{\partial p}\right)=\frac{\partial f}{\partial u}\frac{\partial u}{\partial p}+\frac{\partial f}{\partial p}=\left(p_{2}+2p_{3}u\right)\frac{\partial u}{\partial p}+\left(\begin{array}{ccc}
1 & u & u^{2}\end{array}\right)
\]
for the initial condition $\left.\frac{\partial u}{\partial p}\right|_{t=0}=\frac{\partial u_{0}}{\partial p}=0$.
This is an inhomogeneous system of three coupled \emph{linear} ODEs,
which might look more conventional if we simply transpose both sides:
\[
\frac{\partial}{\partial t}\underbrace{\left(\begin{array}{c}
\frac{\partial u}{\partial p_{1}}\\
\frac{\partial u}{\partial p_{2}}\\
\frac{\partial u}{\partial p_{3}}
\end{array}\right)}_{(\partial u/\partial p)^{T}}=\left(p_{2}+2p_{3}u\right)\left(\begin{array}{c}
\frac{\partial u}{\partial p_{1}}\\
\frac{\partial u}{\partial p_{2}}\\
\frac{\partial u}{\partial p_{3}}
\end{array}\right)+\left(\begin{array}{c}
1\\
u\\
u^{2}
\end{array}\right).
\]
The fact that this depends on our ``forward'' solution $u(p,t)$
makes it not so easy to solve by hand, but a computer can solve it numerically
with no difficulty.  On a computer, we would probably solve for $u$ and $\partial u/ \partial p$\emph{simultaneously} by combining the two ODEs into a single ODE with~4 components:
\[
\frac{\partial}{\partial t}\left(\begin{array}{c}
u \\
(\partial u/\partial p)^{T}
\end{array}\right) =
\left(\begin{array}{c}
p_{1}+p_{2}u+p_{3}u^{2} \\
\left(p_{2}+2p_{3}u\right) (\partial u/\partial p)^{T} +\left(\begin{array}{c}
1\\
u\\
u^{2}
\end{array}\right)
\end{array}\right).
\]
Given $\partial u/\partial p$, we can then plug this into the chain rule for
$G$:
\[
\nabla_{p}G=2\int_{0}^{T}\left[u(p,t)-u_{*}(t)\right]\frac{\partial u}{\partial p}^{T}\,dt
\]
(again, an integral that a computer could evaluate numerically).

\subsubsection{Reverse mode}

In reverse mode, we have an adjoint solution $v(t)\in\mathbb{R}$
(the same shape as $u$) which solves our adjoint equation 
\[
 \frac{\partial v}{dt}=\left(\frac{\partial g}{\partial u}\right)^{T} -\left(\frac{\partial f}{\partial u}\right)^{T}v =2\left[u(p,t)-u_{*}(t)\right] - \left(p_{2}+2p_{3}u\right)v
\]
with a \emph{final} condition $v(T)=0.$ Again, a computer can solve
this numerically without difficulty (given the numerical ``forward''
solution $u$) to find $v(t)$ for $t\in[0,T]$. Finally, our gradient
is the integrated  product:
$$
\nabla_{p}G = -\int_{0}^{T}
\left(
\begin{array}{c} 1 \\ u \\ u^{2} \end{array}
\right)
v\,dt\:.
$$

Another useful exercise is to consider a $G$ that takes the form of a summation:
\begin{problem}
\label{prob:discrete-data}
Suppose that $G(p)$ takes the form of a sum of $K$ terms:
$$
G(p) = \sum_{k=1}^{K} g_k(p,u(p,t_k))
$$
for times $t_k \in (0, T)$ and functions $g_k(p,u)$.  For example, this could arise in least-square fitting of  experimental data $u_*(t_k)$ at $K$ discrete times, with $g_k(u(p,t_k)) = \Vert u_*(t_k) -u(p,t_k)\Vert^2 $ measuring the squared difference between $u(p,t_k)$ and the measured data at time $t_k$.
\begin{enumerate}
    \item Show that such a $G(p)$ can be expressed as a special case of our formulation in this chapter, by defining our function $g(u,t)$ using a sum of Dirac delta functions $\delta(t - t_k)$.
    \item Explain how this affects the adjoint solution $v(t)$: in particular, how the introduction of delta-function terms on the right-hand side of $dv/dt$ causes $v(t)$ to have a sequence of discontinuous jumps.  (In several popular numerical ODE solvers, such discontinuities can be incorporated via discrete-time ``callbacks''.)
    \item Explain how these delta functions  may also introduce a summation into the computation of $\nabla_p G$, but only if $g_k$ depends explicitly on~$p$ (not just via~$u$).
\end{enumerate}
\end{problem}

\subsection{Further reading}

A classic reference on reverse/adjoint differentiation of ODEs (and
generalizations thereof), using notation similar to that used today
(except that the adjoint solution $v$ is denoted $\lambda(t)$, in an
homage to Lagrange multipliers), is Cao~et~al.~(2003) (\url{https://doi.org/10.1137/S1064827501380630}), and a more recent review article is Sapienza~et~al.~(2024) (\url{https://arxiv.org/abs/2406.09699}).
See also the SciMLSensitivity.jl package (\url{https://github.com/SciML/SciMLSensitivity.jl})
for sensitivity analysis with Chris Rackauckas's amazing DifferentialEquations.jl
software suite for numerical solution of ODEs in Julia. There is a
nice 2021 YouTube lecture on adjoint sensitivity of ODEs (\url{https://youtu.be/k6s2G5MZv-I}),
again using a similar notation. A discrete version of this process
arises for recurrence relations, in which case one obtains a reverse-order
``adjoint'' recurrence relation as described in MIT course notes
by S.~G.~Johnson (\url{https://math.mit.edu/~stevenj/18.336/recurrence2.pdf}).

The differentiation methods in this chapter (e.g.~for $\partial u/\partial p$ or $\nabla_p G$) are derived assuming that the ODEs are solved exactly: given the exact ODE for $u$, we derived an exact ODE for the derivative. On a computer, you will solve these forward and adjoint ODEs approximately, and in consequence the resulting derivatives will only be approximately correct (to the tolerance specified by your ODE solver).  This is known as a {\bf differentiate-then-discretize} approach, which has the advantage of simplicity (it is independent of the numerical solution scheme) at the expense of slight inaccuracy (your approximate derivative will not exactly predict the first-order change in your approximate solution $u$).   The alternative is a {\bf discretize-then-differentiate} approach, in which you first approximate (``discretize'') your ODE into a discrete-time recurrence formula, and then \emph{exactly} differentiate the recurrence.  This has the advantage of exactly differentiating your approximate solution, at the expense of complexity (the derivation is specific to your discretization scheme).  Various authors discuss these tradeoffs and their implications, e.g.~in chapter~4 of M.~D.~Gunzburger's \href{https://doi.org/10.1137/1.9780898718720.ch4}{{\it Perspectives in Flow Control and Optimization}} (2002) or in papers like \href{https://dl.acm.org/doi/10.1007/s00158-013-1024-4}{Jensen~et~al.~(2014)}.

\pagebreak

\section{Calculus of Variations}
In this lecture, we will apply our derivative machinery to a new type of input: neither scalars, nor column vectors, nor matrices, but rather the \textbf{inputs will be functions} $u(x)$, which form a perfectly good vector space (and can even have norms and inner products).\footnote{Being fully mathematically rigorous with vector spaces of functions requires a lot of tedious care in specifying a well-behaved set of functions, inserting annoying caveats about functions that differ only at isolated points, and so forth. In this lecture, we will mostly ignore such technicalities---we will implicitly assume that our functions are integrable, differentiable, etcetera, as needed.  The subject of \emph{functional analysis} exists to treat such matters with more care.}  It turns out that there are lots of amazing applications for differentiating with respect to \emph{functions}, and the resulting techniques are sometimes called the ``calculus of variations'' and/or ``Frech{\'{e}}t'' derivatives.

\subsection{Functionals: Mapping functions to scalars}

\begin{example}
For example, consider functions $u(x)$ that map $x\in [0,1] \to u(x) \in \R$.  We may then define the function $f$:
\[
f(u) = \int_0^1 \sin (u(x)) \,\mathrm{d} x.
\]
Such a function, mapping an input \emph{function} $u$ to an output \emph{number}, is sometimes called a ``functional.''  What is $f'$ or $\nabla f$ in this case?
\end{example}

\begin{figure}
    \centering
    \includegraphics[width=0.6\textwidth]{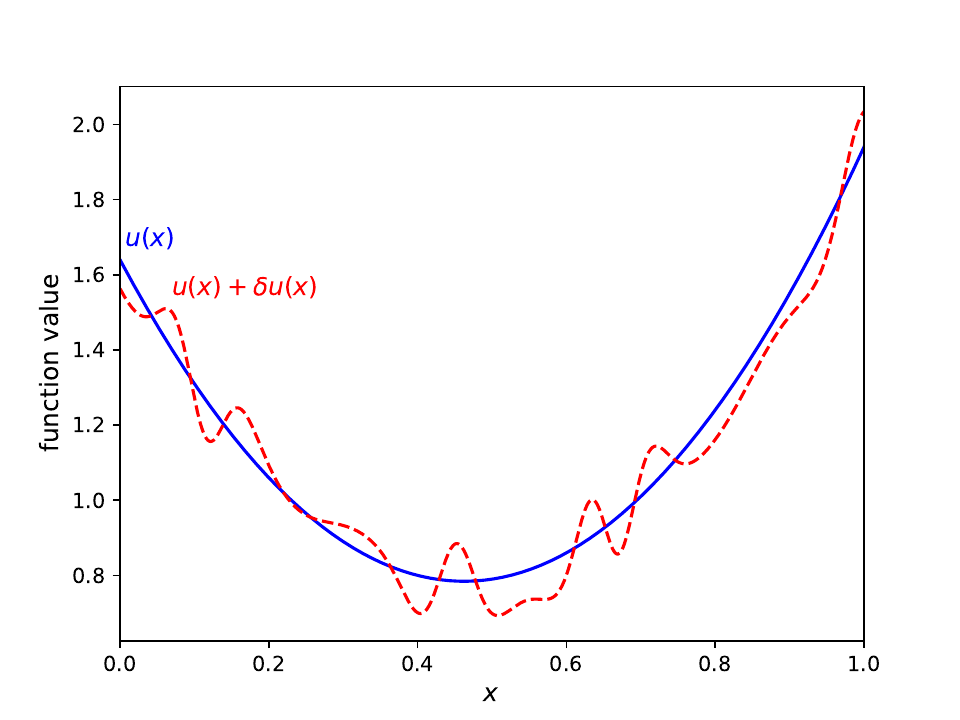}
    \caption{If our $f(u)$'s inputs $u$ are \emph{functions} $u(x)$ (e.g., mapping $[0,1] \mapsto \mathbb{R}$), then the essence of differentiation is linearizing $f$ for small perturbations $\delta u(x)$ that are themselves functions, in the limit where $\delta u(x)$ becomes arbitrarily small.  Here, we show an example of a $u(x)$ and a perturbation $u(x)+\delta u(x)$.}
    \label{fig:du}
\end{figure}

Recall that, given any function $f$, we always define the derivative as a linear operator $f'(u)$ via the equation:
\[
\d f = f (u + \d u ) - f(u) = f'(u) [\d u] \, ,
\]
where now $\d u$ denotes an arbitrary ``small-valued'' \emph{function} $\d u(x)$ that represents a small change in $u(x)$, as depicted in Fig.~\ref{fig:du} for the analogous case of a non-infinitesimal $\delta u(x)$.
Here, we may compute this via linearization of the integrand: 
\begin{align*}
    \d f &= f(u+ \d u) - f(u) \\ 
     &=\int_0^1 \sin (u(x) + \d u(x)) - \sin (u(x)) \, \d x \\
     &= \int_0^1 \cos (u(x)) \, \d u(x) \, \d x = f'(u) [\d u] \, ,
\end{align*}
where in the last step we took $\d u(x)$ to be arbitrarily small\footnote{Technically, it only needs to be small ``almost everywhere'' since jumps that occur only at isolated points don't affect the integral.} so that we could linearize $\sin(u + \d u)$ to first-order in $\d u(x)$.  That's it, we have our derivative $f'(u)$ as a perfectly good linear operation acting on $\d u$!

\subsection{Inner products of functions}

In order to define a gradient $\nabla f$ when studying such ``functionals'' (maps from functions to $\R$), it is natural to ask if there is an inner product on the input space. In fact, there are perfectly good ways to define inner products of functions! Given functions $u(x), v(x)$ defined on $x\in [0,1]$, we could define a ``Euclidean'' inner product:
\[
\langle u, v \rangle = \int_0^1 u(x) v(x) \,\mathrm dx.
\]
Notice that this implies 
\[
\lVert u \rVert := \sqrt{\langle u, u \rangle} = \sqrt{\int_0^1 u(x)^2 dx} \, .
\]

Recall that the gradient $\nabla f$ is \emph{defined} as whatever we take the inner product of $\d u$ with to obtain $\d f$. 
 Therefore, we obtain the gradient as follows: 
\[
\d f = f'(u)[\d u]  = \int_0^1 \cos(u(x)) \,\d u(x)\, \d x  = \langle \nabla f, \d u \rangle \implies \nabla f = \cos (u(x)) \, .
\]
The two infinitesimals $du$ and $dx$ may seem a bit disconcerting, but if this is confusing you can just think of the $du(x)$ as a small non-infinitesimal function $\delta u(x)$ (as in Fig.~\ref{fig:du}) for which we are dropping higher-order terms.

The gradient $\nabla f$ is just another function, $\cos(u(x))$!  As usual, $\nabla f$ has the same ``shape'' as $u$.

\begin{remark}
It might be instructive here to compare the gradient of an integral, above, with a discretized version where the integral is replaced by a sum.   If we have
$$
f(u) = \sum_{k=1}^n \sin(u_k) \Delta x \,
$$
where $\Delta x = 1/n$,
for a vector $u \in \mathbb{R}^n$, related to our previous $u(x)$ by $u_k = u(k\Delta x)$, which can be thought of as a ``rectangle rule'' (or Riemann sum, or Euler) approximation for the integral.  Then,
$$
\nabla_u f = \begin{pmatrix} \cos(u_1) \\ \cos(u_2) \\ \vdots \end{pmatrix} \Delta x  \, .
$$
Why does this discrete version have a $\Delta x$ multiplying the gradient, whereas our continuous version did not?  The reason is that in the continuous version we effectively included the $dx$ in the definition of the inner product $\langle u, v \rangle$ (which was an integral).  In discrete case, the ordinary inner product (used to define the conventional gradient) is just a sum without a $\Delta x$.  However, if we define a \emph{weighted} discrete inner product
$\langle u, v \rangle = \sum_{k=1}^n u_k v_k \Delta x$, then, according to Sec.~\ref{sec:generalvectorspaces}, this changes the definition of the gradient, and in fact will remove the $\Delta x$ term to correspond to the continuous version.
\end{remark}

\subsection{Example: Minimizing arc length}

We now consider a more tricky example with an intuitive geometric interpretation.

\begin{example}
    Let $u$ be a differentiable function on $[0,1]$ and consider the functional 
    \[
    f(u) = \int_0^1 \sqrt{1+ u'(x)^2}\,\d x.
    \]
    Solve for $\nabla f$ when $u(0) = u(1) = 0.$
\end{example}

Geometrically, you learned in first-year calculus that this is simply the \textbf{length of the curve} $u(x)$ from  $x=0$ to $x=1$. To differentiate this, first notice that ordinary single-variable calculus gives us the linearization
\[
\d \left( \sqrt{1 + v^2} \right) = \sqrt{1+ (v + \d v)^2} - \sqrt{1 + v^2} = \left( \sqrt{1 + v^2} \right)' \d v =  \frac{v}{\sqrt{1+ v^2}} \d v \, .
\]
Therefore, 
\begin{align*}
   \d f &= f(u+ \d u) - f(u) \\
   &= \int_0^1 \left( \sqrt{1 + (u+\d u)'^2} - \sqrt{1 + u'^2} \right) \d x \\
   &= \int_0^1 \frac{u'}{\sqrt{1 + u'^2}} \, \d u'  \d x.
\end{align*}
However, this is a linear operator on $\d u'$ and not (directly) on $\d u$.  Abstractly, this is fine, because $\d u'$ is itself a linear operation on $\d u$, so we have $f'(u)[\d u]$ as the composition of two linear operations.  However, it is more revealing to rewrite it explicitly in terms of $\d u$, for example in order to define $\nabla f$. To accomplish this, we can apply \emph{integration by parts} to obtain 
\[
f'(u) [\d u] = \int_0^1 \frac{u'}{\sqrt{1 + u'^2}} \, \d u' \d x = \left. \frac{u'}{\sqrt{1 + u'^2}} \d u \right|_0^1 - \int_0^1 \left(\frac{u'}{\sqrt{1 + u'^2}}\right)' \, \d u \, \d x \, .
\]

Notice that up until now we did not need utilize the ``boundary conditions'' $u(0) = u(1) = 0$ for this calculation. However, if we want to restrict ourselves to such functions $u(x)$, then our perturbation $\d u$ cannot change the endpoint values, i.e.~we must have $\d u(0) = \d u(1) = 0$.
(Geometrically, suppose that we want to find the  $u$ that minimizes arc length between $(0,0)$ and $(1,0)$, so that we need to fix the endpoints.) This implies that the boundary term in the above equation is zero. 
Hence, we have that 
\[
\d f = -\int_0^1 \underbrace{\left(\frac{u'}{\sqrt{1 + u'^2}}\right)'}_{\nabla f} \, \d u \, \d x = \langle \nabla f , \d u \rangle \, .
\]

Furthermore, note that the $u$ that minimizes the functional $f$ has the property that $\left. \nabla f \right|_u = 0$. Therefore, for a $u$ that minimizes the functional $f$ (the \emph{shortest curve}), we must have the following result: 
\begin{align*}
    0 = \nabla f &= -\left(\frac{u'}{\sqrt{1 + u'^2}}\right)' \\
    &= -\frac{u'' \sqrt{1 + u'^2} - u' \frac{u'' u '}{\sqrt{1 + u'^2}}}{ 1 + u'^2} \\
    &= -\frac{u'' ( 1 + u'^2) - u'' u'^2}{(1 + u'^2)^{3/2}} \\
    &= -\frac{u''}{( 1+ u'^2)^{3/2}}.
\end{align*}
Hence, $\nabla f = 0 \implies u''(x) = 0 \implies u(x) = ax + b$ for constants $a,b$; and for these boundary conditions $a=b=0$. In other words, $u$ is the horizontal straight line segment!

Thus, we have recovered the familiar result  that straight line segments in $\R^2$ are the shortest curves between two points!

\begin{remark}
    Notice that the expression $\frac{u''}{(1+ u'^2)^{3/2}}$ is the formula from multivariable calculus for the curvature of the curve defined by $y= u(x)$. It is not a coincidence that the gradient of arc length is the (negative) curvature, and the minimum arc length occurs for zero gradient = zero curvature.
\end{remark}

\subsection{Euler--Lagrange equations}

This style of calculation is part of the subject known as the \textbf{calculus of variations}.  Of course, the final answer in the example above (a straight line) may have been obvious, but a similar approach can be applied to many more interesting problems.  We can generalize the approach as follows:
\begin{example}
    Let $f(u) = \int_a^b F(u, u', x) \, \d x$ where $u$ is a differentiable function on $[a,b]$. Suppose the endpoints of $u$ are fixed (i.e.~its values at $x=a$ and $x=b$ are constants). Let us calculate $\d f$ and $\nabla f$.
\end{example}

We find:
\begin{align*}
    \d f &= f(u + \d u) - f(u) \\
    &= \int_a^b \left(\frac{\partial F}{\partial u} \d u  + \frac{\partial F}{\partial u'} \d u' \right) \d x \\
    &= \underbrace{\frac{\partial F}{\partial u'} \d u \bigr|_a^b}_{= 0} + \int_a^b \left( \frac{\partial F}{\partial u} - \left(\frac{\partial F}{\partial u'}\right)'\right) \, \d u \, \d x \, ,
\end{align*}
where we used the fact that $\d u = 0$ at $a$ or $b$ if the endpoints $u(a)$ and $u(b)$ are fixed. 
Hence, 
\[
\nabla f = \frac{\partial F}{\partial u} - \left(\frac{\partial F}{\partial u'}\right)',
\]
which equals zero at extremum. Notice that $\nabla f = 0$ yields a second-order differential equation in $u$, known as the \href{https://en.wikipedia.org/wiki/Euler%E2%80%93Lagrange_equation}{Euler--Lagrange equations}! 

\begin{remark}
The notation $\partial F / \partial u'$ is a notoriously confusing aspect of the calculus of variations---what does it mean to take the derivative ``with respect to $u'$'' while holding $u$ fixed?   A more explicit, albeit more verbose, way of expressing this is to think of $F(u,v,x)$ as a function of three \emph{unrelated} arguments, for which we only substitute $v=u'$ \emph{after} differentiating with respect to the second argument $v$:
$$
\frac{\partial F}{\partial u'} = \left. \frac{\partial F}{\partial v} \right|_{v=u'} \, .
$$
\end{remark}

There are many wonderful applications of this idea.  For example, search online for information about the ``\href{https://mathworld.wolfram.com/BrachistochroneProblem.html}{brachistochrone problem}'' (animated \href{https://en.wikipedia.org/wiki/Brachistochrone_curve}{here}) and/or the ``\href{https://en.wikipedia.org/wiki/Stationary-action_principle}{principle of least action}''.  Another example is a \href{https://en.wikipedia.org/wiki/Catenary}{catenary} curve, which minimizes the potential energy of a hanging cable.  A classic textbook on the topic is {\it Calculus of Variations} by Gelfand and Fomin.

\pagebreak

\section{Derivatives of Random Functions}
These notes are from a guest lecture by Gaurav Arya in IAP 2023.

\subsection{Introduction}

In this class, we've learned how to take derivatives of all sorts of crazy functions.
Recall one of our first examples:
\begin{equation}
    f(A) = A^2,
\end{equation}
where $A$ is a matrix.
To differentiate this function, we had to go back to the drawing board, and
ask: 
\begin{question}
If we perturb the input slightly, how does the output change?
\label{question:perturb}
\end{question}
To this end, we wrote down something like:
\begin{equation}
    \delta f = (A+\delta A)^2 - A^2 = A (\delta A) + (\delta A) A + \underbrace{(\delta A)^2}_{\text{neglected}}.
\end{equation}
We called $\delta f$ and $\delta A$ \emph{differentials} in the limit where $\delta A$ became arbitrarily small. 
We then had to ask: 
\begin{question}
What terms in the differential can we neglect?
\label{question:neglect}
\end{question}
We decided that $(\delta A)^2$ should be neglected, 
justifying this by the fact
that $(\delta A)^2$ is ``higher-order''. We were left with the derivative operator $\delta A \mapsto 
A(\delta A) + (\delta A) A$: the best possible \emph{linear} approximation
to $f$ in a neighbourhood of $A$. At a high level, the main challenge here was dealing with
complicated input and output spaces: $f$ was matrix-valued, and also matrix-accepting. 
We had to ask ourselves: in this case, what should the notion of a derivative even mean?

In this lecture, we will face a similar challenge, but with an even weirder type of function.
This time, the output of our function will be \emph{random}. Now, we need to 
revisit the same questions. If the output is random,
how can we describe its response to a change in the input?
And how can we form a useful notion of derivative?

\subsection{Stochastic programs}

More precisely, we will consider random, or \emph{stochastic}, functions $X$ with real input
$p \in \mathbb{R}$ and real-valued random-variable output. As a map, we can write $X$ as
\begin{equation}
    p \mapsto X(p),
\end{equation}
where $X(p)$ is a random variable.  (To keep things simple, we'll take $p \in \mathbb{R}$ and $X(p) \in \mathbb{R}$ in this chapter, though of course they could be generalized to other vector spaces as in the other chapters. For now, the randomness is complicated enough to deal with.)

The idea is that we can only \emph{sample} from $X(p)$, according to some distribution of numbers with probabilities that depend upon $p$.   One simple example would be sampling real numbers uniformly (equal probabilities) from the interval $[0,p]$. As a more complicated example, suppose $X(p)$ follows the \emph{exponential distribution} with scale $p$, corresponding to randomly sampled real numbers $x \ge 0$ whose probability decreases proportional to $e^{-x/p}$.  This can be denoted $X(p) \sim \operatorname{Exp}(p)$, and implemented in Julia by:

\begin{minted}{jlcon}
julia> using Distributions

julia> sample_X(p) = rand(Exponential(p))
sample_X (generic function with 1 method)
\end{minted}
We can take a few samples:
\begin{minted}{jlcon}
julia> sample_X(10.0)
1.7849785709142214

julia> sample_X(10.0)
4.435847397169775

julia> sample_X(10.0)
0.6823343897949835

julia> mean(sample_X(10.0) for i = 1:10^9) # mean = p
9.999930348291866
\end{minted}

If our program gives a different output each time, what could a useful notion of derivative be?
Before we try to answer this, let's ask \emph{why} we might want to take a derivative.
The answer is that we may be very interested in \emph{statistical properties} of random functions,
i.e.~values that can be expressed using \emph{averages}.  Even if a function is stochastic, its \emph{average} (``expected value''), assuming the average exists, can be a deterministic function of its parameters that has a conventional derivative.

So, why not take the average \emph{first}, and then take the ordinary derivative of this average?  This simple approach works for very basic stochastic functions (e.g.~the exponential distribution above has expected value $p$, with derivative $1$), but runs into practical difficulties for more complicated distributions (as are commonly implemented by large computer programs working with random numbers).
\begin{remark}
    It is often much easier to produce an ``unbiased estimate'' $X(p)$ of a statistical quantity than to compute it exactly. (Here, an unbiased estimate means that $X(p)$ averages out to our statistical quantity of interest.)
\end{remark}

For example, in deep learning, the ``variational autoencoder'' (VAE) is a very common architecture
that is inherently stochastic. It is easy to get a stochastic \emph{unbiased estimate} of the loss function by running a
random simulation $X(p)$: the loss function $L(p)$ is then the ``average'' value of $X(p)$, denoted by the \emph{expected value} $\EE[X(p)]$.
However, computing the loss $L(p)$ exactly would require integrating over all possible outcomes, which usually is  impractical.
Now, to train the VAE, we also need to differentiate $L(p)$, i.e.~differentiate $\EE[X(p)]$ with respect to $p$!

Perhaps more intuitive examples can be found in the physical sciences, where randomness may be baked into
your model of a physical process.
In this case, it's hard to get around the fact that you need to deal with stochasticity! 
For example, you may have two particles that interact with an \emph{average} rate of $r$.
But in reality, the times when these interactions actually occur follow a stochastic process. 
(In fact, the time until the first interaction might be exponentially distributed, with scale $1/r$.)
And if you want to (e.g.) fit the parameters of your stochastic model to real-world data, it's once again very useful to have derivatives.

If we can't compute our statistical quantity of interest exactly, it seems unreasonable to assume we
can compute its derivative exactly. However, we could hope to stochastically \emph{estimate} its derivative.
That is, if $X(p)$ represents the full program that produces an unbiased estimate of our statistical quantity, 
here's one property we'd definitely like our notion of derivative to have: 
we should be able to construct from it an unbiased gradient estimator\footnote{For more discussion of these concepts, see (e.g.) the review article ``Monte Carlo gradient estimation in machine learning'' (2020) by Mohamed \textit{et al.} (\url{https://arxiv.org/abs/1906.10652}).} $X'(p)$ satisfying
\begin{equation}
    \EE[X'(p)] = \EE[X(p)]' = \frac{\partial \EE[X(p)]}{\partial p}.
\end{equation}
Of course, there are infinitely many such estimators. For example, given any estimator $X'(p)$ we can add any other random variable that has zero average without changing the  expectation value. But in practice there are two additional considerations: (1) we want $X'(p)$ to be easy to compute/sample (about as easy as $X(p)$), and (2) we want the \emph{variance} (the ``spread'') of $X'(p)$ to be small enough that we don't need too many samples to estimate its average accurately (hopefully no worse than estimating $\EE[X(p)]$).

\subsection{Stochastic differentials and the reparameterization trick}

Let's begin by answering our first question (Question \ref{question:perturb}): how does $X(p)$ respond to a change in $p$?
Let us consider a specific $p$ and write down a \emph{stochastic differential}, taking a small but non-infinitesimal $\delta p$ to avoid thinking about infinitesimals for now:
\begin{equation}
    \delta X(p) = X(p + \delta p) - X(p),
\end{equation}
where $\delta p$ represents an arbitrary small change in $p$.
What sort of object is $\delta X(p)$?

Since we're subtracting two random variables, it ought to itself be a random variable.
However, $\delta X(p)$ is still not fully specified! We have only specified 
the marginal distributions of $X(p)$ and $X(p+\delta p)$: to be able to subtract the two,
we need to know their \emph{joint distribution}.

One possibility is to treat $X(p)$ and $X(p + \delta p)$ as independent. 
This means that $\delta X(p)$ would be constructed as the difference of independent samples.
Let's see how samples from $\delta X(p)$ would look like in this case! 
\begin{minted}{jlcon}
julia> sample_X(p) = rand(Exponential(p))
sample_X (generic function with 1 method)

julia> sample_δX(p, δp) = sample_X(p + δp) - sample_X(p)
sample_δX (generic function with 1 method)

julia> p = 10; δp = 1e-5;

julia> sample_δX(p, δp)
-26.000938718875904

julia> sample_δX(p, δp)
-2.6157162001718092

julia> sample_δX(p, δp)
6.352622554495474

julia> sample_δX(p, δp)
-9.53215951927184

julia> sample_δX(p, δp)
1.2232268930932104
\end{minted}
We can observe something a bit worrying: even for a very tiny $\delta p$ (we chose $\delta p = 10^{-5}$),
$\delta X(p)$ is still fairly large: essentially as large as the original random variables.
This is not good news if we want to construct a derivative from $\delta X(p)$: we would rather see its magnitude 
getting smaller and smaller with $\delta p$, like in the non-stochastic case.  Computationally, this will make it very difficult to determine $\EE[X(p)]'$ by averaging \texttt{sample\_δX(p, δp) / δp} over many samples: we'll need a huge number of samples because the \emph{variance}, the ``spread'' of random values, is huge for small~$\delta p$.

Let's try a different approach. It is natural to think of $X(p)$ for all $p$ as forming a \emph{family} of random variables, all 
defined on the same \emph{probability space}.
A probability space, with some simplification, is a sample space $\Omega$, with a probability distribution $\mathbb{P}$ defined 
on the sample space.
From this point of view, each $X(p)$ can be expressed as a function $\Omega \to \mathbb{R}$. 
To sample from a particular $X(p)$, we can imagine drawing a random $\omega$ from $\Omega$ according to $\PP$, and then
plugging this into $X(p)$, i.e.~computing $X(p)(\omega)$.  (Computationally, this is how most distributions are actually implemented: you start with a primitive pseudo-random number generator for a very simple distribution,\footnote{Most computer hardware cannot generate numbers that are actually random, only numbers that \emph{seem} random, called ``pseudo-random'' numbers.  The design of these random-seeming numeric sequences is a subtle subject, steeped in number theory, with a long history of mistakes.  A famous ironic quotation in this field is (Robert Coveyou, 1970): ``Random number generation is too important to be left to chance.''} e.g.~drawing values $\omega$ uniformly from $\Omega = [0,1)$, and then you build other distributions on top of this by transforming $\omega$ somehow.) Intuitively, all of the ``randomness'' resides in the probability space, and crucially 
$\mathbb{P}$ does not depend on $p$: as $p$ varies, $X(p)$ just becomes a different 
\emph{deterministic} map on $\Omega$.

The crux here is that all the $X(p)$ functions now depend on a shared source of randomness: the random draw of $\omega$.
This means that $X(p)$ and $X(p+\delta p)$ have a nontrivial joint distribution: what does it look like?

For concreteness, let's study our exponential random variable $X(p) \sim \operatorname{Exp}(p)$ from above.  Using the ``inversion sampling'' parameterization,
it is possible to choose $\Omega$ to be $[0,1)$ and $\mathbb{P}$ to be the uniform distribution over
$\Omega$; for any distribution, we can construct $X(p)$ to be a corresponding nondecreasing function over $\Omega$ (given by the inverse of $X(p)$'s cumulative probability distribution). 
Applied to $X(p) \sim \operatorname{Exp}(p)$, the inversion method gives $X(p)(\omega) = -p \log{(1-\omega)}$.  This is implemented below, and is a theoretically equivalent way
of sampling $X(p)$ compared with the opaque \texttt{rand(Exponential(p))} function we used above:
\begin{minted}{julia}
julia> sample_X2(p, ω) = -p * log(1 - ω)
sample_X2 (generic function with 1 method)

julia> # rand() samples a uniform random number in [0,1)
julia> sample_X2(p) = sample_X2(p, rand()) 
sample_X2 (generic function with 2 methods)

julia> sample_X2(10.0)
8.380816941818618

julia> sample_X2(10.0)
2.073939134369733

julia> sample_X2(10.0)
29.94586208847568

julia> sample_X2(10.0)
23.91658360124792
\end{minted}
Okay, so what does our joint distribution look like?
\begin{figure}[t]
    \centering
    \includegraphics[width=0.4\textwidth]{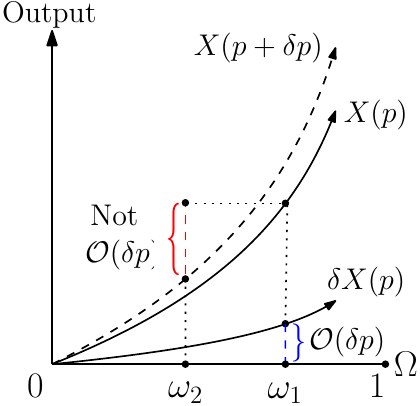}
    \caption{For $X(p) \sim \operatorname{Exp}(p)$ parameterized via the inversion method, we can write $X(p)$, $X(p+\delta p)$, and
    $\delta X(p)$ as functions from $\Omega = [0,1] \to \mathbb{R}$, defined on a probability space with
    $\mathbb{P} = \operatorname{Unif}(0,1)$.}
    \label{fig:exp}
\end{figure}
As shown in Figure~\ref{fig:exp}, we can plot $X(p)$ and $X(p+\delta p)$ as functions over $\Omega$.
To sample the two of them jointly, we use the \emph{same} choice of $\omega$:
thus, $\delta X(p)$ can be formed by subtracting the two functions \emph{pointwise} at each $\Omega$. 
Ultimately, $\delta X(p)$ is itself a random variable over the same probability space, sampled in the same way: 
we pick a random $\omega$ according to $\mathbb{P}$, and evaluate $\delta X(p)(\omega)$, using the function $\delta X(p)$ depicted above.
Our first approach with independent samples is depicted in red in Figure~\ref{fig:exp}, while our second approach is in blue.
We can now see the flaw of the independent-samples approach: the $\mathcal{O}(1)$-sized ``noise'' from the independent samples washes out
the $\mathcal{O}(\delta p)$-sized ``signal''.

What about our second question (Question \ref{question:neglect}): how can actually take the limit of $\delta p \to 0$ and compute the derivative?
The idea is to differentiate $\delta X(p)$ at each fixed sample $\omega \in \Omega$. 
In probability theory terms, we take the limit of random variables $\delta X(p) / \delta p$
as $\delta p \to 0$:
\begin{equation}
    X'(p) = \lim_{\delta p \to 0} \frac{\delta X(p)}{\delta p}.
\end{equation}
For $X(p) \sim \operatorname{Exp}(p)$ parameterized via the inversion method, we get:
\begin{equation}
    X'(p)(\omega) = \lim_{\delta p \to 0} \frac{-\delta p\log{(1-\omega)}}{\delta p} = -\log{(1-\omega)}.
\end{equation}
Once again, $X'(p)$ is a random variable over the same probability space.
The claim is that $X'(p)$ is the notion of derivative we were looking for! Indeed, $X'(p)$ is itself in fact
a valid gradient estimator:
\begin{equation}
    \EE[X'(p)] = \EE\left[ \lim_{\delta p \to 0} \frac{\delta X(p)}{\delta p}\right] \stackrel{?}{=} 
\lim_{\delta p \to 0} \frac{\EE[\delta X(p)]}{\delta p} = 
    \frac{\partial \EE[X(p)]}{\partial p}.
    \label{eq:exchange} 
 \end{equation} 
Rigorously, one needs to justify the interchange of limit and expectation in the above.
In this chapter, however, we will be content with a crude empirical justification:
\begin{minted}{julia}
julia> X′(p, ω) = -log(1 - ω)
X′ (generic function with 1 method)

julia> X′(p) = X′(p, rand())
X′ (generic function with 2 methods)

julia> mean(X′(10.0) for i in 1:10000)
1.011689946421105
\end{minted}
So $X'(p)$ does indeed average to 1, which makes sense since the expectation of $\operatorname{Exp}(p)$ is $p$, which has derivative 1 for any choice of $p$. However, 
the crux is that this notion of derivative 
also works for more complicated random variables that can be formed via \emph{composition} of simple ones such as an exponential random variable. In fact, it turns out to obey the same chain rule as usual!


Let's demonstrate this. Using the dual numbers introduced in Chapter~\ref{sec:AD}, we can differentiate the expectation of the square of a sample 
from an exponential distribution \emph{without} having an analytic expression for this quantity.
(The expression for $X'$ we derived is already implemented as a dual-number rule in Julia by the \texttt{ForwardDiff.jl} package.)
The primal and dual values of the outputted dual number are samples from the joint distribution of $(X(p), X'(p))$.
\begin{minted}{julia}
julia> using Distributions, ForwardDiff: Dual

julia> sample_X(p) = rand(Exponential(p))^2
sample_X (generic function with 1 method)

julia> sample_X(Dual(10.0, 1.0)) # sample a single dual number!
Dual{Nothing}(153.74964559529033,30.749929119058066)

julia> # obtain the derivative!
julia> mean(sample_X(Dual(10.0, 1.0)).partials[1] for i in 1:10000)
40.016569793650525
\end{minted}
Using the ``reparameterization trick'' to form a gradient estimator, as we have done here, is a fairly old idea. It is also called 
the ``pathwise'' gradient estimator. Recently,
it has become very popular in machine learning due to its use in VAEs [e.g.~Kingma~\& Welling (2013): \url{https://arxiv.org/abs/1312.6114}], and lots of resources
can be found online on it. Since composition simply works by the usual chain rule, it also works in reverse mode, and
can differentiate functions far more complicated than the one above!

\subsection{Handling discrete randomness}

So far we have only considered a continuous random variable. Let's see how the picture changes for a discrete random 
variable! Let's take a simple Bernoulli variable $X(p) \sim \operatorname{Ber}(p)$, which is 1 with probability $p$ and
0 with probability $1-p$.
\begin{minted}{jlcon}
julia> sample_X(p) = rand(Bernoulli(p))
sample_X (generic function with 1 method)

julia> p = 0.5
0.6

julia> sample_X(δp) # produces false/true, equivalent to 0/1
true


julia> sample_X(δp)
false

julia> sample_X(δp)
true
\end{minted}
The parameterization of a Bernoulli variable is shown in Figure 2.
Using the inversion method once again, the parameterization of a Bernoulli variable looks like a step function:
for $\omega < 1-p$, $X(p)(\omega) = 0$, while for $\omega \geq 1-p$, $X(p)(\omega) = 1$.

Now, what happens when we perturb $p$? Let's imagine perturbing $p$ by a positive amount $\delta p$.
As shown in Figure 2, something qualitatively very different has happened here. At nearly every $\omega$ except a small region
of probability $\delta p$, the output does not change. Thus, the quantity $X'(p)$ we defined in the previous subsection 
(which, strictly speaking, was defined by an "almost-sure" limit that neglects regions of probability 0)
is 0 at every $\omega$: after all, for every $\omega$, there exists small enough $\delta p$ such that $\delta X(p)(\omega) = 0$.

\begin{figure}[t]
    \centering
    \begin{subfigure}[b]{0.4\textwidth}
        \centering 
        \includegraphics[height=4.5cm]{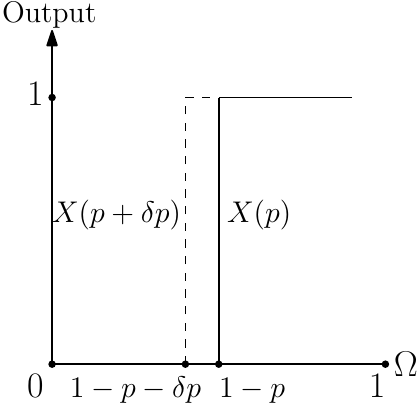}
        \label{fig:ber1} 
    \end{subfigure} 
    \begin{subfigure}[b]{0.4\textwidth}
        \centering
        \includegraphics[height=4.5cm]{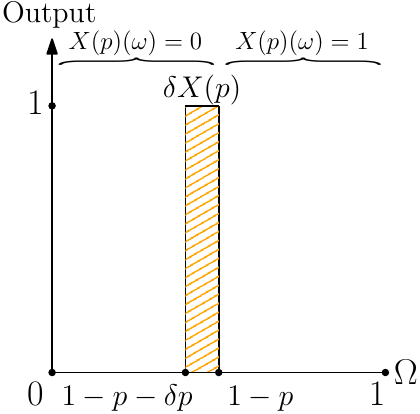}
        \label{fig:ber2}
    \end{subfigure}
    \caption{For $X(p) \sim \operatorname{Ber}(p)$ parameterized via the inversion method, 
    plots of $X(p)$, $X(p+\delta p)$, and $\delta X(p)$ as functions $\Omega: [0, 1] \to \mathbb{R}$. }
\end{figure}

However, there is certainly an important derivative contribution to consider here. 
The expectation of a Bernoulli is $p$, so we would expect the derivative to be 1: but $\EE[X'(p)] = \EE[0] = 0$. What has gone wrong is that, although $\delta X(p)$ is 0 with tiny probability, the value of $\delta X(p)$
on this region of tiny probability is 1, which is \emph{large}. In particular, it does not approach 0 as $\delta p$
approaches 0. Thus, to develop a notion of derivative of $X(p)$, we need to somehow capture these large jumps with ``infinitesimal''
probability. 

A recent (2022) publication (\url{https://arxiv.org/abs/2210.08572}) by the author of this chapter (Gaurav Arya), together with Frank~Sch\"afer, Moritz~Schauer, and Chris~Rackauckas, worked to extend the above ideas to develop a notion of 
``stochastic derivative'' for discrete randomness, implemented by a software package called \texttt{StochasticAD.jl} that performs automatic differentiation of such stochastic processes. 
It generalizes the idea of dual numbers to stochastic \emph{triples}, which include a third component to capture exactly these large jumps. For example, the stochastic triple of a Bernoulli variable might look like:
\begin{minted}{jlcon}
julia> using StochasticAD, Distributions
julia> f(p) = rand(Bernoulli(p)) # 1 with probability p, 0 otherwise
julia> stochastic_triple(f, 0.5) # Feeds 0.5 + δp into f
StochasticTriple of Int64:
0 + 0ε + (1 with probability 2.0ε)
\end{minted}
Here, $\delta p$ is denoted by $\upepsilon$,  imagined to be an ``infinitesimal unit'', so that the above triple indicates a flip from 0 to 1 with probability that has derivative $2$. 

However, many aspects of these problems are still difficult, and there are a lot of improvements awaiting future developments! If you're
interested in reading more, you may be interested in the paper and our package linked above, as well as the 2020 review article by Mohamed \textit{et~al.} (\url{https://arxiv.org/abs/1906.10652}), which is a great survey
of the field of gradient estimation in general.

At the end of class, we considered a differentiable random walk example with \texttt{StochasticAD.jl}. Here it is!





\begin{minted}{jlcon}
julia> using Distributions, StochasticAD

julia> function X(p)
           n = 0
           for i in 1:100
               n += rand(Bernoulli(p * (1 - (n+i)/200)))
           end
           return n
       end
X (generic function with 1 method)

julia> mean(X(0.5) for _ in 1:10000) # calculate E[X(p)] at p = 0.5
32.6956

julia> st = stochastic_triple(X, 0.5) # sample a single stochastic triple at p = 0.5
StochasticTriple of Int64:
32 + 0δp + (1 with probability 74.17635818221052δp)

julia> derivative_contribution(st) # derivative estimate produced by this triple
74.17635818221052

julia> # compute d/dp of E[X(p)] by taking many samples
julia> mean(derivative_contribution(stochastic_triple(f, 0.5)) for i in 1:10000)
56.65142976168479
\end{minted}

\pagebreak

\section[Second Derivatives, Bilinear Maps,
and Hessian Matrices]{Second Derivatives, Bilinear Maps, \\
and Hessian Matrices}
\label{sec:hessians}
In this chapter, we apply the principles of this course to \emph{second} derivatives, which are conceptually just derivatives of derivatives but turn out to have many interesting ramifications.  We begin with a (probably) familiar case of scalar-valued functions from multi-variable calculus, in which the second derivative is simply a matrix called the \emph{Hessian}.  Subsequently, however, we will show that similar principles can be applied to more complicated input and output spaces, generalizing to a notion of $f''$ as a \emph{symmetric bilinear map}.

\subsection{Hessian matrices of scalar-valued functions}
\label{sec:Hessian-scalar}

Recall that for a function $f(x) \in \mathbb{R}$ that maps column vectors $x \in \mathbb{R}^n$ to scalars ($f: \mathbb{R}^n \mapsto \mathbb{R}$), the first derivative $f'$ can be expressed in terms of the familiar gradient $\nabla f = (f')^T$ of multivariable calculus:
$$
\nabla f = \begin{pmatrix} \frac{\partial f}{\partial x_1} \\ \frac{\partial f}{\partial x_2} \\ \vdots \\ \frac{\partial f}{\partial x_n}
\end{pmatrix} \, .
$$
If we think of $\nabla f$ as a new (generally nonlinear) function mapping $x \in \mathbb{R}^n \mapsto \nabla f \in \mathbb{R}^n$, then its derivative is an $n \times n$ Jacobian matrix (a linear operator mapping vectors to vectors), which we can write down explicitly in terms of \emph{second} derivatives of $f$:
$$
  (\nabla f)' =
\begin{pmatrix}
   \frac{\partial^2 f}{ \partial x_1^2} & \cdots & \frac{\partial^2 f}{\partial x_n \partial x_1} \\
        \vdots  &\ddots & \vdots \\
        \frac{\partial^2 f}{\partial x_1 \partial x_n} & \cdots & \frac{\partial^2 f}{\partial x_n \partial x_n}
    \end{pmatrix} = H \, .
$$
This matrix, denoted here by $H$, is known as the \textbf{Hessian} of $f$, which has entries:
\[
H_{i,j} = \frac{\partial^2 f}{ \partial x_j \partial x_i} = \frac{\partial^2 f}{ \partial x_i \partial x_j} = H_{j,i}  \, .
\]
The fact that you can take partial derivatives in either order is a familiar fact from multivariable calculus (sometimes called the ``symmetry of mixed derivatives'' or ``equality of mixed partials''), and means that the Hessian is a \emph{symmetric matrix} $H = H^T$.  (We will later see that such symmetries arise very generally from the construction of second derivatives.)

\begin{example}
For $x \in \mathbb{R}^2$ and the function $f(x) = \sin(x_1) + x_1^2 x_2^3$, its gradient is
$$
\nabla f = \begin{pmatrix} \cos(x_1) + 2x_1 x_2^3 \\
3 x_1^2 x_2^2
\end{pmatrix} \, ,
$$
and its Hessian is
$$
H = (\nabla f)' = \begin{pmatrix} -\sin(x_1) + 2x_2^3 & 6 x_1 x_2^2 \\
6 x_1 x_2^2 & 6 x_1^2 x_2 \end{pmatrix} = H^T \, .
$$
\end{example}

If we think of the Hessian as the Jacobian of $\nabla f$, this tells us that $H \, dx$ predicts the change in $\nabla f$ to first order:
$$
d(\nabla f) = \evalat{\nabla f}{x+dx} - \evalat{\nabla f}{x} = H \, dx \, .
$$
Note that $\evalat{\nabla f}{x+dx}$ means $\nabla f$ evaluated at $x+dx$, which is very different from $df = (\nabla f)^T dx$, where we act $f'(x)=(\nabla f)^T$ \emph{on} $dx$.

Instead of thinking of $H$ of predicting the \emph{first}-order change in $\nabla f$, however, we can also think of it as predicting the \emph{second}-order change in $f$, a \textbf{quadratic approximation} (which could be viewed as the first three terms in a multidimensional Taylor series):
$$
f(x+\delta x) = f(x) + (\nabla f)^T \, \delta x + \frac{1}{2} \delta x^T \, H \, \delta x + o(\Vert \delta x \Vert^2) \, ,
$$
where both $\nabla f$ and $H$ are evaluated at $x$, and we have switched from an infinitesimal $dx$  to a finite change $\delta x$ so that we emphasize the viewpoint of an \emph{approximation} where terms higher than second-order in $\Vert \delta x \Vert$ are dropped.   You can derive this in a variety of ways, e.g. by taking the derivative of both sides with respect to $\delta x$ to reproduce $\evalat{\nabla f}{x+\delta x} = \evalat{\nabla f}{x} + H \, \delta x + o(\delta x)$: a quadratic approximation for $f$ corresponds to a linear approximation for $\nabla f$.   Related to this equation, another useful (and arguably more fundamental) relation that we can derive (and \emph{will} derive much more generally below) is:
$$
dx^T H dx' = f(x + \d x + \d x') + f(x) - f(x + \d x) - f(x + \d x') = f''(x)[dx,dx'] \, 
$$
where $dx$ and $dx'$ are two independent ``infinitesimal'' directions and we have dropped terms of higher than second order.   This formula is very suggestive, because it uses $H$ to map \emph{two} vectors into a \emph{scalar}, which we will generalize below into the idea of a \textbf{bilinear map} $f''(x)[dx,dx']$.  This formula is also obviously symmetric with respect to interchange of $dx$ and $dx'$ --- $f''(x)[dx,dx'] = f''(x)[dx',dx]$ --- which will lead us once again to the symmetry $H=H^T$ below.

\begin{remark}
Consider the Hessian matrix versus other Jacobian matrices.  The Hessian matrix expresses the \emph{second} derivative of a scalar-valued multivariate function, and is always square and symmetric.  A Jacobian matrix, in general, expresses the \emph{first} derivative of a vector-valued multivariate function, may be non-square, and is rarely symmetric.  (However, the Hessian matrix \emph{is} the Jacobian of the $\nabla f$ function!)
\end{remark}

\subsection{General second derivatives: Bilinear maps}

Recall, as we have been doing throughout this class, that we define the derivative of a function $f$ by a linearization of its change $\d f$ for a small (``infinitesimal'') change $\d x$ in the input: 
\[
\d f = f(x + \d x) - f(x) = f'(x) [\d x] \, ,
\]
implicitly dropping higher-order terms.
If we similarly consider the second derivative $f''$ as simply the same process applied to $f'$ instead of $f$, we  obtain the following formula, which is easy to write down but will take some thought to interpret: 
\[
d f' = f'(x + \d x') - f'(x) = f''(x) [\d x'].
\]
(Notation: $\d x'$ is not some kind of derivative of $\d x$; the prime simply denotes a \emph{different}
arbitrary small change in~$x$.) What kind of ``thing'' is $df'$?   Let's consider a simple concrete example:

\begin{example}
Consider the following function $f(x): \mathbb{R}^2 \mapsto \mathbb{R}^2$ mapping two-component vectors $x \in \mathbb{R}^2 $ to two-component vectors $f(x) \in \mathbb{R}^2$:
$$
f(x) = \begin{pmatrix} x_1^2 \sin(x_2) \\ 5x_1 - x_2^3
\end{pmatrix} \, .
$$
Its first derivative is described by a $2\times 2$ Jacobian matrix:
$$
f'(x) = \begin{pmatrix}
2x_1 \sin(x_2) & x_1^2 \cos(x_2) \\
5 & -3x_2^2
\end{pmatrix}
$$
that maps a small change $dx$ in the input vector $x$ to the corresponding small change $df = f'(x)dx$ in the output vector $f$.

What is $df' = f''(x)[dx']$?  It must take a small change $dx' = (dx_1', dx_2')$ in $x$ and return the first-order change $df' = f'(x+dx')-f'(x)$ in our Jacobian matrix $f'$.  If we simply take the differential of each entry of our Jacobian (a function from vectors $x$ to matrices $f'$), we find:
$$
df' = 
\begin{pmatrix}
2\,dx_1' \sin(x_2) + 2x_1 \cos(x_2) \,dx_2' & 2 x_1 \, dx_1' \cos(x_2) - x_1^2 \sin(x_2)\, dx_2' \\
0 & -6x_2\,dx_2'
\end{pmatrix}
= f''(x)[dx']
$$
That is, $df'$ is a $2\times 2$ matrix of ``infinitesimal'' entries, of the same shape as $f'$.

From this viewpoint, $f''(x)$ is a linear operator acting on vectors $dx'$ and outputting $2\times 2$ matrices $f''(x)[dx']$, but this is one of the many cases where it is easier to write down the linear operator as a ``rule'' than as a ``thing'' like a matrix.  The ``thing'' would have to either be some kind of ``three-dimensional matrix'' or we would have to ``vectorize'' $f'$ into a ``column vector'' of 4 entries in order to write its $4\times 4$ Jacobian, as in Sec.~\ref{sec:kronecker} (which can obscure the underlying structure of the problem).

Furthermore, since this $df'$ is a linear operator (a matrix), we can act it on \emph{another} vector $dx = (dx_1, dx_2)$ to obtain:
$$
df' \begin{pmatrix} dx_1 \\ dx_2 \end{pmatrix} =
\begin{pmatrix}
2 \sin(x_2) \,dx_1'\,dx_1 + 2x_1 \cos(x_2) (\,dx_2'\,dx_1 + \, dx_1' \, dx_2) - x_1^2 \sin(x_2)\, dx_2' \, dx_2 \\
-6x_2\,dx_2' \,dx_2
\end{pmatrix} = f''(x)[dx'][dx] \, .
$$
Notice that \emph{this} result, which we will call $f''(x)[dx', dx]$ below, is the ``same shape'' as $f(x)$ (a 2-component vector).  Moreover, it doesn't change if we swap $dx$ and $dx'$: $f''(x)[dx', dx] = f''(x)[dx, dx']$, a key symmetry of the second derivative that we will discuss further below.
\end{example} 

$df'$ is an (infinitesimal) object of the same ``shape'' as $f'(x)$, \emph{not} $f(x)$.  Here, $f'$ is a linear operator,
so its change $\d f'$ must \emph{also} be an (infinitesimal) linear operator (a ``small change'' in a linear operator)
that we can therefore act on an arbitrary
$\d x$ (or $\delta x$), in the form: 
\[
df'[dx] = f''(x) [\d x'][\d x] := f''(x) [\d x', \d x] \, ,
\]
where we combine the two brackets for brevity. This final result $f''(x) [\d x', \d x]$ is the same type of object (vector) as the original output $f(x)$.  This implies that $f''(x)$ is a \textit{bilinear map}: acting on \textbf{two} vectors, and linear in either vector taken individually.  (We will see shortly that the ordering of $dx$ and $dx'$ doesn't matter: $f''(x)[dx',dx]=f''(x)[dx,dx']$.)  

More precisely, we have the following.
\begin{definition}[Bilinear Map]
    Let $U,V,W$ be a vector spaces, not necessarily the same. Then, a \textbf{bilinear map} is a function $B:U\times V \to W$, mapping a $u \in U$ and $v \in V$ to $B[u,v] \in W$, such that we have linearity in both arguments:
    \[
    \begin{cases}
        B[u, \alpha v_1 + \beta v_2] = \alpha B[u,v_1] + \beta B[u,v_2] \\
        B[\alpha u_1 + \beta u_2, v] = \alpha B[u_1,v] + \beta B[u_2,v]
    \end{cases}
    \]
     for any scalars $\alpha, \beta$,

     If $W = \mathbb{R}$, i.e.~the output is a scalar, then it is called a \textbf{bilinear form}.
\end{definition}

\begin{figure}
\begin{center}
\includegraphics[width=0.8\textwidth]{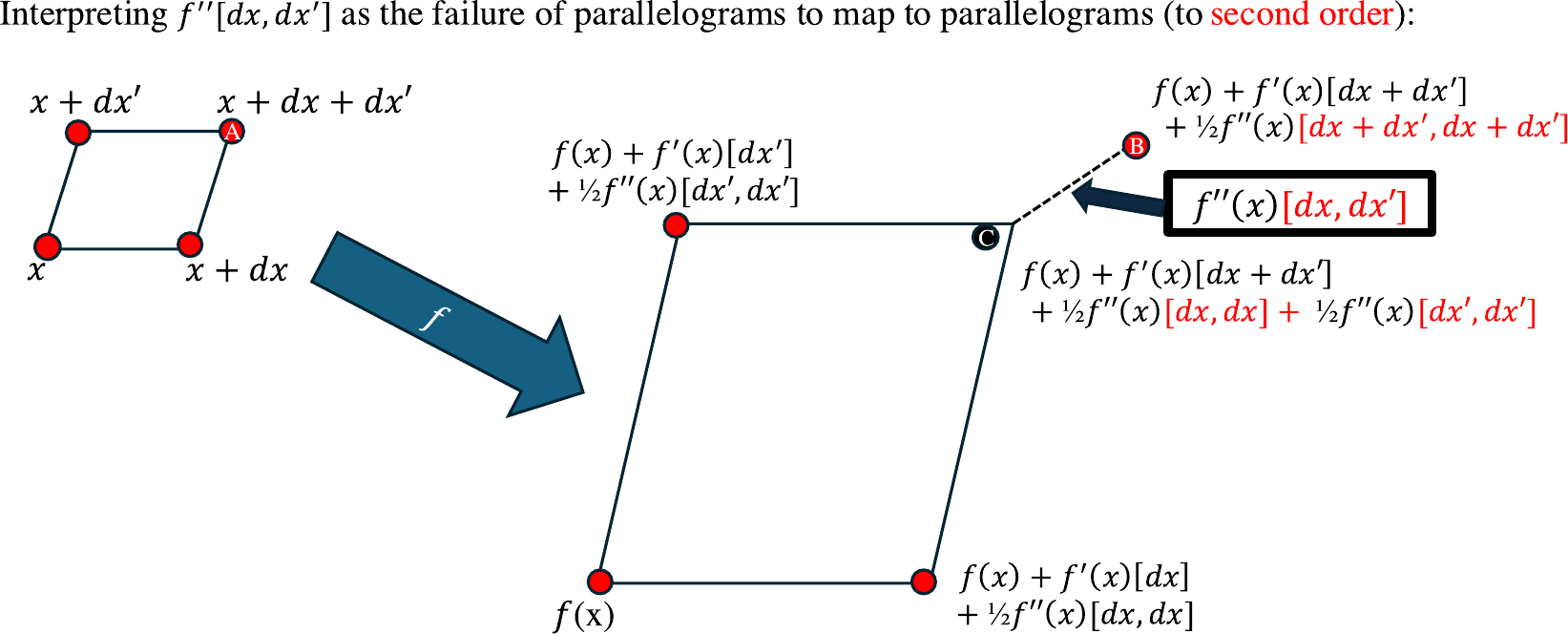}
\end{center}
\caption{Geometric interpretation of $f''(x)[dx,dx']$:
To first order, a function $f$ maps parallelograms to parallelograms.  To second order, however it ``opens'' parallelograms: The deviation from point $B$ (the image of $A$) from point $C$ (the completion of the parallelogram) is the second
derivative $f''(x)[dx,dx']$.  
The symmetry of $f''$ as a bilinear form can be traced back  geometrically 
to the mirror symmetry of the input parallelogram across its diagonal
from $x$ to point~A. \label{fig:second-deriv}
}
\end{figure}

Note that in general, even if $U = V$ (the two inputs $u,v$ are the ``same type'' of vector) we may have $B[u,v] \neq B[v,u]$, but in the case of $f''$ we have something very special that happens. In particular, we can show that $f''(x)$ is a \textit{symmetric bilinear map}, meaning 
\[
f''(x) [\d x', \d x] = f''(x) [\d x, \d x']
\]
for any $\d x$ and $\d x'$.  Why? Because, applying the definition of $f''$ as giving the change in $f'$ from $\d x'$, and then the definition of $f'$ as giving the change in $f$ from $\d x$, we can  re-order terms to obtain: 
\begin{align*}
    f''(x) [\d x', \d x] &= f'(x + 
    \d x') [\d x] - f'(x) [\d x] \\
    &= \left(f(x+ \underbrace{\d x' + \d x}_{= \d x + \d x'}) - f(x+ \d x')\right) - \left(f(x+ \d x) - f(x)\right) \\ 
    &= \boxed{f(x + \d x + \d x') + f(x) - f(x + \d x) - f(x + \d x')} \\
    &= \left(f(x + \d x + \d x')  - f(x + \d x)\right) - \left(f(x + \d x') - f(x)\right) \\
    &= f'(x + 
    \d x) [\d x'] - f'(x) [\d x'] \\
    &= f''(x) [\d x, \d x'] \, 
\end{align*}
where we've boxed the middle formula for $f''$ which emphasizes its symmetry in a natural way.
(The basic reason why this works is that the ``$+$'' operation is always \emph{commutative} for any vector space. A geometric interpretation is depicted in Fig.~\ref{fig:second-deriv}.)

\begin{example}
    Let's review the familiar example from multivariable calculus, $f: \R^n \to \R$. That is, $f(x)$ is a scalar-valued function of a column vector $x\in \R^n$.  What is $f''$?
\end{example}
Recall that 
\[
f'(x) = (\nabla f)^T \implies f'(x) [\d x] = \text{scalar } \d f = (\nabla f)^T \d x. 
\]
Similarly, 
\begin{align*}
    f''(x) [\d x', \d x] &= \text{scalar from two vectors, linear in both} \\
    &= \d x'^T H \d x \, ,
\end{align*}
where $H$ must be exactly the $n \times n$ matrix  \textbf{Hessian matrix} introduced in Sec.~\ref{sec:Hessian-scalar}, since an expression like $\d x'^T H \d x$ is the most general possible bilinear form mapping two vectors to a scalar.  Moreover, since we now know that $f''$ is always a \emph{symmetric} bilinear form, we must have:
\begin{align*}
f''(x) [\d x', \d x] &=\d x'^T H \d x \\
&= f''(x) [\d x, \d x'] = \d x^T H \d x' =
     (\d x^T H \d x')^T  \qquad (\mathrm{scalar} = \mathrm{scalar}^T)\\
    &= \d x'^T H^T \d x 
\end{align*}
for all $dx$ and $dx'$.
This implies that $H = H^T$: the Hessian matrix is symmetric.  As discussed in Sec.~\ref{sec:Hessian-scalar}, we already knew this from multi-variable calculus.  Now, however, this ``equality of mixed partial derivatives'' is simply a special case of $f''$ being a symmetric bilinear map.

As an example, let's consider a special case of the general formula above:

\begin{example}
    Let $f(x) = x^T A x$ for $x \in \R^n$ and $A$ an $n \times n$ matrix. As above, $f(x) \in \R$ (scalar outputs). Compute $f''$.
\end{example}

The computation is fairly straightforward. Firstly, we have that 
\[
f' = (\nabla f)^T = x^T(A + A^T).
\]
This implies that $\nabla f = (A + A^T) x$, a linear function of $x$. Hence, the Jacobian of $\nabla f$ is the Hessian
$f'' = H = A + A^T$. Furthermore, note that this implies 
\begin{align*}
    f(x) &= x^T A x = (x^T A x)^T  \qquad (\mathrm{scalar} = \mathrm{scalar}^T)\\
    &= x^T A^T x \\
    &= \frac{1}{2} (x^T A x + x^T A^T x) = \frac{1}{2} x^T (A + A^T) x\\
    &= \frac{1}{2} x^T H x = \frac{1}{2} f''[x,x] \, ,
\end{align*}
which will turn out to be a special case of the quadratic approximations of Sec.~\ref{sec:Hessian-quadratic} (exact in this example since $f(x)=x^T A x$ is quadratic to start with).

\begin{example}
    Let $f(A) = \det A$ for $A$ an $n\times n$ matrix. Express $f''(A)$ as a rule for $f''(A)[dA,dA']$ in terms of $\d A$ and $\d A'$.
\end{example}

From lecture 3, we have the first derivative 
\[
f'(A) [\d A] = \d f = \det (A) \tr(A^{-1} \d A).
\]
Now, we want to compute the change $d'(df) = d'(f'(A)[dA]) = f'(A+dA')[dA] - f'(A)[dA]$ in this formula, i.e.~the differential (denoted $d'$) where we change $A$ by $dA'$ while treating $dA$ as a constant:
\begin{align*}
    f''(A) [\d A, \d A'] &= \d ' (\det A \tr (A^{-1} \d A)) \\
    &= \det A \tr(A^{-1} \d A') \tr(A^{-1} \d A) - \det A \tr(A^{-1} \,\d A' A^{-1}\, \d A) \\
    &= f''(A) [\d A', \d A]
\end{align*}
where the last line (symmetry) can be derived
explicitly by the cyclic property of the trace (although of course it must be true for any $f''$).  Although $f''$ here is a perfectly good bilinear form acting on matrices $dA,dA'$, it is not very natural to express $f''$ as a ``Hessian matrix.''

If we really wanted to express $f''$ in terms of an explicit Hessian matrix, we could use the the ``vectorization'' approach of Sec.~\ref{sec:kronecker}.
Let us consider, for example, the term $\tr(A^{-1} \,\d A' A^{-1}\, \d A)$ using Kronecker products (Sec.~\ref{sec:kronecker}).  In general, for matrices $X,Y,B,C$:
$$
(\vecm{X})^T (B \otimes C) \vecm{Y} = (\vecm{X})^T \vecm{(CYB^T)} = \tr(X^T CYB^T) = \tr(B^T X^T CY) \, ,
$$
recalling that $(\vecm{X})^T \vecm{Y} = \tr(X^T Y)$ is the Frobenius inner product (Sec.~\ref{sec:generalvectorspaces}).  Thus,
$$
\tr(A^{-1} \,\d A' A^{-1}\, \d A)
= \vecm{(\d A'^T)}^T (A^{-T} \otimes A^{-1}) \vecm{(\d A)} \, .
$$
This is still not quite in the form we want for a Hessian matrix, however, because it involves $\vecm{(\d A'^T)}$ rather than $\vecm{(\d A')}$ (the two vectors are related by a permutation matrix, sometimes called a ``commutation'' matrix). Completing this calculation would be a nice exercise in mastery of Kronecker products, but getting an explicit Hessian seems like a lot of algebra for a result of dubious utility!

\subsection{Generalized quadratic approximation}
\label{sec:Hessian-quadratic}

So how do we ultimately think about $f''$? We know that $f'$ is the linearization/linear approximation of $f(x)$, i.e. 
\[
f(x + \delta x ) = f(x) + f'(x) [\delta x] + o(\lVert \delta x\rVert).
\] 
Now, just as we did for the simple case of Hessian matrices in Sec.~\ref{sec:Hessian-scalar} above, we can use $f''$ to form a \textit{quadratic approximation} of $f(x)$. In particular, one can show that 
\[
f(x+ \delta x) = f(x) + f'(x) [\delta x] + \frac{1}{2} f''(x) [\delta x, \delta x] +  o (\lVert \delta x\rVert^2).
\]
Note that the $\frac{1}{2}$ factor is just as in the Taylor series. To derive this, simply plug the quadratic approximation into 
\[
f''(x) [\d x, \d x'] = f(x + \d x + \d x') + f(x) - f(x + \d x) - f(x + \d x').
\]
and check that the right-hand side reproduces $f''(x)$.   (Note how $dx$ and $dx'$ appear symmetrically in this formula, which reflects the symmetry of $f''$.)

\subsection{Hessians and optimization}

Many important applications of second derivatives, Hessians, and quadratic approximations arise in optimization: minimization (or maximization) of functions $f(x)$.\footnote{Much of machine learning uses only variations on gradient descent, without incorporating Hessian information except implicitly via ``momentum'' terms.  Partly this can be explained by the fact that optimization problems in ML are typically solved only to low accuracy, often have nonsmooth/stochastic aspects, rarely involve nonlinear constraints, and are often very high-dimensional.  This is only a small corner of the wider universe of computational optimization!}

\subsubsection{Newton-like methods}

When searching for a local minimum (or maximum) of a complicated function $f(x)$, a common procedure is to approximate $f(x+\delta x)$ by a simpler ``model'' function for small $\delta x$, and then to optimize this model to obtain a potential optimization step.   For example, approximating $f(x+\delta x) \approx f(x)+f'(x)[\delta x]$ (an affine model, colloquially called ``linear'') leads to gradient descent and related algorithms.  A better approximation for $f(x + \delta x)$ will often lead to faster-converging algorithms, and so a natural idea is to exploit the \emph{second} derivative $f''$ to make a quadratic model, as above, and accelerate optimization.

For unconstrained optimization, minimizing $f(x)$ corresponds to finding a root of the derivative $f' = 0$ (i.e., $\nabla f = 0$), and a \emph{quadratic} approximation for $f$ yields a first-order (affine) approximation $f'(x + \delta x) \approx f'(x) + f''(x)[\delta x]$ for the derivative $f'$. In $\mathbb{R}^n$, this is $\delta(\nabla f) \approx H \delta x$.  So, minimizing a quadratic model is effectively a \emph{Newton step} $\delta x \approx -H^{-1} \nabla f$ to find a root of $\nabla f$ via first-order approximation.   Thus, optimization via quadratic approximations is often viewed as a form of Newton algorithm.  As discussed below, it is also common to employ \emph{approximate} Hessians in optimization, resulting in ``quasi-Newton'' algorithms.

More complicated versions of this idea arise in optimization with constraints, e.g.~minimizing an objective function $f(x)$ subject to one or more nonlinear inequality constraints $c_k(x) \le 0$.   In such cases, there are a variety of methods that take both first and second derivatives into account, such as ``sequential quadratic programming''\footnote{The term ``programming'' in optimization theory does not refer to software engineering, but is rather an anachronistic term for optimization problems.  For example, ``linear programming'' (LP) refers to optimizing affine objectives and affine constraints, while ``quadratic programming'' (QP) refers to optimizing convex quadratic objectives with affine constraints.} (SQP) algorithms that solve a sequence of ``QP'' approximations involving quadratic objectives with affine constraints (see e.g.~the book \textit{Numerical Optimization} by Nocedal and Wright, 2006).

There are many technical details, beyond the scope of this course, that must be resolved in order to translate such high-level ideas into practical algorithms.
For example, a quadratic model is only valid for small enough $\delta x$, so there must be some mechanism to limit the step size. One possibility is ``backtracking line search'': take a Newton step $x+\delta x$ and, if needed, progressively ``backtrack'' to $x+\delta x/10, x+\delta x/100, \ldots$ until a  sufficiently decreased value of the objective is found.  Another commonplace idea is a ``trust region'': optimize the model with the constraint that $\delta x$ is sufficiently small, e.g.~$\Vert \delta x \Vert \le s$ (a spherical trust region), along with some rules to adaptively enlarge or shrink the trust-region size ($s$) depending on how well the model predicts $\delta f$.  There are many variants of Newton/SQP-like algorithms depending on the choices made for these and other details.

\subsubsection{Computing Hessians}

In general, finding $f''$ or the Hessian is often computationally expensive in higher dimensions. If $f(x): \R^n\to\R$, then the Hessian, $H$, is an $n\times n$ matrix, which can be huge if $n$ is large---even storing $H$ may be prohibitive, much less computing it. When using automatic differentiation (AD), Hessians are often computed by a \emph{combination} of forward and reverse modes (Sec.~\ref{sec:forward-over-reverse}), but AD does not circumvent the fundamental scaling difficulty for large~$n$.

Instead of computing $H$ explicitly, however, one can instead \emph{approximate} the Hessian in various ways; in the context of optimization, approximate Hessians are found in ``quasi-Newton'' methods such as the famous ``BFGS'' algorithm and its variants. One can also derive efficient methods to compute Hessian--vector products $Hv$ without computing $H$ explicitly, e.g.~for use in Newton--Krylov methods.   (Such a product $Hv$ is equivalent to a directional derivative of $f'$, which is efficiently computed by ``forward-over-reverse'' AD as in Sec.~\ref{sec:forward-over-reverse}.)

\subsubsection{Minima, maxima, and saddle points}

Generalizing the rules you may recall from single- and multi-variable calculus, we can use the second derivative to determine whether an extremum is a minimum, maximum, or saddle point. Firstly, an extremum of a scalar function $f$ is a point $x_0$  such that $f'(x_0) = 0$. That is, 
\[
f'(x_0) [\delta x] = 0
\] 
for \textit{any} $\delta x$. Equivalently, 
\[
\nabla f \bigr|_{x_0} = f'(x_0)^T = 0.
\]

Using our quadratic approximation around $x_0$, we then have that 
\[
f(x_0 + \delta x) = f(x_0) + \underbrace{f'(x_0) [\delta x]}_{=0} + \frac{1}{2} f''(x_0) [\delta x, \delta x] + o(\lVert \delta x\rVert^2).
\]
The definition of a local minimum $x_0$ is that $f(x_0 + \delta x ) > f(x_0)$ for any $\delta x \neq 0$ with $\lVert \delta x\rVert$ sufficiently small. To achieve this at a point where $f' = 0$, it is enough to have $f''$ be a positive-definite quadratic form:
\[
f''(x_0) [\delta x, \delta x] >0 \text{ for all } \delta x \ne 0 \iff \textbf{positive-definite } f''(x_0) \, .
\]

For example, for inputs $x\in \R^n$, so that $f''$ is a real-symmetric $n \times n$ Hessian matrix, $f''(x_0) = H(x_0) = H(x_0)^T$, this corresponds to the usual criteria for a positive-definite matrix: 
\[
f''(x_0) [\delta x, \delta x] = \delta x^T H(x_0) \delta x >0 \text{ for all } \delta x \ne 0 \iff H(x_0) \text{ positive-definite } \iff \text{all eigenvalues of } H(x_0) > 0.
\]

In first-year calculus, one often focuses in particular on the 2-dimensional case, where $H$ is a $2\times 2$ matrix.  In the $2\times 2$ case, there is a simple way to check the signs of the two eigenvalues of $H$, in order to check whether an extremum is a minimum or maximum: the eigenvalues are both positive if and only if $\det (H) >0$ and $\tr(H)  >0$, since $\det (H) = \lambda_1 \lambda_2$ and $\tr(H) = \lambda_1 + \lambda_2$. In higher dimensions, however, one needs more complicated techniques to compute eigenvalues and/or check positive-definiteness, e.g.~as discussed in MIT courses 18.06 (Linear Algebra) and/or 18.335 (Introduction to Numerical Methods).  (In practice, one typically checks positive-definiteness by performing a form of Gaussian elimination, called a Cholesky factorization, and checking that the diagonal ``pivot'' elements are $> 0$, rather than by computing eigenvalues which are much more expensive.)

Similarly, a point $x_0$ where $\nabla f = 0$ is a local \emph{maximum}
 if $f''$ is negative-definite, or equivalently if the eigenvalues of the Hessian are all negative. Additionally, $x_0$ is a \emph{saddle} point if $f''$ is indefinite, i.e. the eigenvalues include both positive and negative values. However, cases where some eigenvalues are zero are more complicated to analyze; e.g.~if the eigenvalues are all $\ge 0$ but some are $=0$, then whether the point is a minimum depends upon higher derivatives.

\subsection{Further Reading}

All of this formalism about ``bilinear forms'' and so forth may seem like a foray into abstraction for the sake of abstraction.  Can't we always reduce things to ordinary matrices by choosing a basis (``vectorizing'' our inputs and outputs)?  However, we often don't \emph{want} to do this for the same reason that we often prefer to express first derivatives as linear operators rather than as explicit Jacobian matrices.   Writing linear or bilinear operators as explicit matrices, e.g.~$\vecm(A\,dA + dA\,A) = (I\otimes A + A^T \otimes I)\vecm(dA)$ as in Sec.~\ref{sec:kronecker}, often disguises the underlying structure of the operator and introduces a lot of algebraic complexity for no purpose, as well as being potentially computationally costly (e.g.~exchanging small matrices $A$ for large ones $I \otimes A$). 

As we discussed in this chapter, an important generalization of quadratic operations to arbitrary vector spaces come in the form of \href{https://en.wikipedia.org/wiki/Bilinear_map}{bilinear maps} and \href{https://en.wikipedia.org/wiki/Bilinear_form}{bilinear forms}, and there are many textbooks and other sources discussing these ideas and variations thereof. For example, we saw that the second derivative can be seen as a \href{https://en.wikipedia.org/wiki/Symmetric_bilinear_form}{symmetric bilinear form}. This is closely related to a \href{https://en.wikipedia.org/wiki/Quadratic_form}{quadratic form} $Q[x]$, which what we get by plugging the same vector twice into a symmetric bilinear form $B[x,y]=B[y,x]$, i.e.~$Q[x] = B[x,x]$.  (At first glance, it may seem like $Q$ carries ``less information'' than $B$, but in fact this is not the case.  It is easy to see that one can recover $B$ from $Q$ via $B[x,y] = (Q[x+y] - Q[x-y])/4$, called a ``polarization identity.'') For example, the $f''(x) [\delta x, \delta x]/2$ term that appears in quadratic approximations of $f(x+ \delta x)$ is a quadratic form. The most familiar multivariate version of $f''(x)$ is the \href{https://en.wikipedia.org/wiki/Hessian_matrix}{Hessian matrix} when $x$ is a column vector and $f(x)$ is a scalar, and Khan Academy has an elementary \href{https://www.khanacademy.org/math/multivariable-calculus/applications-of-multivariable-derivatives/quadratic-approximations/a/quadratic-approximation}{introduction to quadratic approximation}.

\href{https://en.wikipedia.org/wiki/Definite_matrix}{Positive-definite} Hessian matrices, or more generally \href{https://en.wikipedia.org/wiki/Definite_quadratic_form}{definite quadratic forms} $f''$, appear at extrema ($f' =0$) of scalar-valued functions $f(x)$ that are local minima. There are a lot \href{http://www.columbia.edu/~md3405/Unconstrained_Optimization.pdf}{more formal treatments} of the same idea, and conversely Khan Academy has the \href{https://www.khanacademy.org/math/multivariable-calculus/applications-of-multivariable-derivatives/optimizing-multivariable-functions/a/second-partial-derivative-test}{simple 2-variable version} where you can check the sign of the $2\times 2$ eignevalues just by looking at the determinant and a single entry (or the trace). There's a nice \href{https://math.stackexchange.com/questions/2285282/relating-condition-number-of-hessian-to-the-rate-of-convergence}{stackexchange discussion} on why an \href{https://nhigham.com/2020/03/19/what-is-a-condition-number/}{ill-conditioned} Hessian tends to make steepest descent converge slowly. Some Toronto \href{https://www.cs.toronto.edu/~rgrosse/courses/csc421_2019/slides/lec07.pdf}{course notes on the topic} may also be useful.

Lastly, see for example these Stanford notes on \href{https://web.stanford.edu/class/ee364b/lectures/seq_notes.pdf}{sequential quadratic optimization} using trust regions (Section 2.2), as well as the 18.335 \href{https://github.com/mitmath/18335/blob/spring21/notes/BFGS.pdf}{notes on BFGS quasi-Newton methods}. The fact that a quadratic optimization problem in a sphere has \href{https://en.wikipedia.org/wiki/Strong_duality}{strong duality}, and hence is efficiently solvable, is discussed in Section 5.2.4 of the \href{https://web.stanford.edu/~boyd/cvxbook/}{\textit{Convex Optimization} book}. There has been a lot of work on \href{https://en.wikipedia.org/wiki/Hessian_automatic_differentiation}{automatic Hessian computation}, but for large-scale problems you may only be able to compute Hessian--vector products efficiently in general, which are equivalent to a directional derivative of the gradient and can be used (for example) for \href{https://en.wikipedia.org/wiki/Newton%E2%80%93Krylov_method}{Newton--Krylov methods}.

The Hessian matrix is also known as the ``curvature matrix"
especially in optimization. If we have a scalar function $f(x)$ of $n$ variables, its ``graph'' is the set of points $(x,f(x))$ in $R^{n+1}$; we call the last dimension the ``vertical" dimension. At a ``critical point'' $x$ (where $\nabla f = 0$), then $v^T H v$  is the ordinary curvature sometimes taught in first-year calculus, of the curve obtained by intersecting the graph with the plane in the direction $v$ from $x$ and the vertical (the ``normal section''). The determinant of $H$, sometimes known as the Hessian determinant, yields the Gaussian curvature.

A closely related idea is the derivative of the unit normal.  For a graph as in the 
preceding paragraph we may assume that $f(x)=x^THx/2$ to second order.  It is easy
to see that at any point $x$ the tangents have the form $(dx, f'(x)[dx])=(dx,x^THdx)$
and the normal is then $(Hx,1)$.  Near $x=0$ this a unit normal to second order, and
its derivative is $(Hdx,0)$.  Projecting onto the horizontal, we see that the Hessian
is the derivative of the unit normal.  This is called the ``shape operator" in differential
geometry.

\todo{ The Hessian matrix is also known as the ``curvature matrix"
especially in the fields of optimization.
If we have a scalar function $f(x)$ of $n$ variables the graph of a function is the set of points $(x,f(x))$ in $R^{n+1}$. 
We call the last dimension the ``vertical" dimension.
At a critical point $x$ (meaning 0 gradient), we have $v^THv$ 
is the ordinary curvature taught in some calculus classes of the curve obtained by intersecting the graph with the plane in the direction $v$ from $x$ and the vertical.  (The ``normal section").
The determinant, sometimes known as the Hessian determinant,
is the Gaussian curvature (up to sign because of the embedding??).
Alternatively, one can take a tangent plane to be the horizontal
(effectively rotating the graph) and then the Hessian has
the same properties as above at any point.
Yet another equivalent way to say this, perhaps a bit more elegantly,
is to consider the magnitude of the derivative of the unit normal
with resepect to tangents.  (The derivative
is a function from tangents to tangents.)
TODO: Confirm that the scalar curvature is the sum of the eigenvalues taken two at a time, and write down the Riemannian and Ricci curvature, and place in another chapter?  I think the Ricci curvature is related
to sums of two by two determinants of the Hessian, in other words, it is related to the second compound matrix of the Hessian.
Also in some other chapter some day, maybe go so far as talking
about ...?}

\pagebreak
\section{Derivatives of Eigenproblems}
\subsection{Differentiating on the Unit Sphere}

\begin{figure}[h]
    \centering
    \includegraphics[width=0.3\textwidth]
    {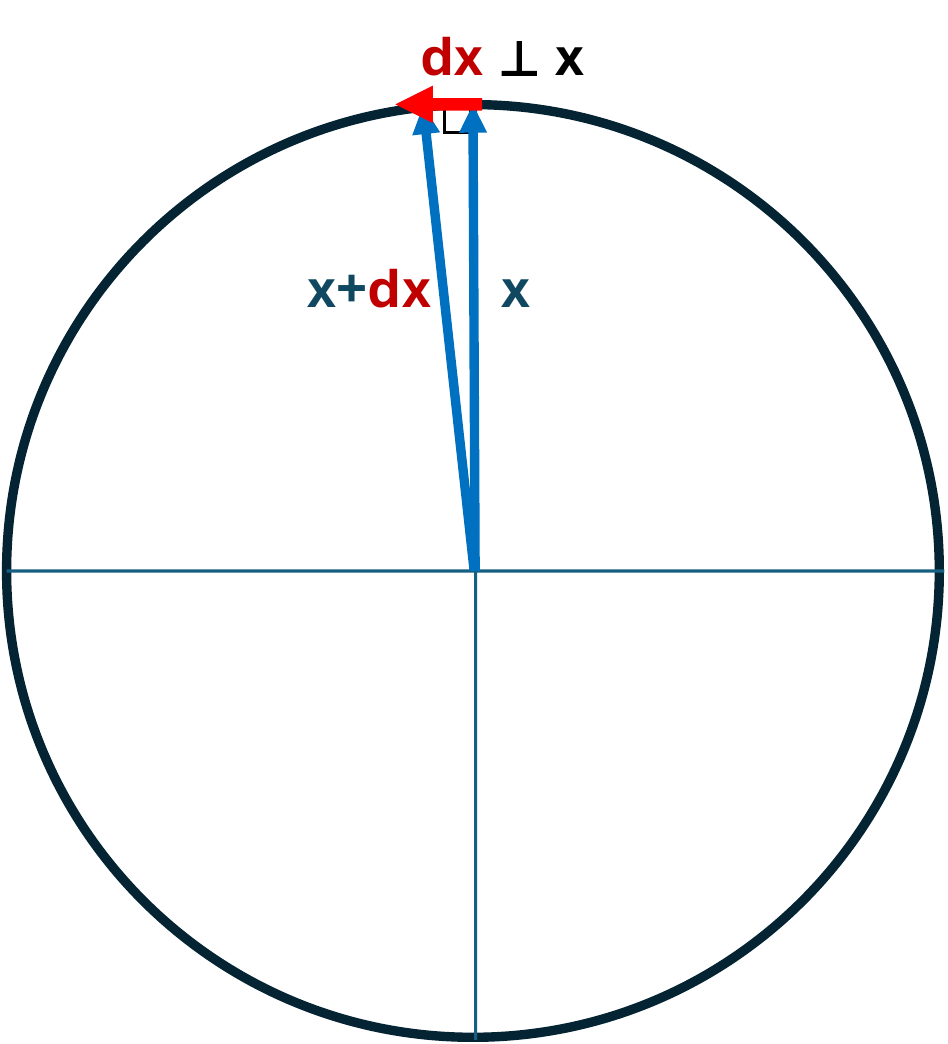}
    \caption{Differentials on a sphere ($x^T x = 1$): the differential $dx$ is constrained to be perpendicular to $x$.}
    \label{fig:tangentcircle}
\end{figure}

Geometrically, we know that velocity vectors (equivalently, tangents) on the sphere are orthogonal to the radii. Out differentials say this algebraically, since given $x\in \mathbb{S}^n$ we have $x^T x = 1$, this implies that 
\[
2x^T \d x = \d(x^T x) = \d (1) = 0.
\]
In other words, at the point $x$ on the sphere (a radius, if you will), $\d x$, the linearization of the constraint of moving along the sphere satisfies $\d x \perp x$. This is our first example where we have seen the infinitesimal perturbation $\d x$ being constrained. See Figure \ref{fig:tangentcircle}.

\subsubsection{Special Case: A Circle}

Let us simply consider the unit circle in the plane where $x = (\cos \theta, \sin \theta)$ for some $\theta \in [0,2\pi)$. Then, 
\[
x^T \d x = (\cos \theta, \sin \theta) \cdot (-\sin \theta, \cos \theta) d\theta = 0.
\]
Here, we can think of $x$ as ``extrinsic'' coordinates, in that it is a vector in $\R^2$. On the other hand, $\theta$ is an ``intrinsic'' coordinate, as every point on the circle is specified by one $\theta$.

\subsubsection{On the Sphere}

You may remember that the rank-1 matrix $x x^T$, for any unit vector $x^T x = 1$, is a \textbf{projection matrix} (meaning that it is equal to its square and it is symmetric) which projects vectors onto their components in the direction of $x$.  Correspondingly, $I - x x^T$ is also a projection matrix, but onto the directions \emph{perpendicular} to $x$: geometrically, the matrix removes components in the $x$ direction. In particular, if $x^T \d x= 0$, then $(I - x x^T) \d x = \d x.$ It follows that if $x^T \d x = 0$ and $A$ is a symmetric matrix, we have 
\begin{align*}
    \d \left(\frac{1}{2} x^T A x\right) &= (A x)^T \d x \\
    &= x^T A (\d x) \\
    &= x^T A ( I - x x^T) \d x  \\
    &= ((I - x x^T) A x )^T \d x.
\end{align*}
In other words, $(I - x x^T)A x$ is the gradient of $\frac{1}{2} x^T A x$ \textit{on the sphere.}

So what did we just do? To obtain the gradient on the sphere, we needed (i) a linearization of the function that is correct on tangents, and (ii) a direction that \textit{is} tangent (i.e.~satisfies the linearized constraint). Using this, we obtain the gradient of a general scalar function on the sphere:

\begin{theorem}
    Given $f: \mathbb{S}^n \to \R$, we have 
    \[
    \d f = g(x)^T \d x = (( I - x x^T) g(x))^T \d x.
    \]
\end{theorem}

The proof of this is precisely the same as we did before for $f(x) = \frac{1}{2} x^T A x$.

\subsection{Differentiating on Orthogonal Matrices}

Let $Q$ be an orthogonal matrix. Then, computationally (as is done in the Julia notebook), one can see that $Q^T \d Q$ is an anti-symmetric matrix (sometimes called skew-symmetric). 
\begin{definition}
A matrix $M$ is anti-symmetric if $M = - M^T$. Note that all anti-symmetric matrices thus have zeroes on their diagonals.
\end{definition}

In fact, we can prove that $Q^T \d Q$ is anti-symmetric.
\begin{theorem}
    Given $Q$ is an orthogonal matrix, we have that $Q^T \d Q$ is anti-symmetric.
\end{theorem}

\begin{proof}
    The constraint of being orthogonal implies that $Q^T Q = I$. Differentiating this equation, we obtain 
    \[
    Q^T \d Q + \d Q^T \, Q = 0 \implies Q^T \d Q = - (Q^T \d Q)^T.
    \]
    This is precisely the definition of being anti-symmetric.
\end{proof}

Before we move on, we may ask what the dimension of the ``surface'' of orthogonal matrices is in $\R^{n^2}$.

When $n = 2$, all orthogonal matrices are rotations and reflections, and rotations have the form 
\[
Q = \begin{pmatrix}
    \cos \theta & \sin \theta \\ 
    - \sin \theta & \cos \theta
\end{pmatrix}.
\]
Hence, when $n=2$ we have one parameter.

When $n = 3$, airplane pilots know about ``roll, pitch, and yaw'', which are the three parameters for the orthogonal matrices when $n =3.$ In general, in $\R^{n^2}$, the orthogonal group has dimension $n(n-1)/2$.

There are a few ways to see this.
\begin{itemize}
    \item Firstly, orthogonality $Q^T Q = I$ imposes $n(n+1)/2$ constraints, leaving $n(n-1)/2$ free parameters.
    \item When we do $QR$ decomposition, the $R$ ``eats'' up $n(n+1)/2$ of the parameters, again leaving $n(n-1)/2$ for $Q$.
    \item Lastly, If we think about the symmetric eigenvalue problem where $S = Q \Lambda Q^T$, $S$ has $n(n+1) /2$ parameters and $\Lambda$ has $n$, so $Q$ has $n(n-1)/2$.
\end{itemize}

\subsubsection{Differentiating the Symmetric Eigendecomposition}

Let $S$ be a symmetric matrix, $\Lambda$ be diagonal containing eigenvalues of $S$, and $Q$ be orthogonal with column vectors as eigenvectors of $S$ such that $S = Q \Lambda Q^T$.  [For simplicity, let's assume that the eigenvalues are ``simple'' (multiplicity~1); repeated eigenvalues turn out to greatly complicate the analysis of perturbations because of the ambiguity in their eigenvector basis.] Then, we have 
\[
\d S = \d Q \, \Lambda Q^T + Q \, \d \Lambda \, Q^T + Q \Lambda \d Q^T,
\]
which may be written as 
\[
Q^T \d S \, Q = Q^T \d Q \Lambda - \Lambda Q^T \d Q + \d \Lambda.
\]

As an exercise, one may check that the left and right hand sides of the above are both symmetric. This may be easier if one looks at the diagonal entries on their own, as there $(Q^T \d S \, Q)_{ii} = q_i^T \d S \, q_i$. Since $q_i$ is the $i$th eigenvector, this implies $q_i^T \d S\, q_i = \d \lambda_i.$   (In physics, this is sometimes called the ``Hellman--Feynman'' theorem, or non-degenerate first-order eigenvalue-perturbation theory.)

Sometimes we think of a curve of matrices $S(t)$ depending on a parameter such as time. If we ask for $\frac{\d \lambda_i }{\d t}$, this implies it is thus equal to $q_i^T \frac{ \d S(t)}{\d t} q_i.$ So how can we get the gradient $\nabla \lambda_i$ for one of the eigenvalues? Well, firstly, note that 
\[
\tr (q_i q_i^T)^T \d S) = \d \lambda_i \implies \nabla \lambda_i = q_i q_i^T.
\]

What about the eigenvectors? Those come from off diagonal elements, where for $i \neq j,$
\[
(Q^T \d S \, Q)_{ij} = \left(Q^T \frac{\d Q}{
\d t}\right)_{ij} (\lambda_j - \lambda_i).
\]
Therefore, we can form the elements of $Q^T \frac{\d Q}{ \d t}$, and left multiply by $Q$ to obtain $\frac{\d Q}{\d t}$ (as $Q$ is orthogonal).

It is interesting to get the second derivative of eigenvalues when moving along a line in symmetric matrix space. For simplicity, suppose $\Lambda$ is diagonal and $S(t) = \Lambda + t E.$ Therefore, differentiating
\[
\frac{\d \Lambda}{\d t} = \diag\left(Q^T \frac{\d S(t)}{\d t} Q\right),
\]
we get 
\[
\frac{\d^2 \Lambda}{\d t^2} = \diag \left(Q^T \frac{\d^2 S(t)}{\d t^2} Q\right) + 2 \diag \left(Q^T \frac{\d S(t)}{\d t} \frac{\d Q}{\d t}\right).
\]
Evaluating this at $Q = I$ and recognizing the first term is zero as we are on a line, we have that 
\[
\frac{\d^2 \Lambda}{\d t^2} = 2 \diag \left(E \cdot  \frac{\d Q}{\d t}\right),
\]
or 
\[
\frac{\d^2 \Lambda}{\d t^2} = 2 \sum_{k \neq i} E_{ik}^2/(\lambda_i - \lambda_k).
\]
Using this, we can write out the eigenvalues as a Taylor series: 
\[
\lambda_i (\epsilon) = \lambda_i + \epsilon E_{ii} + \epsilon^2 \sum_{k \neq i} E_{ik}^2/(\lambda_i - \lambda_k) + \dots.
\]
(In physics, this is known as second-order eigenvalue perturbation theory.)
\pagebreak

\section{Where We Go From Here}
There are many topics that we did not have time to cover, even in 16 hours of lectures. If you came into this class thinking that taking derivatives is easy and you already learned everything there is to know about it in first-year calculus, hopefully we've convinced you that it is an enormously rich subject that is impossible to exhaust in a single course. Some of the things it might have been nice to include are:
\begin{itemize}
    \item When automatic differentiation (AD) hits something it cannot handle, you may have to write a custom Jacobian--vector product (a ``Jvp,'' ``frule,'' or ``pushforward'') in forward-mode, and/or a custon row vector--Jacobian product (a ``vJp,'' ``rrule,'' ``pullback,'' or ``Jacobian$^T$-vector product'') in reverse-mode. In Julia with Zygote AD, this is done using \href{https://github.com/JuliaDiff/ChainRulesCore.jl}{the ChainRules packages}. In Python with JAX, this is done with \href{https://jax.readthedocs.io/en/latest/_autosummary/jax.custom_jvp.html}{jax.custon\_jvp} and/or \href{https://jax.readthedocs.io/en/latest/_autosummary/jax.custom_vjp.html}{jax.custon\_vjp} respectively. In principle, this is straightforward, but the APIs can take some getting used to because the of the generality that they support.
    \item For functions $f(z)$ with complex arguments $z$ (i.e. complex vector spaces), you cannot take ``ordinary'' complex derivatives whenever the function involves the conjugate $\overline{z}$, for example, $|z|, \Re(z),$ and $\Im(z)$. This \textit{must} occur if $f(z)$ is purely real-valued and not constant, as in optimization problems involving complex-number calculations. One option is to write $z = x+iy$ and treat $f(z)$ as a two-argument function $f(x,y)$ with real derivatives, but this can be awkward if your problem is ``naturally'' expressed in terms of complex variables (for instance, the \href{https://en.wikipedia.org/wiki/Frequency_domain}{Fourier frequency domain}). A common alternative is the ``CR calculus'' (or ``Wirtinger calculus''), in which you write 
    \[
    \d f = \left(\frac{\partial f}{\partial z}\right) \d z + \left(\frac{\partial f}{\partial \overline{z}}\right)\d \overline{z},
    \]
    as if $z$ and $\overline{z}$ were independent variables. This can be extended to gradients, Jacobians, steepest-descent, and Newton iterations, for example. A nice review of this concept can be found in these \href{https://arxiv.org/abs/0906.4835}{UCSD course notes} by K. Kreuz Delgado.
    \item Many, many more derivative results for matrix functions and factorizations can be found in the literature, some of them quite tricky to derive. For example, a number of references are listed in this \href{https://github.com/JuliaDiff/ChainRules.jl/issues/117}{GitHub issue for the ChainRules package}.
    \item Another important generalization of differential calculus is to derivatives on curved manifolds and differential geometry, leading to the \href{https://en.wikipedia.org/wiki/Exterior_derivative}{exterior derivative}.
    \item When differentiating eigenvalues $\lambda$ of matrices $A(x)$, a complication arises at eigenvalue crossings (where multiplicity $k>1$). Here, the eigenvalues and eigenvectors usually cease to be differentiable. More generally, this problem arises for any \href{https://en.wikipedia.org/wiki/Implicit_function}{implicit function} with a repeated root. In this case, one option is use an expanded definition of sensitivity analysis called a \textbf{generalized gradient} (a $k\times k$ matrix-valued linear operator $G(x) [\d x]$ whose \textit{eigenvalues} are the perturbations $\d \lambda$. See for example \href{https://doi.org/10.1006/jfan.1995.1117}{Cox (1995)}, \href{https://doi.org/10.1007/BF01742705}{Seyranian \textit{et al.} (1994)}, and \href{https://doi.org/10.1016/j.laa.2022.04.019}{Stechlinski (2022)}. Physicistss call a related idea ``degenerate perturbation theory.'' A recent formulation of similar ideas is called the \textbf{lexicographic directional derivative}. See for example \href{https://doi.org/10.1007/s10107-005-0633-0}{Nesterov (2005)} and \href{https://doi.org/10.1080/10556788.2017.1374385}{Barton \textit{et al.} (2017)}.

    Sometimes, optimization problems involving eigenvalues can be reformulated to avoid this difficulty by using \href{https://en.wikipedia.org/wiki/Semidefinite_programming}{SDP constraints}. See for example \href{http://doi.org/10.1364/OE.22.022632}{Men \textit{et al.} (2014)}.

    For a \href{https://en.wikipedia.org/wiki/Defective_matrix}{defective matrix} the situation is worse: even the generalized derivatives blow up because $\d \lambda$ can be proportional to (e.g.) the square root of the perturbation $\lVert \d A\rVert$ (for an eigenvalue with algebraic multiplicity $=2$ and geometric multiplicity $=1$).

    \item Famous generalizations of differentation are the ``\href{https://en.wikipedia.org/wiki/Distributional_derivative}{distributional}'' and ``\href{https://en.wikipedia.org/wiki/Weak_derivative}{weak}'' derivatives. For example, to obtain \href{https://en.wikipedia.org/wiki/Dirac_delta_function}{Dirac delta ``functions''} by differentiating discontinuities. This requires changing not only the definition of ``derivative,'' but also changing the definition of \textit{function}, as reviewed at an elementary level in these \href{https://math.mit.edu/~stevenj/18.303/delta-notes.pdf}{MIT course notes}.
\end{itemize}

\end{document}